\documentclass{aims}
\usepackage{amsmath}
\usepackage{mathabx}
  \usepackage{paralist}
\usepackage{graphicx}  \usepackage{epstopdf}
 \usepackage[colorlinks=true]{hyperref}
\hypersetup{urlcolor=blue, citecolor=red}

\def\currentvolume{X}
 \def\currentissue{X}
  \def\currentyear{200X}
   \def\currentmonth{XX}
    \def\ppages{X--XX}
     \def\DOI{10.3934/xx.xx.xx.xx}

\newtheorem{theorem}{Theorem}[section]

\newtheorem{lemma}[theorem]{Lemma}
\newtheorem{proposition}{Proposition}

\theoremstyle{definition}

\newtheorem{remark}{Remark}

\def\debproof{ {\it Proof.} }
\def\finproof{\hfill {\small $\Box$} \\}

\newcommand{\CC}{\mathbb{C}}
\newcommand{\RR}{\mathbb{R}}

\newcommand{\EE}{\mathbb{E}}

\newcommand{\chemin}{.}

\title{Mean Field Model for Collective Motion Bistability}
\author[Josselin Garnier, George Papanicolaou, and Tzu-Wei Yang]{}

	

\email{josselin.garnier@polytechnique.edu}
\email{papanicolaou@stanford.edu}
\email{yangx953@umn.edu}

\thanks{$^*$ Corresponding author}

\begin{document}
\maketitle

\centerline{\scshape Josselin Garnier$^*$}
\medskip
{\footnotesize
 \centerline{Centre de Math\'ematiques Appliqu\'ees, Ecole Polytechnique}
 \centerline{91128 Palaiseau Cedex, France}
} 

\medskip

\centerline{\scshape George Papanicolaou}
\medskip
{\footnotesize
 \centerline{Department of Mathematics, Stanford University}
 \centerline{Stanford, CA 94305, USA}
}

\medskip

\centerline{\scshape Tzu-Wei Yang}
\medskip
{\footnotesize
 \centerline{School of Mathematics, University of Minnesota}
 \centerline{Minneapolis, MN 55455, USA }
}

\bigskip

 \centerline{(Communicated by the associate editor name)}

\begin{abstract}
	We consider the Czir\'ok model for collective motion of locusts along a one-dimensional torus.  
	In the model, each agent's velocity locally interacts with other agents' velocities in the system, 
	and there is also exogenous randomness to each agent's velocity. 
	The interaction tends to create the alignment of collective	motion.  
	By analyzing the associated nonlinear Fokker-Planck equation, 
	we obtain the condition for the existence of stationary order states and the conditions for their linear stability. 
	These conditions depend on the noise level, which should be strong enough, and on the interaction between the agent's velocities,
	which should be neither too small, nor too strong.
	We carry out the fluctuation analysis of the interacting system and describe the large deviation principle to calculate 
	the transition probability from one order state to the other. Numerical simulations confirm our analytical findings. 
\end{abstract}

\section{Introduction}

Collective motion has become an emerging topic in biology and statistical physics and has a large potential for the applications in other fields, for example, 
opinion dynamics in social science \cite{Motsch2014} and autonomous robotic networks in engineering \cite{Savkin2003}.  Collective motion describes the systemic 
behavior of a large group of animals, insects or bacteria.  In a collective motion model, the system has a large number of individuals, expressed as their locations 
and velocities, $(x_i, u_i)$, $i=1,\ldots,N$.  The movement of each individual is affected by the movements of other individuals in the system, mostly attractive 
interactions.  Most biologically realistic and mathematically interesting models suppose \textit{local} interactions, that is, any individual is affected 
only by its neighbors within a limited range.  In addition, a collective motion model generally adds independent, internal or external randomness to the individuals'  movements.

Two approaches are widely used to model collective motion: 
self-propelled particle (SPP) models and coarse-grained (CG) models.  In an SPP model (also known as an individual-based model), 
 the system is viewed as a large number of coupled stochastic processes
describing the evolutions of the individuals in the system.  
In an SPP model we can easily build an elaborated interaction between an individual and its neighbors to match actual biological observations.
However, an SPP model consists of numerous nonlinear equations and is difficult to analyze.  To obtain the analytical tractability, further simplifications
are usually required \cite{Yates2009,Escudero2010}; otherwise the system behavior has to be verified by numerical 
simulations.  On the other hand, a CG model (also known as a macroscopic or mean-field model) assumes infinitely many individuals in the system and considers the evolution of the distribution of the individuals.  
Most CG models in the literature lead to hydrodynamic equations so that the system behavior can be analyzed.  
CG models have some limitations.  Although some CG models \cite{Bertin2006,Degond2007,Ha2008,Bertin2009,Mishra2010} are indeed constructed from microscopic dynamics, it is generally difficult to relate the models to actual individuals' properties \cite{Ariel2015}. 
Another limitation of  CG models is that noise present at the microscopic level
can be averaged out by coarse graining.
As a result most hydrodynamic equations only describe the deterministic evolution and cannot capture and analyze some 
important stochastic events like transitions from one order state to another order state observed 
for instance in \cite{Buhl2006}.
It is possible to incorporate noise in the CG model in a  phenomenological way \cite{Bertin2013}
but it is important to model noise at the macroscopic level in a consistent way and to clarify what type of noise
in the CG model is generated by noise in the SSP model.

In the recent years, there is a huge amount of literature concerning the modeling of collective motion so we mention only a few papers guiding our own work.  
Vicsek \textit{et al.} \cite{Vicsek1995} propose the scalar noise model that is currently the most well-known SPP model.  The Czir\'ok model \cite{Czirok1999} is a one-dimensional model on a torus. It is a simplified version of
the Vicsek model and is used to describe the collective motion of locusts for instance
(see \cite{Vicsek2012,Ariel2015} and references therein).
Buhl \textit{et al.} \cite{Buhl2006} 
experimentally and numerically study the switching of the alignments of locusts and find that the probability of switching is exponentially decaying with the number of locusts.  
In \cite{Yates2009,Escudero2010} the authors assume global and uniform interactions in the Czir\'ok model and they 
derive the nonlinear Fokker-Planck equations for the average velocity and consider overall stochastic behavior. 
Some interesting features, like the existences of leaders or predators are considered in 
\cite{Couzin2005} and \cite{Reynolds2009}, respectively.  
Other various individual-based models related to the Czir\'ok model are discussed in 
\cite{OLoan1999,Chate2008,Romanczuk2012,Ariel2014}, and \cite{Toner1995,Edelstein-Keshet1998,Toner1998,Bertin2006,Bertin2009,Mishra2010,Ihle2011,Topaz2012,Solon2013,Caussin2014,Degond2015,Solon2015,Solon2015b} construct 
and analyze CG models of collective motion. 
These papers show that the linear stability analysis of the homogeneous solutions is a simple but powerful technique,
and that the formation of ``bands" (regions with high density of particles moving in the same direction) can be observed in numerical simulations
and related to traveling solutions of CG equations.

One of the most important goals of the collective motion modeling is to mathematically reproduce and analyze the \textit{ordering phenomenon} where the individuals 
show a coherence in their movements, such as swarms of locusts, flocks of birds, or schools of fish. Such an ordering phenomenon is called an order state of the 
system and there are generally multiple order states that represent different directions and rotations.  The counterpart of the order states is the disorder state, 
where the system does not show such a pattern.  From the perspective of statistical physics, in this paper we address several interesting questions that are also 
considered in many of the literature:
\begin{enumerate}
	\item Existence and number of order states.
	
	\item Quantities that control the phase transition between the disorder state and the order states.
	
	\item Transition probability from one order state to another order state.
\end{enumerate} 

Our contributions in this paper are the following. 
We interpret the Czir\'ok model as an interacting particle system in the McKean-Vlasov framework \cite{Dawson1983,Dawson1987}.
In this sense, we characterize the state of the system by the empirical probability measure 
of the locations and velocities of all the individuals, and as the number 
of individuals goes to infinity, the random, empirical measure converges to a deterministic,
probability measure whose density is the solution of a nonlinear Fokker-Planck equation.  
Therefore, when the number of individuals is large, we can consider the Fokker-Planck equation instead of all the individuals.  
In this way, we are also able to provide a clear connection between the SPP model and the CG model.  
This procedure was proved to be very efficient to study the stationary states and their stability properties 
in the case of opinion dynamics models \cite{Garnier2016}. It turns out that it is still efficient to address collective motion models,
although the framework is different: In the opinion dynamics  model addressed in \cite{Garnier2016},
the  Fokker-Planck equation is non-degenerate while it is degenerate in the collective motion model that we addresse here,
as the external noise affects only the velocities and not the locations. As a result noise plays a very different role
which strongly affects the conditions for the existence and stability of the stationary states.
For instance, the noise has to be large enough to ensure the existence of the non-trivial stationary states 
in the opinion dynamics  model. In dramatic contrast, we will see that the existence of the order states in the collective motion model
does not depend on the noise level and that their stability depends in a nontrivial way of the noise level.
Noise has to be large enough to ensure the linear stability of the order states, but large noise also triggers
frequent transitions from one order state to another one.

By analyzing the Fokker-Planck equation obtained in the Czir\'ok model, 
we find that the system may have one disorder state and two order states 
(clockwise and counterclockwise rotations at constant velocity).  
We find the necessary and sufficient condition for the existence of the order states. 
This existence condition quantitatively describes the 
phase transition between the disorder state and the order states.  
We then perform the linear stability analysis to examine the stabilities of the order and disorder 
states, and we find that the order states are stable when the external noises are sufficiently strong. 
The noise threshold value depends in a nontrivial way of the interaction between individuals,
the smallest threshold value is reached when the interaction is neither too weak, nor too strong.
Such a phenomenon that the randomness of movement improves 
the alignment of collective motion coincides with the observation in \cite{Yates2009}. 
When noise is small, the order states are unstable and the complex modulational instability 
generates moving clusters whose velocities and sizes can be identified.
Finally, when the number of individuals is large but finite, the empirical 
measure is still stochastic and has a small probability to transit from one order state to the other one.  
This small probability can be described by the large 
deviation principle. It increases with the noise level. 
It is exponentially decaying as a function of the number of individuals; 
this fact is also observed in \cite{Buhl2006,Yates2009}.  
We confirm our analysis by extensive numerical simulations.

The paper is organized as follows.  The interacting particle model and its mean-field limit is presented in section~\ref{sec:model}.  The equilibria of the 
nonlinear Fokker-Planck equation are analyzed in section~\ref{sec:equilibrium}.  In section~\ref{sec:linear stability}, the linearized stability analysis of the 
nonlinear Fokker-Planck equation is performed.  We provide the fluctuation analysis and the large deviation principle in sections~\ref{sec:fluctuation} and 
\ref{sec:large deviation}, respectively.  In section~\ref{sec:numerics}, we verify our analytical findings by numerical simulations.  
We briefly consider the nonsymmetric
Czir\'ok model with normalized influence functions in section~\ref{sec:nonsymmetry}.  
We end with a brief summary and conclusions in section~\ref{sec:conclusion}.

\section{Model and Mean Field Limit}
\label{sec:model}

We consider the Czir\'ok model, a system of $N$ agents moving along the torus $[0,L]$.  The position $x_i$ and velocity $u_i$ of particle $i$ satisfy the following 
stochastic differential equations: for $i=1,\ldots,N$,
\begin{equation}
	\label{eq:model}
	\begin{aligned}
		dx_i &= u_i dt , \\
		du_i & = \left[ G( \left<u\right>_i)-u_i \right]dt+\sigma dW_i(t),
	\end{aligned}
\end{equation}
where $\{W_i(t)\}_{i=1}^N$ are independent Brownian motions and $\left<u\right>_i$ is a weighted average of the velocities $\{u_j\}_{j=1}^N$:
\begin{equation}
	\label{def:aveu}
	\left<u\right>_i = \frac{1}{N} \sum_{j=1}^N u_j\phi( \|x_j-x_i \|),
\end{equation}
with the weights depending on the distance $\|\cdot\|$ on the torus between the position $x_i$ and the positions of the other agents (i.e. $\|x_j-x_i\|=\min(|x_j-x_i|,L-|x_j-x_i|)$).
Here $\phi(x)$ is a nonnegative influence function normalized so that
\begin{equation}
	\label{eq:phi0}
	\frac{1}{L} \int_0^L \phi(\|x\|) dx = 1 ,
\end{equation}
and $G(u)$ is an odd and smooth function. Note that there is no loss of generality in assuming (\ref{eq:phi0})
since we can always rescale $G$ to ensure that this condition is satisfied.
As we will see in this paper the interesting configuration is when $u\mapsto u-G(u)$ derives from an even double-well potential,
with two symmetric global minima and one saddle point at zero.
For instance, we can think at
\begin{equation}
	\label{eq:exG1}
	G(u) = 2 \tanh(u) 
\end{equation}
and then $u-G(u)$ derives from the double-well 
potential $u^2 /2- 2 \log (\cosh(u)) $, or 
\begin{equation}
	\label{eq:exG2}
	G(u)=2u-u^3
\end{equation} 
and then $u-G(u)$ derives from the double-well potential $ u^4/4 - u^2/2$. 

We define the empirical probability measure:
\begin{equation}
	\label{eq:empirical measure}
	\mu_N(t,dx,du) = \frac{1}{N} \sum_{i=1}^N \delta_{(x_i(t),u_i(t))}(dx,du).
\end{equation}
Let us assume that $\mu_N(0,dx,du)$ converges to a deterministic measure  $\bar{\rho}(x,u)dxdu$ as $N\to\infty$.
This happens in particular when the positions and velocities of the agents at the inital time $t=0$ are
independent and identically distributed with the distribution with density $\bar{\rho}(x,u)$.
Then, as $N\to\infty$, $\mu_N(t,dx,du) $ weakly converges to the 
deterministic measure $\rho(t,x,u)dxdu$ whose density is the solution of the nonlinear Fokker-Planck equation 
\begin{equation}
	\label{eq:nonlinear Fokker-Planck}
	\frac{\partial \rho}{\partial t} = 
	- u\frac{\partial \rho}{\partial x}
	- \frac{\partial}{\partial u}\left\{\left[G\left(\iint u'\phi(\|x'\|) \rho(t,x-x',u') du'dx' \right) - u \right] \rho\right\}
	+ \frac{1}{2} \sigma^2 \frac{\partial^2 \rho}{\partial u^2} ,
\end{equation}
starting from $\rho(t=0,x,u)=\bar{\rho}(x,u)$.  Note that the Fokker-Planck equation (\ref{eq:nonlinear Fokker-Planck}) is degenerate because there is no diffusion 
in $x$.  The convergence proof is standard and can be found in \cite{Dawson1983,Gartner1988}.

\section{Stationary Analysis for Order and Disorder States}
\label{sec:equilibrium}

We look for a stationary equilibrium, that is to say, a probability density function stationary in time:
\begin{equation}
	\label{eq:stationary}
	- u\frac{\partial \rho}{\partial x}
	- \frac{\partial}{\partial u}\left\{\left[G\left(\iint u'\phi(\|x'\|) \rho(x-x',u') du'dx' \right) - u \right] \rho\right\}
	+ \frac{1}{2} \sigma^2 \frac{\partial^2 \rho}{\partial u^2}
	= 0 ,
\end{equation}
with periodic boundary conditions in $x$.

\begin{proposition}
	\label{prop:equilibrium}
	The probability-density-valued solutions $\rho(x,u)$ to (\ref{eq:stationary}) have the following form:
	\begin{equation}
		\label{eq:equilibrium}
		\rho_\xi(x,u) = \frac{1}{L} F_\xi(u),  \quad F_\xi(u)= \frac{1}{\sqrt{\pi\sigma^2}} \exp
		\Big(- \frac{(u-\xi)^2}{\sigma^2} \Big).
	\end{equation}
	They are uniform in space, Gaussian in velocity, and their mean velocity $\xi$ satisfies the compatibility condition:
	\begin{equation}
		\label{eq:compatibility condition}
		\xi = G( \xi) .
	\end{equation}
\end{proposition}
There are, therefore, as many stationary equilibria as there are solutions to the compatibility equation (\ref{eq:compatibility condition}).

\noindent
\debproof
It is straightforward to check that (\ref{eq:equilibrium}) is a solution of (\ref{eq:stationary}). 
Reciprocally, let $\rho$ be a solution of (\ref{eq:stationary}).
	After integrating (\ref{eq:stationary}) from $u=-\infty$ to $u=+\infty$, we obtain
	\begin{equation*}
		-\frac{\partial}{\partial x}\int_{-\infty}^\infty u \rho(x,u) du = 0,
	\end{equation*}
	which shows that $\int_{-\infty}^\infty u \rho(x,u) du = \xi/L$ is a constant  independent of $x$.
	Then Eq.~(\ref{eq:stationary}) can be written as a linear equation
	\begin{equation}
		\label{eq:simplified stationary}
		- u\frac{\partial \rho}{\partial x}
		- \frac{\partial}{\partial u}\left\{[G(\xi) - u]\rho\right\}
		+ \frac{1}{2} \sigma^2 \frac{\partial^2 \rho}{\partial u^2}
		= 0 .
	\end{equation}
	Let $F(u) = \int_0^L \rho(x,u)dx$. Then $F(u)$ satisfies
	\begin{equation}
		\label{eq:u stationary}
		- \frac{\partial}{\partial u}\{[G(\xi) - u] F(u)\}
		+ \frac{1}{2} \sigma^2 \frac{\partial^2 F(u)}{\partial u^2}
		= 0  ,
	\end{equation}
        which shows that $F(u)$ is a Gaussian density function:
	\begin{equation}
		F(u) = \frac{1}{\sqrt{\pi\sigma^2}}\exp\Big(-\frac{1}{\sigma^2}(u-G(\xi))^2\Big).
	\end{equation}
	The condition $\xi = L\int u \rho(x,u) du$ requires that 
	\begin{equation}
		\xi = \frac{1}{L}\int_0^L \xi dx = \int_0^L \int u \rho(x,u) du dx = \int u F(u) du = G(\xi).
	\end{equation}
	To show that $\rho(x,u)$ is uniform in $x$, we first note that $\rho(x,u)$ is periodic in $x$. We expand it as:
	\begin{equation}
		\rho(x,u) = \sum_{k=-\infty}^\infty \rho_k(u) e^{- i 2 \pi k x /L},\quad \rho_k(u) =\frac{1}{L}  \int_0^L \rho(x,u) e^{i 2 \pi k x /L} dx.
	\end{equation}
	From (\ref{eq:simplified stationary}) and (\ref{eq:compatibility condition}), $\rho_k$ satisfies
	\begin{equation}
		\label{eq: eqn for rho_k}
		\frac{i2\pi k}{L} u\rho_k
		- \frac{\partial}{\partial u} \left[ (\xi-u) \rho_k\right]
		+ \frac{1}{2} \sigma^2 \frac{\partial^2 \rho_k}{\partial u^2}
		= 0 .
	\end{equation}
	We take a Fourier transform in $u$:
	\begin{equation}
		\hat{\rho}_k(\eta) =  \int_{-\infty}^\infty \rho_k(u)e^{-i \eta u} du .
	\end{equation}
	From (\ref{eq: eqn for rho_k}), $\hat{\rho}_k$ satisfies the ordinary differential equation
	\begin{equation}
		-\Big( \frac{2\pi k}{L}  + \eta \Big) \frac{\partial \hat{\rho}_k}{\partial \eta}  
		= \Big( i \xi \eta +\frac{1}{2} \sigma^2 \eta^2\Big) \hat{\rho}_k .
	\end{equation}
This equation can be solved. Let $k<0$. For any $\eta \in (-\infty,-2 \pi k/ L)$ we have
\begin{equation*}
	|\hat{\rho}_k(\eta)| = |\hat{\rho}_k(0)| \left|1+\frac{L\eta}{2\pi k}\right|^{- \frac{2\pi^2 k^2\sigma^2 }{L^2}}  \exp \Big( -\frac{\sigma^2 \eta^2}{4}
	+\frac{\pi k \sigma^2\eta}{L} \Big) ,
\end{equation*}
which goes to $+\infty$ as $\eta \nearrow -2 \pi k / L$ if $\hat{\rho}_k(0) \neq 0$.
However $|\hat{\rho}_k(\eta)| \leq \int |\rho_k(u)|du \leq (1/L)\int \int_0^L \rho(x,u) dx du =1/L$
is uniformly bounded, which imposes $\hat{\rho}_k(0)=0$ and therefore $\hat{\rho}_k(\eta)=0$ for any $\eta \in (-\infty,-2 \pi k/L)$.
Similarly, for any $\eta \in (-2 \pi k/ L,\infty)$ we have
\begin{equation*}
	|\hat{\rho}_k(\eta)| = |\hat{\rho}_k(-4\pi k /L)| \Big|1+\frac{L\eta}{2\pi k}\Big|^{- \frac{2\pi^2 k^2\sigma^2}{L^2}}  \exp \Big( \big(-\frac{\sigma^2 \eta}{4}
	+\frac{2\pi k \sigma^2}{L} \big) \big(\eta+ \frac{4\pi k }{L} \big) \Big) ,
\end{equation*}
which goes to $+\infty$ as $\eta \searrow -2 \pi k / L$ if $\hat{\rho}_k(-4\pi k /L) \neq 0$.
The boundedness of $\hat{\rho}_k$ imposes $\hat{\rho}_k(-4\pi k/L)=0$ and therefore $\hat{\rho}_k(\eta)=0$ for any $\eta \in (-2 \pi k/L,\infty)$.
Moreover, $\hat{\rho}_k$ is continuous by the Riemann-Lebesgue lemma and therefore $ \hat{\rho}_k(\eta) = 0$ for any $\eta\in \RR$ and $k<0$.
Since $\hat{\rho}_{-k}(-\eta) =\overline{\hat{\rho}_k(\eta)}$,
we also find that  $ \hat{\rho}_k(\eta) = 0$ for any $\eta\in \RR$ and $k>0$.
Therefore $\rho_k(u)=0$ for any $u \in \RR$ and $k \neq 0$ and the stationary solution is equal to
$\rho(x,u) =  \rho_0(u) $.
\finproof

When $G$ is such that $u-G(u)$ derives from an even double-well potential, such as the two examples (\ref{eq:exG1}) and (\ref{eq:exG2}),
there are three $\xi$ satisfying the compatibility condition (\ref{eq:compatibility condition}): $0$ and $\pm\xi_e$, with $\xi_e>0$
the unique positive solution of $G(\xi)=\xi$.  
In the experiment, 
the equilibrium $\rho_0(x,u)$ is the disorder state and $\rho_{\pm\xi_e}(x,u)$ are the two order states of the locusts marching on the torus, the clockwise and 
counterclockwise rotations.  The existence of the nonzero solution $\pm\xi_e$ and hence the stationary states $\rho_{\pm\xi_e}(x,u)$ is independent of 
the value of the noise strength $\sigma$.  This observation is unusual and in fact inconsistent with the conclusion in many models in statistical physics and opinion 
dynamics  \cite{Dawson1983,Garnier2016}.
We will see, however, that the noise strength plays a crucial role in the stability of the stationary states.

\section{Linear Stability Analysis}
\label{sec:linear stability}
In this section, we assume that  the compatibility condition (\ref{eq:compatibility condition}) has solutions 
and we use the linear stability analysis to study the stability of the stationary states.
Let $\xi$ be such that $G(\xi)=\xi$ and consider 
\begin{equation}
	\rho(t,x,u) = \rho_\xi(x,u)+\rho^{(1)}(t,x,u) = \frac{1}{L}F_\xi(u) + \rho^{(1)}(t,x,u),
\end{equation}
for small perturbation $\rho^{(1)}$.  
By linearizing the nonlinear Fokker-Planck equation (\ref{eq:nonlinear Fokker-Planck}) we find that $\rho^{(1)}$ satisfies
\begin{align}
	\label{eq:linearized Fokker-Planck}
	\frac{\partial \rho^{(1)}}{\partial t} &= 
	-u\frac{\partial \rho^{(1)}}{\partial x} 
	-\frac{\partial}{\partial u}\left[(\xi-u) \rho^{(1)}\right]\\
	\notag
	&\quad - \frac{1}{L}G'(\xi)\left[\iint u' \phi(\|x'\|) \rho^{(1)}(t, x-x', u') du'dx' \right] F'_\xi(u)
	+ \frac{1}{2} \sigma^2 \frac{\partial^2 \rho^{(1)}}{\partial u^2}.
\end{align}
Our goal is to clarify under which circumstances the linearized system is stable.

\subsection{Expressions of the linearized modes}
The perturbation  $\rho^{(1)}$ is periodic in $x$. We expand it as:
\begin{equation}
\label{eq:expand:rho1k}
	\rho^{(1)}(t,x,u) = \sum_{k=-\infty}^\infty \rho^{(1)}_k(t,u) e^{- i 2 \pi k x /L},\quad 
	\rho^{(1)}_k(t,u) = \frac{1}{L}  \int_0^L \rho^{(1)}(t,x,u) e^{i 2 \pi k x /L} dx  .
\end{equation}
The mode $\rho^{(1)}_k$ satisfies
\begin{align}
	\label{eq:modelink}
	\frac{\partial\rho_k^{(1)}}{\partial t} &= \frac{i 2\pi k}{L} u \rho_k^{(1)} - \frac{\partial}{\partial u}\left[(\xi-u) \rho_k^{(1)}\right]\\
	&\quad - G'(\xi)\phi_k\left[\int u'\rho_k^{(1)}(t,u') du'\right] F'_\xi(u) + \frac{1}{2} \sigma^2 \frac{\partial^2 \rho_k^{(1)}}{\partial u^2},\notag
\end{align}
where the discrete Fourier coefficients 
\begin{equation}
\phi_k = \frac{1}{L}  \int_0^L \phi(\|x\|) e^{i 2 \pi k x /L} dx
\end{equation}
are real-valued (because $\|L-x\|=\|x\|$) and bounded by one (because $|\phi_k| \leq \int_0^L \phi(\|x\|)dx/L=1$).
We note that the equations are uncoupled in $k$ and the system is linearly stable if all modes are stable. 
We then take a Fourier transform in $u$:
\begin{equation}
	\hat{\rho}_k^{(1)}(t,\eta) =  \int_{-\infty}^\infty \rho_k^{(1)}(t,u)e^{-i \eta u} du.
\end{equation}
From (\ref{eq:modelink}), $\hat{\rho}_k^{(1)}$ satisfies the following first-order partial differential equation:
\begin{equation}
	\frac{\partial}{\partial t} \hat{\rho}_k^{(1)} + \Big(\frac{2\pi k}{L} + \eta\Big) \frac{\partial}{\partial \eta} \hat{\rho}_k^{(1)}
	= \Big(-i\xi\eta - \frac{1}{2}\sigma^2 \eta^2\Big) \hat{\rho}_k^{(1)}
	+ G'(\xi) \phi_k \Big[\frac{\partial}{\partial \eta} \hat{\rho}_k^{(1)}(t,0)\Big] \eta \hat{F}_\xi,
	\label{eq:hatrhok}
\end{equation}
starting from an initial condition that is assumed to satisfy 
\begin{equation*}
	\| \hat{\rho}_k^{(1)}(0,\cdot) \|_{L^1}+
	\| \hat{\rho}_k^{(1)}(0,\cdot) \|_{L^\infty}+
	\|\partial_\eta \hat{\rho}_k^{(1)}(0,\cdot) \|_{L^\infty}  
	<\infty .
\end{equation*}
We say that the $k$th-order mode is stable if 
\begin{equation*}
	\sup_{t \geq 0} \big\{
	\| \hat{\rho}_k^{(1)}(t,\cdot) \|_{L^1} +
	\| \hat{\rho}_k^{(1)}(t,\cdot) \|_{L^\infty} +
	\|\partial_\eta \hat{\rho}_k^{(1)}(t,\cdot) \|_{L^\infty}  
	\big\}
	<\infty .
\end{equation*}
By the method of characteristics, we can obtain the following implicit expression for $\hat{\rho}_k^{(1)}(t,\eta)$ that 
is useful for the stability analysis:
\begin{lemma}
	\label{lma:characteristics}
	Let $w_k(t) = {\partial_\eta} \hat{\rho}_k^{(1)}(t,0)$. $\hat{\rho}_k^{(1)}(t,\eta)$ has the following expression:
	\begin{align}
		\label{eq:characteristics}
		\hat{\rho}_k^{(1)}(t,\eta) &= e^{-g_k(t,\eta)} \hat{\rho}_k^{(1)} \Big(0 , e^{-t} \eta - D_k(1-e^{-t})\Big)\\
		&\quad + G'(\xi)\phi_k \int_0^t e^{-g_k(t-s,\eta)} w_k(s) H_\xi \Big(e^{-(t-s)} \eta - D_k(1-e^{-(t-s)})\Big) ds,\notag
	\end{align}
	where $D_k = 2\pi k/L$, $H_\xi(\eta)=\eta \hat{F}_\xi(\eta)$ and 
	\begin{align}
		\label{def:gk}
		g_k(t,\eta) &= \frac{\sigma^2 }{4} (\eta+D_k)^2 (1-e^{-2t}) + (i\xi -\sigma^2 D_k) (\eta + D_k)(1-e^{-t})\\
		&\quad + \Big( -i\xi D_k + \frac{\sigma^2}{2}D_k^2\Big)t.\notag
	\end{align}
\end{lemma}

\noindent
\debproof
	The function $\hat{\rho}_k^{(1)}(t,\eta) $ satisfies (\ref{eq:hatrhok}).
	We consider the mapping $\eta(t)=\eta_0 e^t - {2\pi k}/{L}$ so that $\eta'(t) = \eta(t) + {2\pi k}/{L}$. By the method of characteristics, 
	\begin{equation*}
		\frac{d}{d t} \hat{\rho}_k^{(1)}(t, \eta(t)) 
		= \Big( -i\xi\eta - \frac{1}{2}\sigma^2 \eta^2 \Big)\hat{\rho}_k^{(1)} (t, \eta)
		+ G'(\xi) \phi_k \Big[\frac{\partial}{\partial \eta} \hat{\rho}_k^{(1)}(t, 0)\Big] \eta \hat{F}_\xi(\eta).
	\end{equation*}
	Then by using the method of integrating factors, we have
	\begin{align*}
		\hat{\rho}_k^{(1)}(t, \eta(t)) &= e^{-\int_0^t g(s)ds} \hat{\rho}_k^{(1)}(0, \eta_0 - 2\pi k/L)\\
		&\quad + G'(\xi) \phi_k \int_0^t e^{-\int_s^t g(w)dw} \Big[\frac{\partial}{\partial \eta} \hat{\rho}_k^{(1)}(s, 0)\Big] 
		\eta(s) \hat{F}_\xi(\eta(s))ds, 
	\end{align*}
	where $g(t) = i\xi\eta(t) + \frac{1}{2}\sigma^2 \eta^2(t)$.
	We note that $\int_s^t g(w) dw = g_k(t-s,\eta(t))$, where $g_k$ is defined by (\ref{def:gk})
	and $D_k = 2\pi k/L$.  By letting $H_\xi(\eta)=\eta \hat{F}_\xi(\eta)$ and $w_k(s)={\partial_\eta} \hat{\rho}_k^{(1)}(s, 0)$, 
	$\hat{\rho}_k^{(1)}(t,\eta)$ can be written as (\ref{eq:characteristics}).
\finproof

\subsection{Stability of the $0$th-order mode}
In this section we find a necessary and sufficient condition for the stability of the mode $k=0$.
For $k=0$,  $\rho_0^{(1)}$ satisfies
\begin{equation}
	\label{eq:modelin0}
	\frac{\partial\rho_0^{(1)}}{\partial t}
	= -\frac{\partial}{\partial u}\left[(\xi-u) \rho_0^{(1)}\right]
	- G'(\xi)\left[\int u'\rho_0^{(1)}(t,u')du'\right] F'_\xi(u) + \frac{1}{2} \sigma^2 \frac{\partial^2 \rho_0^{(1)}}{\partial u^2}.
\end{equation}
We use the method of moments.
We let $m_{0,0}^{(1)}$, $m_{0,1}^{(1)}$ denote the zeroth and first moments of $\rho_0^{(1)}$, respectively:
\begin{equation}
	m_{0,0}^{(1)}(t) = \int \rho^{(1)}_0(t,u) du,\quad m_{0,1}^{(1)}(t) = \int u\rho^{(1)}_0(t,u) du.
\end{equation}
From (\ref{eq:modelin0}), we get
\begin{equation*}
	\frac{\partial}{\partial t}m_{0,0}^{(1)} = 0, \quad
	\frac{\partial}{\partial t}m_{0,1}^{(1)} = \xi m_{0,0}^{(1)} - m_{0,1}^{(1)} + G'(\xi)  m_{0,1}^{(1)}.
\end{equation*}
The first moment $m_{0,1}^{(1)}$ is stable if and only if $G'(\xi) < 1$.
More exactly, if $G'(\xi)  = 1$ then the first moment grows linearly in time.  If $G'(\xi)  <1$  then the first moment is bounded.
Once $m_{0,1}^{(1)}$ has been shown to bounded, then the solution $\rho_0^{(1)}$ of (\ref{eq:modelin0}) can be seen as the solution of a  linear parabolic equation which is smooth and bounded as explained in the following proposition.

\begin{proposition}
	\label{prop:stablity of mode 0}
	The $0$th-order mode is stable if and only if $G'(\xi)  < 1$.
\end{proposition}

\noindent
\debproof
	We know that $m_{0,1}^{(1)}(t)$ is uniformly bounded in $t$ if and only if $G'(\xi)  < 1$. 
	From (\ref{eq:characteristics}) in Lemma \ref{lma:characteristics} and noting that 
	$w_0(t) = {\partial_\eta} \hat{\rho}_0^{(1)}(t,0) = -im_{0,1}^{(1)}(t)$, 
	we have the following expression for $\hat{\rho}_0^{(1)}(t,\eta)$:
	\begin{align}
	\label{eq:expressrho0p}
		\hat{\rho}_0^{(1)}(t,\eta) &= \hat{\rho}_0^{(1)}(0,\eta e^{-t}) \exp\Big( -i\xi\eta(1-e^{-t}) - \frac{1}{4}\sigma^2 \eta^2(1-e^{-2t}) \Big) \\
		&\quad  -i\int_0^t F_0(\eta e^{-s} ) m_{0,1}^{(1)}( t-s) \exp\Big( -i\xi\eta(1-e^{-s}) - \frac{1}{4}\sigma^2 \eta^2(1-e^{-2s})\Big) ds ,\notag
	\end{align}
	with $F_0(\eta) =  G'(\xi)  \eta \hat{F}_\xi(\eta)$.  We find that, for any $t\geq 1$, 
	\begin{align*}
		|\hat{\rho}_0^{(1)}(t,\eta)|&\leq \| \hat{\rho}_0^{(1)}( 0,\cdot) \|_{L^\infty} e^{ - \frac{1}{4}\sigma^2 \eta^2 (1-e^{-2}) }  \\
		&\quad + \|m_{0,1}^{(1)}\|_{L^\infty} \int_0^t |F_0(\eta e^{-s} ) |\exp \Big(  -\frac{1}{4}\sigma^2 \eta^2(1-e^{-2s})\Big) ds \\
		&\leq \| \hat{\rho}_0^{(1)}( 0,\cdot) \|_{L^\infty} e^{ - \frac{1}{4}\sigma^2 \eta^2 (1-e^{-2})} + \| m_{0,1}^{(1)}\|_{L^\infty} 
		\int_0^1 |F_0(\eta e^{-s} ) | ds  \\
		&\quad + \| m_{0,1}^{(1)}\|_{L^\infty} |G'(\xi)| \| \hat{F}_\xi\|_{L^\infty} \int_{L^1}^t |\eta| e^{-s} 
		\exp\Big(  -\frac{1}{4}\sigma^2 \eta^2(1-e^{-2})\Big) ds\\
		&\leq \| \hat{\rho}_0^{(1)}( 0,\cdot) \|_{L^\infty} e^{ - \frac{1}{4}\sigma^2 \eta^2(1-e^{-2}) } 
		+ \| m_{0,1}^{(1)}\|_{L^\infty}     \int_0^1 |F_0(\eta e^{-s} ) | ds  \\
		&\quad + \| m_{0,1}^{(1)}\|_{L^\infty} |G'(\xi)| \| \hat{F}_\xi\|_{L^\infty} |\eta| e^{ - \frac{1}{4}\sigma^2 \eta^2(1-e^{-2}) } .
	\end{align*}
	Therefore, there exists a constant $C$ that depends only on $G$, $\xi$ and $\sigma$ such that
	\begin{equation*}
		\| \hat{\rho}_0^{(1)}(t,\cdot)\|_{L^1}+ \| \hat{\rho}_0^{(1)}(t,\cdot)\|_{L^\infty}\leq C .
	\end{equation*} 
	We can proceed in the same way after differentiating with respect to $\eta$ 
	the expression (\ref{eq:expressrho0p}) of $\hat{\rho}_0^{(1)}$ to show that $\|\partial_\eta \hat{\rho}_0^{(1)}(t,\cdot)\|_{L^\infty}  $ is bounded uniformly in $t$.
\finproof

	When $u-G(u)$ derives from an even double-well potential,
	such as the examples (\ref{eq:exG1}) and (\ref{eq:exG2}), then, as soon as 
	there exist nonzero solutions $\pm\xi_e$ to the compatibility equation $\xi = G(\xi)$, we have 
	$G'(\pm\xi_e)  < 1$ because $\partial_u(u-G(u))_{u=\pm \xi_e}>0$.  Therefore, from Proposition \ref{prop:equilibrium} and 
	Proposition \ref{prop:stablity of mode 0}, the existence of the order states $\rho_{\pm\xi_e}$ in (\ref{eq:nonlinear Fokker-Planck}) is equivalent to
	the stability of the $0$-th mode $\rho_0^{(1)}$ in (\ref{eq:modelin0}).  
	In addition, the condition that the equation $\xi = G(\xi)$ has the nonzero solutions $\pm\xi_e$ implies $G'(0)  > 1$ 
	 because $\partial_u(u-G(u))_{u=0}<0$, and therefore the disorder state $\rho_0$ is an unstable equilibrium to  the nonlinear Fokker-Planck equation (\ref{eq:nonlinear Fokker-Planck}) when the order states $\rho_{\pm\xi_e}$ exist.

\subsection{Sufficient condition for the stability of the $k$th-order modes}
In this section we find a sufficient condition for the stability of the $k$th-order modes for $k\neq0$.
We show that all the nonzero modes are stable for sufficiently large $\sigma$.
From (\ref{eq:characteristics}) in Lemma \ref{lma:characteristics}, $w_k(t) = {\partial_\eta} \hat{\rho}_k^{(1)}(t,0)$ satisfies
\begin{equation}
	\label{eq:w_k}
	w_k(t) = \psi_k(t) + \int_0^t R_k(t-s) w_k(s) ds,
\end{equation}
where
\begin{align}
	\beta_k (t) &= D_k(1-e^{-t}),\\
	\label{def:psik}
	\psi_k(t) &= \left( -\partial_\eta g_k(t,0) \hat{\rho}_k^{(1)} (0, -\beta_k) +\partial_\eta  \hat{\rho}_k^{(1)} (0 , -\beta_k) e^{-t} \right) e^{-g_k(t,0)},\\
	\label{eq:R_k}
	R_k(t) &= G'(\xi) \phi_k  \left[-\partial_\eta g_k(t,0) H_\xi(-\beta_k(t)) + H_\xi'(-\beta_k(t))e^{-t}\right]e^{-g_k(t,0)},
\end{align}
The strategy is to show that  $\psi_k \in L^\infty(0,\infty)$  and 
that $\int_0^\infty |R_k(t)|dt<1$, and thus 
$\|w_k\|_{L^\infty} \leq (1-\|R_k\|_{L^1})^{-1}\|\psi_k\|_{L^\infty}$ 
by (\ref{eq:w_k}).
Once it is known that $w_k$ is uniformly bounded, it is not difficult to show that the $k$th mode is stable
by the inspection of the formula (\ref{eq:characteristics}).

\begin{lemma}
	\label{prop:stablity of w_k}
\begin{enumerate}
\item
$\psi_k$ defined by (\ref{def:psik})  belongs to $L^\infty[0,\infty)$.
\item
Let $R_k$ be defined by (\ref{eq:R_k}). 
If $|G'(\xi)  \phi_k|<1$ for all nonzero $k$, then $\int_0^\infty |R_k(t)|dt<1$ for all nonzero $k$
if $\sigma$ is large enough.
\end{enumerate}
\end{lemma}
\debproof
See Appendix \ref{app:A}.
\finproof

\begin{proposition}
	\label{prop:stablity of mode k}
	If $|G'(\xi) \phi_k|<1$  for all nonzero $k$ and $\sigma$ is large enough, then all the $k$th-order modes for $k \neq 0$ are stable.\\
\end{proposition}
\debproof
	First we write $\hat{\rho}_k^{(1)}(t,\eta)$ as (\ref{eq:characteristics}). By Lemma \ref{prop:stablity of w_k}, $w_k\in L^\infty[0,\infty)$ and thus we can show that 
	$\| \hat{\rho}_k^{(1)}(t,\cdot)\|_{L^1}+\| \hat{\rho}_k^{(1)}(t,\cdot)\|_{L^\infty}+\| \partial_\eta \hat{\rho}_k^{(1)}(t,\cdot)\|_{L^\infty}$ is 
	uniformly bounded by the estimate similar to 
	the proof of Proposition \ref{prop:stablity of mode 0}.
\finproof

\subsection{Necessary and sufficient condition for the instability of the $k$th-order modes}
\label{subsec:instaspa}%
The linear stability analysis that we have presented reveals that the stability
of the $k$-th mode ($k\neq 0$) is determined by the behavior of the function $w_k$ solution of (\ref{eq:w_k}).
The $k$-th mode is stable if and only if the function $w_k$ is bounded.
We have shown that the condition $\int_0^\infty |R_k(t)|dt <1$
is a sufficient condition for the stability of the $k$-th mode, which is not, however, necessary.
It is possible to give a necessary and sufficient condition for the  instability of the $k$-th mode using the Fourier-Laplace transform of (\ref{eq:w_k}).
This condition is that 
\begin{equation}
C_k = \Big\{ \gamma \in \CC \mbox{ s.t. } {\rm Re}(\gamma)>0 \mbox{ and } \int_0^\infty R_k(t)e^{-\gamma t} dt =1 \Big\}
\end{equation}
is non-empty.
If this happens, then 
the Fourier-Laplace transform of $w_k$:
$$
{\mathcal W}_k(\gamma) = \int_0^\infty w_k(t) e^{-\gamma t} dt
$$
blows up when the complex frequency $\gamma$ goes to an element in $C_k$, 
so we can conclude that $|w_k(t)|$ grows exponentially with the growth rate
$\gamma_r(k)=\sup \{ {\rm Re}(\gamma), \, \gamma \in C_k\}$.
If there exists a unique $\gamma(k) \in C_k$ with ${\rm Re}(\gamma(k))=\gamma_r(k)$,
then we can predict that the $k$-th mode grows like $\exp(  \gamma(k) t)$.
The real part $\gamma_r(k)$ is the growth rate of the $k$-th mode,
the imaginary part $\gamma_i(k)$ describes the dynamical behavior of the unstable mode
as follows.

The mode with the largest growth rate is the one that drives the instability.
If we denote 
$k_{\rm max} = {\rm argmax}_k ( \gamma_r(k) )$, 
then the instability has the form (up to a multiplicative constant):
\begin{equation}
\label{eq:predictvel}
\exp [ - i ( 2\pi k_{\rm max}  x /L - \gamma_i(k_{\rm max} ) t) ] \exp [ \gamma_r(k_{\rm max} ) t] .
\end{equation}
For a fixed time, 
its spatial profile has the form of a periodic modulation with the spatial period $L/ k_{\rm max}$.
This modulation moves in time at the velocity $\gamma_i(k_{\rm max} ) L/(2\pi k_{\rm max})$
and grows  exponentially in time with the rate $ \gamma_r(k_{\rm max} )$.
This shows that the initial uniform distribution of the positions is not stable and that spatial clustering appears.

\begin{figure}
	\centering
	{\bf a)}
	\includegraphics[width=0.38\linewidth]{\chemin/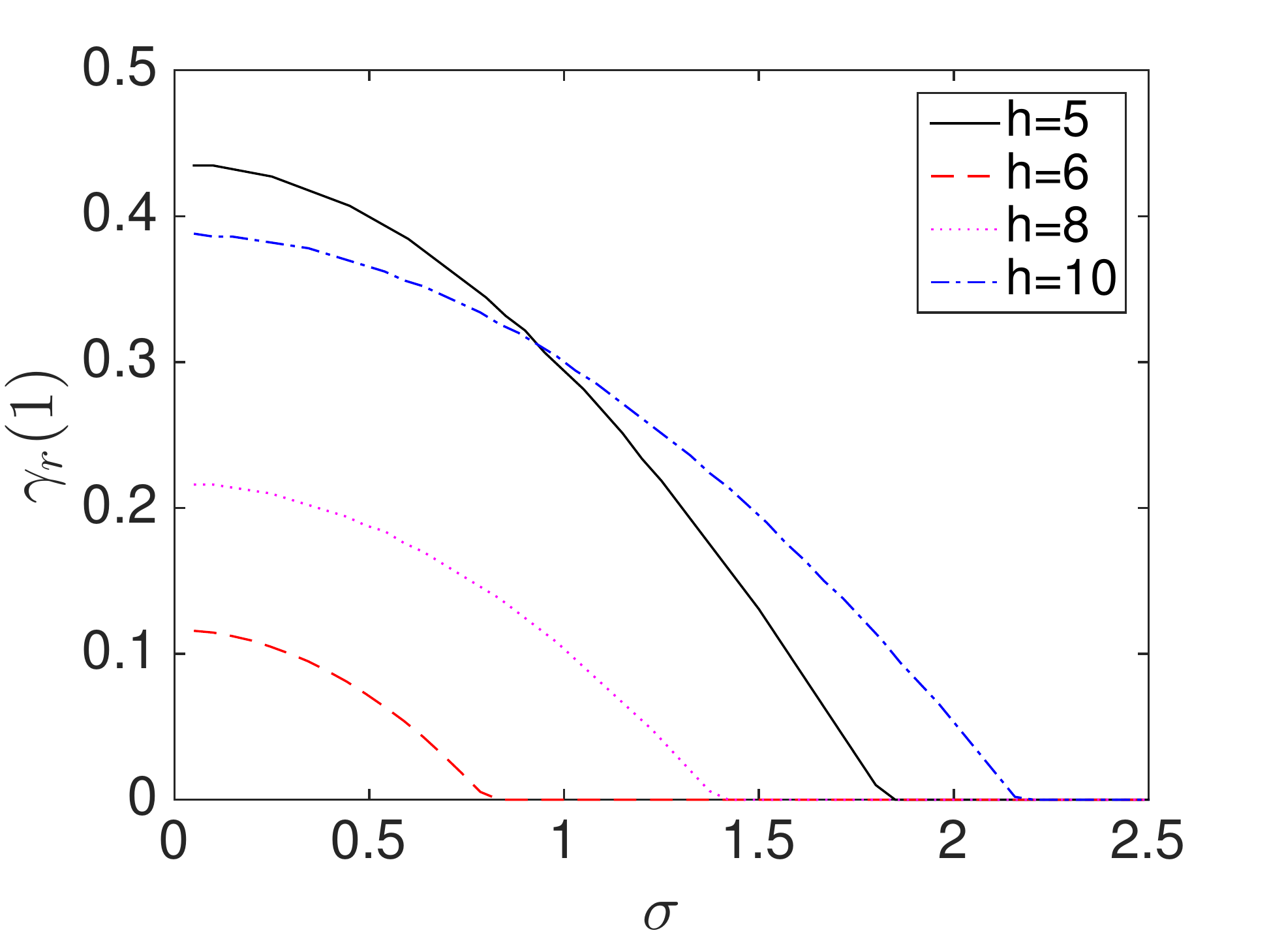}
	{\bf b)}
	\includegraphics[width=0.38\linewidth]{\chemin/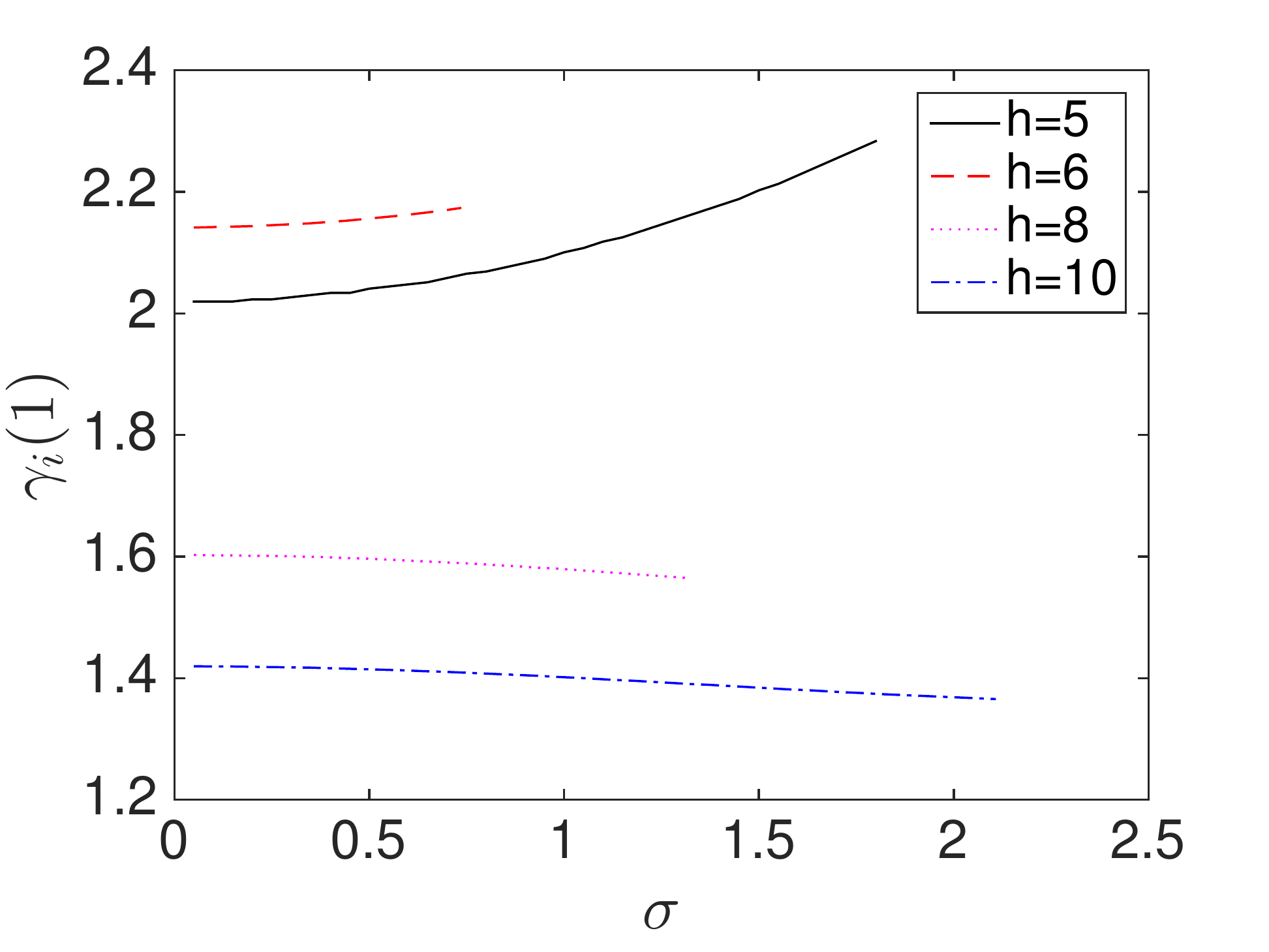}

	{\bf c)}
	\includegraphics[width=0.38\linewidth]{\chemin/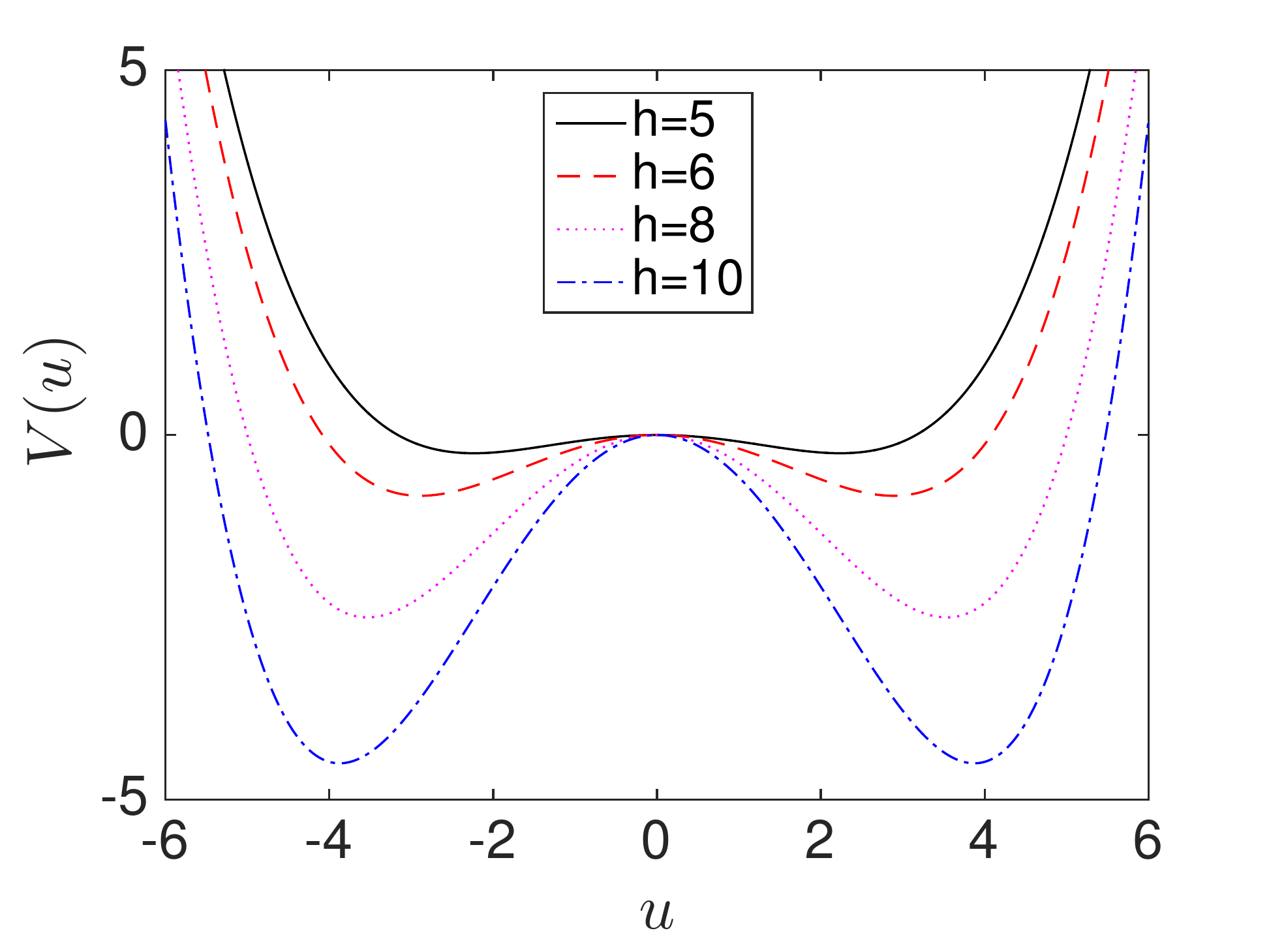}

	\caption{The real (a) and imaginary (b) parts of the complex growth rate of the first mode (the most unstable mode) 
	as a function of $\sigma$ for $h=5$ (solid black), $h=6$ (dashed red), $h=8$ (magenta dotted), and $h=10$ (blue dot-dashed). 
	Here $\phi$ is given by (\ref{def:phi}) and $L=10$. $G$ is given by (\ref{def:G}) and it derives from the double-well potential plotted in picture c.
	We address the linear stability of the order state $\rho_{\xi_e}$. One can see that the threshold value for the noise level $\sigma$  to ensure linear stability is $1.8$ for $h=5$, $0.85$ for $h=6$, $1.4$ for $h=8$, and $2.2$ for $h=10$. The most stable situation is the one corresponding to $h=6$.}
	\label{fig:growthrate1}
\end{figure}

In Figure \ref{fig:growthrate1} we plot the growth rate of the most unstable mode 
in the situation addressed in the numerical section \ref{sec:numerics}.
The double-well potential $G$ defined by (\ref{def:G})
is parameterized by $h$ which determines the depths of the wells and one can see 
that the linear stability is optimal when $h$ is nor too small, neither too large.
This will be confirmed by the numerical simulations presented below.

\subsection{Conclusion}
If $\xi=G(\xi)$ and $G'(\xi)  < 1$, then the $0$th-order mode of the stationary state $\rho_\xi$ is stable.
If, additionally, $|G'(\xi) \phi_k| < 1$ for all $k \neq 0$,  
then all the $k$th-order modes for $k\neq 0$ are stable for a sufficiently large $\sigma$. 
Under these conditions the stationary state $\rho_\xi$ is stable.
When $G$ is such that $u-G(u)$ derives from an even double-well potential, such as the two examples (\ref{eq:exG1}) and (\ref{eq:exG2}), this means that the order states 
\begin{equation*}
	\rho_{\pm\xi_e}(x,u) = \frac{1}{L} \frac{1}{\sqrt{\pi\sigma^2}} \exp\Big(- \frac{(u\mp\xi_e)^2}{\sigma^2} \Big),\quad 
	\xi_e = G(\xi_e) >0,
\end{equation*}
are stable equilibria to the nonlinear Fokker-Planck equation (\ref{eq:nonlinear Fokker-Planck}),
while the disorder state 
\begin{equation}
	\rho_0(x,u) = \frac{1}{L} \frac{1}{\sqrt{\pi\sigma^2}} \exp\Big(- \frac{u^2}{\sigma^2} \Big)
\end{equation}
is an unstable equilibrium.
This shows that the noise strength $\sigma$ improves the stability of the order states $\rho_{\pm\xi_e}(x,u)$.
Such a phenomenon that the randomness of the movements establishes the alignment of collective motion is also found in \cite{Yates2009} both numerically and experimentally.
	
Note, however, an interesting phenomenon: 
if $G'(\xi)<1$, then the $0$th-order mode is stable, 
but for some $k \neq 0$  the $k$th-order mode may be unstable.
This instability is very different from the instability that results from violating the condition $G'(\xi)<1$,
because it manifests itself in an instability in the distribution of the positions of the particles which cannot stay uniform. 
This situation happens when $u-G(u)$ derives from a double-well potential and the wells are very deep, as we will see later in the numerical section. The prediction that the order states become unstable when the depths of the wells are too deep
seems counterintuitive and indeed this never happens when the interaction is uniform $\phi \equiv 1$, but it happens when 
the interaction is local.

\section{Fluctuations Analysis}
\label{sec:fluctuation}

In this section, we consider the fluctuation analysis for the interacting system of agents around the stationary state. 
 Recall that $\mu_N$ is the empirical probability measure of the locations and 
velocities (\ref{eq:empirical measure}).
Here we consider a solution $\xi$ to the compatibility equation $\xi = G(\xi)$.
We assume that, at time $0$, the initial positions and velocities of the agents are independent and identically distributed 
with the distribution with density $\rho_\xi(x,u)$.
Therefore $\mu_N(0,dx,du)$ converges to the stationary state 
$\rho_{\xi}(x,u)dxdu$ in (\ref{eq:equilibrium}) as $N\to\infty$. 
Moreover, $\mu_N(t,dx,du)$ converges to $\rho_{\xi}(x,u)dxdu$ for all $t$ as $N\to\infty$ by (\ref{eq:nonlinear Fokker-Planck}).

We define, for any test function  $f(x,u)$ and any measure $X(dx,du)$ (in the space ${\mathcal S}'$ of tempered Schwartz distributions
on $[0,L]\times \RR$):
$$
\langle X,f\rangle = \int_0^L\int_{-\infty}^\infty f(x,u) X(dx,du)  .
$$ 
We have
$$
\langle \mu_N(t), f \rangle = \frac{1}{N}\sum_{i=1}^N f(x_i(t),u_i(t)) .
$$ 
We define the fluctuation process around $\rho_{\xi}(x,u)dxdu$: 
\begin{equation}
X_N(t,dx,du)=\sqrt{N}[\mu_N(t,dx,du)-\rho_{\xi}(x,u)dxdu],
\end{equation}
and study its dynamics as $N\to\infty$.

First, by central limit theorem, at time $0$, 
$$
\sqrt{N} \Big( \langle \mu_N(0), f \rangle - \langle \rho_\xi dxdu ,f \rangle
\Big) = \frac{1}{\sqrt{N}}\sum_{i=1}^N\big(  f(x_i(0),u_i(0)) -  \langle \rho_\xi dxdu ,f \rangle \big)
$$
converges to a Gaussian random variable with mean zero and variance $\langle \rho_\xi dxdu ,f^2 \rangle -\langle \rho_\xi dxdu ,f \rangle ^2$.
This shows that $X_N(t=0,dx,du)$ converges weakly as $N \to \infty$ to a Gaussian measure $X(t=0,dx,du)$
with mean zero and covariance 
\begin{align}
\nonumber
\EE \big[ \langle X(t=0), f_1 \rangle \langle X(t=0), f_2 \rangle  \big]
=
 \int_0^L\int_{-\infty}^\infty f_1(x,u) f_2(x,u) \rho_\xi (x,u) dx du \\
-
 \int_0^L\int_{-\infty}^\infty  f_1(x,u)  \rho_\xi (x,u) dx du
 \int_0^L\int_{-\infty}^\infty f_2(x,u) \rho_\xi (x,u) dx du , \label{eq:covXN0}
\end{align}
for any test functions $f_1$ and $f_2$.

Second, by It\^o's formula,
\begin{align*}
	d \langle  \mu_N,f\rangle 
	&= \frac{1}{N}\sum_{i=1}^N df(x_i(t),u_i(t))\\
	&= \langle  \mu_N,  u\frac{\partial f}{\partial x}\rangle dt
	+ \langle \mu_N, \big[G \big(\langle  \mu_N(t,dx',du'), u'\phi(\|x'-x\|)\rangle \big) - u \big]\frac{\partial f}{\partial u} \rangle dt\\
	&\quad + \langle  \mu_N, \frac{1}{2}\sigma^2\frac{\partial^2 f}{\partial u^2}\rangle dt
	+  \frac{1}{N}\sum_{i=1}^N \langle \delta_{(x_i,u_i)}(dx,du)  , \sigma\frac{\partial f}{\partial u} \rangle dW_i(t).
\end{align*}
By integration by parts, we can  write 
\begin{align*}
	d\mu_N 
	&= -u \frac{\partial \mu_N}{\partial x} dt 
	- \frac{\partial}{\partial u}\left\{ [G(\langle  \mu_N(t,dx',du'), u'\phi(\|x'-x\|)\rangle) - u] \mu_N\right\} dt\\
	&\quad + \frac{1}{2}\sigma^2 \frac{\partial^2 \mu_N}{\partial u^2}dt
	- \sigma \frac{\partial}{\partial u} \Big[\frac{1}{N}\sum_{i=1}^N \delta_{(x_i,u_i)}(dx,du)dW_i(t)\Big].
\end{align*}
Therefore $X_N$ satisfies
\begin{align}
\notag
	dX_N &= \sqrt{N} d[\mu_N(t,dx,du)-\rho_{\xi}(x,u)dxdu]\\
\notag
	&= -u \frac{\partial X_N}{\partial x} dt
	- \frac{\partial}{\partial u}\left\{ [G(\langle\mu_N(t,dx',du'), u'\phi(\|x'-x\|)\rangle) - u] \sqrt{N}\mu_N\right\} dt\\
\notag
	&\quad + \frac{\partial}{\partial u}[(G(\xi)- u) \sqrt{N}\rho_{\xi}]dt
	+ \frac{1}{2}\sigma^2 \frac{\partial^2 X_N}{\partial u^2}dt\\
	&\quad - \frac{\sigma}{\sqrt{N}}\sum_{i=1}^N \frac{\partial}{\partial u}\delta_{(x_i,u_i)}(dx,du)dW_i(t) .
\notag
\end{align}
	Note that, because 
\begin{equation*}
	\xi=G(\xi) =  G\big( \langle \rho_{\xi}(x',u') dx'du', u'\phi(\|x'-x\|)\rangle \big) ,
\end{equation*}
we have
\begin{align*}
&
G(\langle\mu_N(t,dx',du'), u'\phi(\|x'-x\|)\rangle)\\
&=
G \big(\langle  \rho_{\xi}(x',u') dx'du'+N^{-1/2} X_N(t,dx',du'), u'\phi(\|x'-x\|)\rangle \big) \\
&=
G\big(\xi+ N^{-1/2} \langle  X_N(t,dx',du'), u'\phi(\|x'-x\|)\rangle \big) ,
\end{align*}
so that we can also write that $X_N$ satisfies
\begin{align}
\notag
	dX_N 	
	&= -u \frac{\partial X_N}{\partial x} dt\\
	\notag 
	&\quad 
	- \frac{\partial}{\partial u}\left\{ \big[G \big(\xi+ N^{-1/2}\langle X_N(t,dx',du'), u'\phi(\|x'-x\|)\rangle \big) - u \big] X_N\right\} dt\\
\notag
	&\quad + \frac{\partial}{\partial u}\left\{\sqrt{N} \big[G(\xi)-G\big(\xi+N^{-1/2} \langle X_N(t,dx',du'), u'\phi(\|x'-x\|)\rangle \big) \big]\rho_{\xi}\right\}dt\\
	&\quad + \frac{1}{2}\sigma^2 \frac{\partial^2 X_N}{\partial u^2}dt
	+ \sigma dW^N(t),
\label{eq:XN1}
\end{align}
where $W^N(t)$ is the measure-valued stochastic process
\begin{equation*}
	W^N(t) = - \int_0^t \frac{1}{\sqrt{N}}\sum_{i=1}^N \frac{\partial}{\partial u}\delta_{(x_i(s),u_i(s))}(dx,du)dW_i(s) .
\end{equation*}
The process $W^N(t)$ in ${\mathcal C}([0,\infty), {\mathcal S}')$ is such that, for any test function $f(x,u)$,
\begin{align*}
	\langle W^N(t),f\rangle
	&= -\int_0^t \int_0^L \int_{-\infty}^\infty f(x,u)\frac{1}{\sqrt{N}}\sum_{i=1}^N \frac{\partial}{\partial u}\delta_{(x_i(s),u_i(s))}(dx,du)dW_i(s)\\
	&= \int_0^t \frac{1}{\sqrt{N}}\sum_{i=1}^N \frac{\partial f}{\partial u}(x_i(s),u_i(s)) dW_i(s).
\end{align*}
It is a zero-mean martingale whose quadratic variation is, for any test functions $f_1$ and $f_2$,
\begin{align*}
	 \left[\langle W^N, f_1\rangle, \langle W^N, f_2 \rangle\right]_t
	& = \int_0^{t} \frac{1}{N} \sum_{i=1}^{N} \frac{\partial f_1}{\partial u}(x_i(s),u_i(s))\frac{\partial f_2}{\partial u}(x_i(s),u_i(s)) ds\\
	& = \int_0^{t} \langle\mu_N (s,dx,du), \frac{\partial f_1}{\partial u}\frac{\partial f_2}{\partial u}\rangle ds ,
\end{align*}
which converges to the deterministic process
\begin{align*}
 \left[\langle  W^N,f_1\rangle, \langle  W^N,f_2 \rangle\right]_t	\stackrel{N \to \infty}{\longrightarrow}
 t \, \langle  \rho_{\xi}(x,u) dxdu,\frac{\partial f_1}{\partial u}\frac{\partial f_2}{\partial u}\rangle .
\end{align*}
In other words, as $N\to\infty$, $W^N$ converges to a Brownian field,
and this implies that $X_N$ converges to a generalized Ornstein-Uhlenbeck process, as explained in 
the following proposition.
\begin{proposition}
\label{prop:fluctuation}
Let $\mu_N(t,dx,du) = \frac{1}{N} \sum_{i=1}^N \delta_{(x_i(t),u_i(t))}(dx,du)$ and $\xi$ be a solution to $\xi = G(\xi)$.
We assume that the initial positions and velocities of the agents are independent and identically distributed 
with the distribution with density $\rho_\xi(x,u)$.\\
1)  As $N\to\infty$, $\mu_N(0,dx,du)$ converges to the stationary state $\rho_{\xi}(x,u)dxdu$  and
	$X_N(0,dx,du)$ converges to the Gaussian measure-valued process $X(0,dx,du)$ with mean zero and covariance
	(\ref{eq:covXN0}).\\
2) As $N\to\infty$, $X_N(t,dx,du)$ converges to the Gaussian
	measure-valued process $X(t,dx,du)$ satisfying
	\begin{align}
	\label{eq:linearized Fokker-Planck brown}
		dX &= -u\frac{\partial X}{\partial x}dt - \frac{\partial}{\partial u}[(\xi- u) X]dt\\
		&\quad - G'(\xi) \frac{\partial\rho_{\xi}}{\partial u} \left[\int_0^L\int_{-\infty}^\infty u'\phi(\|x'-x\|)X(t,dx', du')\right] dt
		+ \frac{1}{2} \sigma^2 \frac{\partial^2 X}{\partial u^2}dt + \sigma dW_{\xi},
		\notag
	\end{align}
	where $W_{\xi}(t,x,u)$ is a Gaussian process independent of  $X(0,dx,du)$ with mean zero and covariance
	\begin{align}
	\nonumber
		&\mathbf{Cov}\left( \int_0^L\int_{-\infty}^{\infty} W_{\xi}(t,x,u) f_1(x,u)dxdu, \int_0^L\int_{-\infty}^{\infty} W_{\xi}(t',x,u) f_2(x,u)dxdu\right)\\
		&\quad = \min(t,t')\int_0^L \int_{-\infty}^{\infty} \frac{\partial f_1}{\partial u}(x,u) \frac{\partial f_2}{\partial u}(x,u) \rho_{\xi}(x,u) dxdu ,
	\end{align}
	for any test functions $f_1$ and $f_2$.
\end{proposition}
The convergences hold in the sense of weak convergence on ${\mathcal C}([0, \infty), {\mathcal S}')$ where 
${\mathcal S}'$  denotes the space of tempered Schwartz distributions on $[0,L]\times \RR$.

\noindent
\debproof
	We only provide the main steps of the derivation of Proposition \ref{prop:fluctuation}, 
	which follows the one of \cite[Theorem 4.1.1]{Dawson1983}. We first note that the solution $X_N$ of (\ref{eq:XN1}) 
	is the solution of a martingale problem. To prove the existence of a weak limit, we need to 
	prove that the set $\{X_N\}_{N=1}^\infty$ solution of (\ref{eq:XN1}) is weakly compact and thus the sequence $\{X_N, N=1,2,3,\ldots\}$ has limits.  To prove the uniqueness of the weak limit, we need to prove that any weak limit $X$ is solution of a well-posed martingale problem. 
	Note that the second and third terms in the right-hand side of (\ref{eq:XN1}) satisfy
\begin{equation*}
 \big[G \big(\xi+N^{-1/2} \langle X_N(t,dx',du'), u'\phi(\|x'-x\|)\rangle \big) - u \big] X_N \stackrel{N\to\infty}{\longrightarrow}
(\xi - u )X
\end{equation*}
and
\begin{eqnarray*}
	&&\sqrt{N}[G(\xi)-G(\xi+N^{-1/2} \langle  X_N(t,dx',du'), u'\phi(\|x'-x\|)\rangle)] \\
	&&\stackrel{N\to\infty}{\longrightarrow}
	-G'(\xi)\langle X(t,dx'du'), u'\phi(\|x'-x\|) \rangle,
\end{eqnarray*}
when $X_N$ weakly converges to $X$. This shows that $X$ is solution of the martingale problem 
associated to the Markov diffusion process (\ref{eq:linearized Fokker-Planck brown}).
Then the uniqueness of the limit is ensured 
	by the uniqueness to the solution of the limit martingale problem.  
\finproof

The equation (\ref{eq:linearized Fokker-Planck brown}) for the fluctuation process $X$ involves the same
linearized operator as in (\ref{eq:linearized Fokker-Planck}).
So the linear stability analysis carried out in the previous section can be invoked to claim
that the system (\ref{eq:model}) is linearly stable when 
the initial empirical distribution is initially close to $\rho_\xi$,
where $\xi$ is such that $G(\xi)=\xi$, $G'(\xi) <1$,
$|\phi_k G'(\xi)|<1$ for all $k\neq 0$,
and $\sigma$ is large enough.
When $G$ is such that $u-G(u)$ derives from an even double-well potential, such as the two examples (\ref{eq:exG1}) and (\ref{eq:exG2}), this means that the order states $\rho_{\pm\xi_e}$
are stable equilibria, in the sense that if the initial empirical distribution is close to one of them, then the empirical
distribution remains close to it, with small normal fluctuations described by (\ref{eq:linearized Fokker-Planck brown}) in the limit $N\to \infty$,
while the disorder state $\rho_0$ is an unstable equilibrium, in the sense that 
if the initial empirical distribution is close to it, then the empirical distribution quickly moves away from it.

\section{Large Deviations Analysis}
\label{sec:large deviation}

Here we assume that two order states $\rho_{\pm\xi_e}$ exist and are stable
and that the initial empirical measure $\mu_N(t=0,dx,du)$ converges to $\rho_{\xi_e}$. By the mean-field theory, the empirical measure
(\ref{eq:empirical measure})
converges to $\rho_{\xi_e}(x,u)dxdu$ as $N\to\infty$ and the stability analysis 
ensures that $\mu_N(t,dx,du)$ is 
approximately $\rho_{\xi_e}$  over the time interval $[0,T]$.  
However, when $N$ is large but finite, the empirical measure $\mu_N$ is still stochastic and  
there is a very small probability that, eventhough $\mu_N\approx\rho_{\xi_e}$ at time $t=0$, $\mu_N\approx\rho_{-\xi_e}$ at time $t=T$.  
Here we use the large deviation principle (LDP) to write down the asymptotic probability of a transition from one stable
order state to the other one. 
The empirical measure $\mu_N$ satisfies the large deviation principle in the space of continuous probability-measure-valued processes as the following: for a set $A$ of paths of probability measures $\mu=(\mu(t,dx,du))_{t\in [0,T]}$, we have
\begin{align*}
	- \inf_{\mu\in \mathring{A}} I(\mu ) &\leq \liminf_{N\to\infty} \frac{1}{N} \log\mathbf{P}(\mu_N \in A)\\
	&\leq \limsup_{N\to\infty} \frac{1}{N}\log\mathbf{P}(\mu_N \in A) \leq - \inf_{\mu\in\bar{A}} I(\mu ),
\end{align*}
where $\mathring{A}$ and $\bar{A}$ are the interior and closure of $A$, respectively, with an appropriate topology \cite{Dawson1987}.
The function $I$ is called the rate function:
\begin{equation}
	\label{eq:rate function}
	I(\mu ) = \frac{1}{2\sigma^2} \int_0^T \sup_{f(x,u):\langle \mu(t,\cdot,\cdot), (\frac{\partial f}{\partial u})^2 \rangle \neq 0}
	\frac{\langle \frac{\partial \mu}{\partial t}(t,\cdot,\cdot) - \mathcal{L}_{\mu(t,\cdot,\cdot)}^* \mu(t,\cdot,\cdot), f\rangle^2}{\langle \mu(t,\cdot,\cdot), (\frac{\partial f}{\partial u})^2 \rangle} dt,
\end{equation}
where $\mathcal{L}_\nu^*$ is the differential operator associated to the Fokker-Planck equation:
\begin{equation*}
	\mathcal{L}_\nu^*\mu 
	= - u\frac{\partial \mu}{\partial x}
	- \frac{\partial}{\partial u}\left[\left(G\left(\iint u'\phi(\|x'\|) \nu(x-x',u') du'dx' \right) - u \right) \mu\right]
	+ \frac{1}{2} \sigma^2 \frac{\partial^2 \mu}{\partial u^2}.
\end{equation*}

We are interested in the rare event $A$ which corresponds to the transition from one order state $\rho_{\xi_e}(x,u)dxdu$ to the other one $\rho_{-\xi_e}(x,u)dxdu$.
We therefore assume that the initial conditions correspond to independent agents with the distribution $\rho_{\xi_e}(x,u)dx du$ and we consider the event
\begin{equation}
	\label{eq:event of transitions}
	A = \left\{ (\mu(t,dx,du))_{t \in [0,T]} \, : \, \vvvert \mu(T,dx,du)- \rho_{-\xi_e}(x,u) dx du\vvvert \leq \delta\right\},
\end{equation}
for some small $\delta$, where $\vvvert {\cdot} \vvvert$ stands for the metric of the space of probability measures.  

Roughly speaking, for large but finite $N$, the probability of transitions from one stable state to the other is an exponential function of $N$ whose exponential 
decay rate is governed by the rate function $I$:
\begin{equation}
	\label{eq:transition probability}
	\mathbf{P}(\mu_N \in A) \overset{N\gg 1}{\approx} \exp\left(-N\inf_{\mu\in A} I(\mu)\right).
\end{equation}
The transition probability (\ref{eq:transition probability}) of the switching of the alignments of locusts is experimentally and numerically observed in 
\cite{Buhl2006}.

\begin{remark}
	In this section we have only provided a formal large deviation description, because the classical large deviation result \cite{Dawson1987} requires that the noise must be 
	non-degenerate.  Therefore, the result of \cite{Dawson1987} cannot be directly applied to our model and we are not able to well define the space for our case 
	until the rigorous large deviation is proven.  However, a rigorous large deviation principle could be possibly constructed by the technique in \cite{Budhiraja2012}.
\end{remark}

\section{Numerical Simulations}
\label{sec:numerics}

We use the numerical simulations to illustrate and verify our theoretical analysis.  The spatial domain is the torus $[0,L]=[0,10]$, the influence function is 
\begin{equation}
\label{def:phi}
	\phi(\|x\|)= 5 \times {\bf 1}_{[0, 1]}(\|x\|),
\end{equation}
which is such that $\int_0^L \phi(\|x\|)dx / L=1$ and $\phi_k= {\rm sinc}( \pi k / 5)$, and we let 
\begin{equation}
\label{def:G}
	G(u) = \frac{h+1}{5} u- \frac{h}{125} u^3 ,
\end{equation}
which is such that $u-G(u)$ derives from the potential $V(u)=\frac{4-h}{10} u^2 +\frac{h}{500} u^4$, that is a double-well potential as soon as $h>h_c:=4$.
The parameter $h$ allows us to quantify the depths of the two wells of the double-well potential when $h>h_c$ (see Figure \ref{fig:growthrate1}c).

We discretize the time domain: $t_n=n\Delta t$, $n=0, 1, 2, 3, \ldots$ and use the Euler method to simulate (\ref{eq:model}):
\begin{align*}
	 x_i^{n+1}  -x_i^n &= u_i^n \Delta t , \\
	u_i^{n+1} -u_i^n &= \left[G(\left<u\right>_i^n) - u_i^n \right] \Delta t+\sigma \Delta W_i^n,
\end{align*}
where $\{\Delta W_i^n\}_{i,n}$ are independent Gaussian random variables with mean $0$ and variance $\Delta t$, and $\left<u\right>_i^n$ is a weighted average of the 
velocities with respect to $u_i^n$:
\begin{equation}
	\left<u\right>_i^n = \frac{1}{N} \sum_{j=1}^N u_j^n\phi( \|x_j^n-x_i^n\|).
\end{equation}

\subsection{Existence of the order states}

\begin{figure}
	\centering
	{\bf a)}
	\includegraphics[width=0.38\linewidth]{\chemin/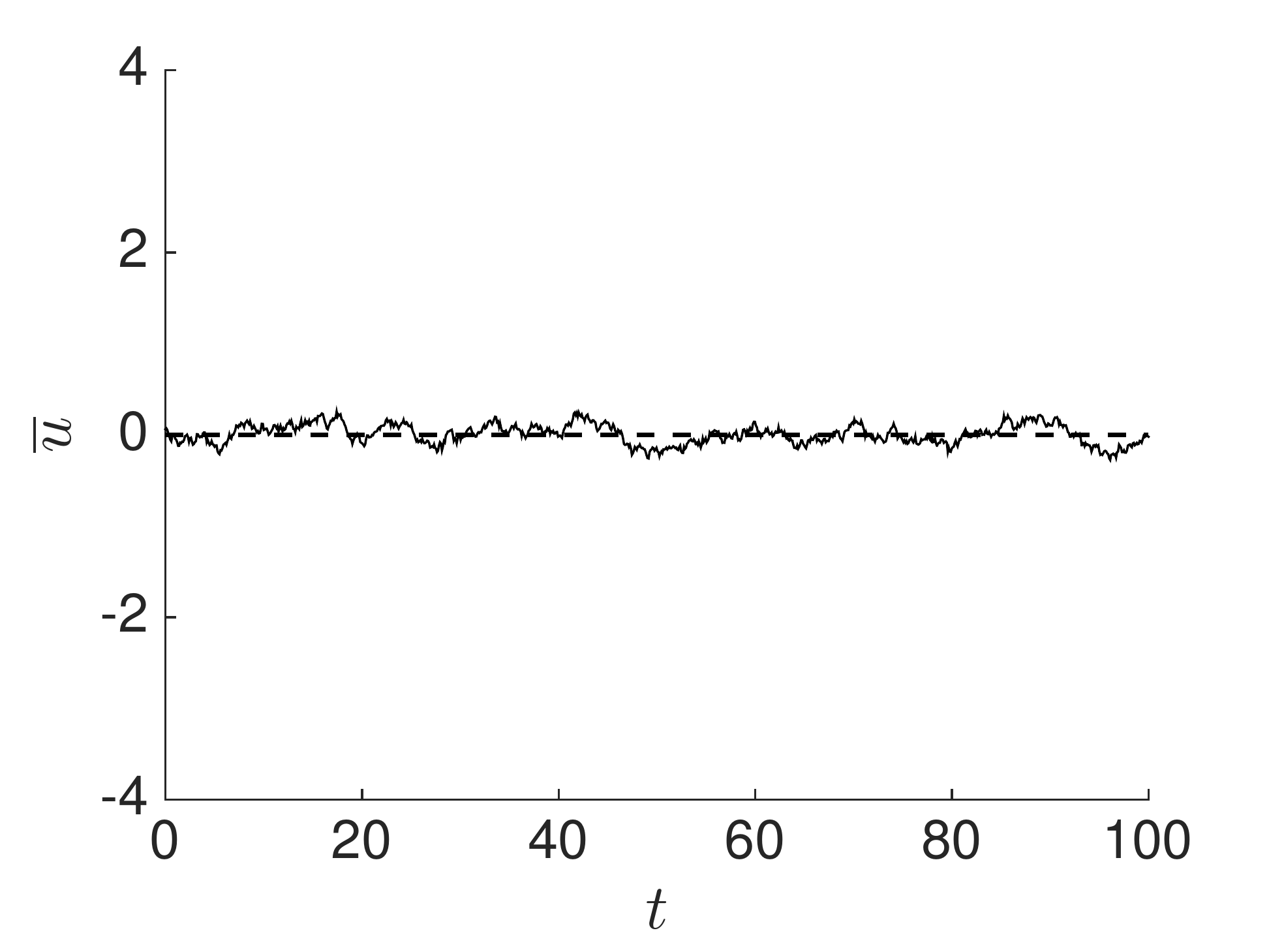}
	{\bf b)}
	\includegraphics[width=0.38\linewidth]{\chemin/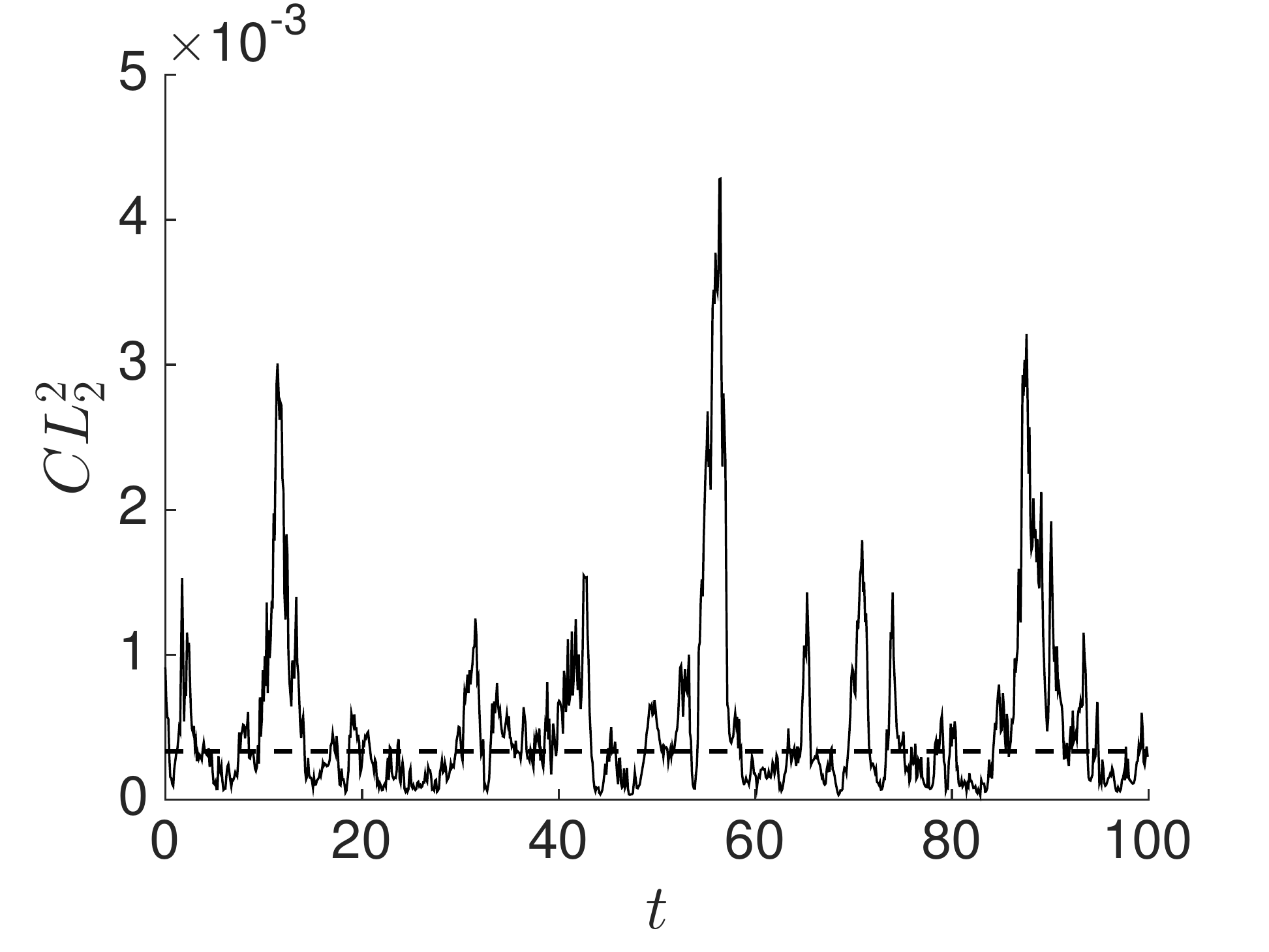}
	
	{\bf c)}
	\includegraphics[width=0.38\linewidth]{\chemin/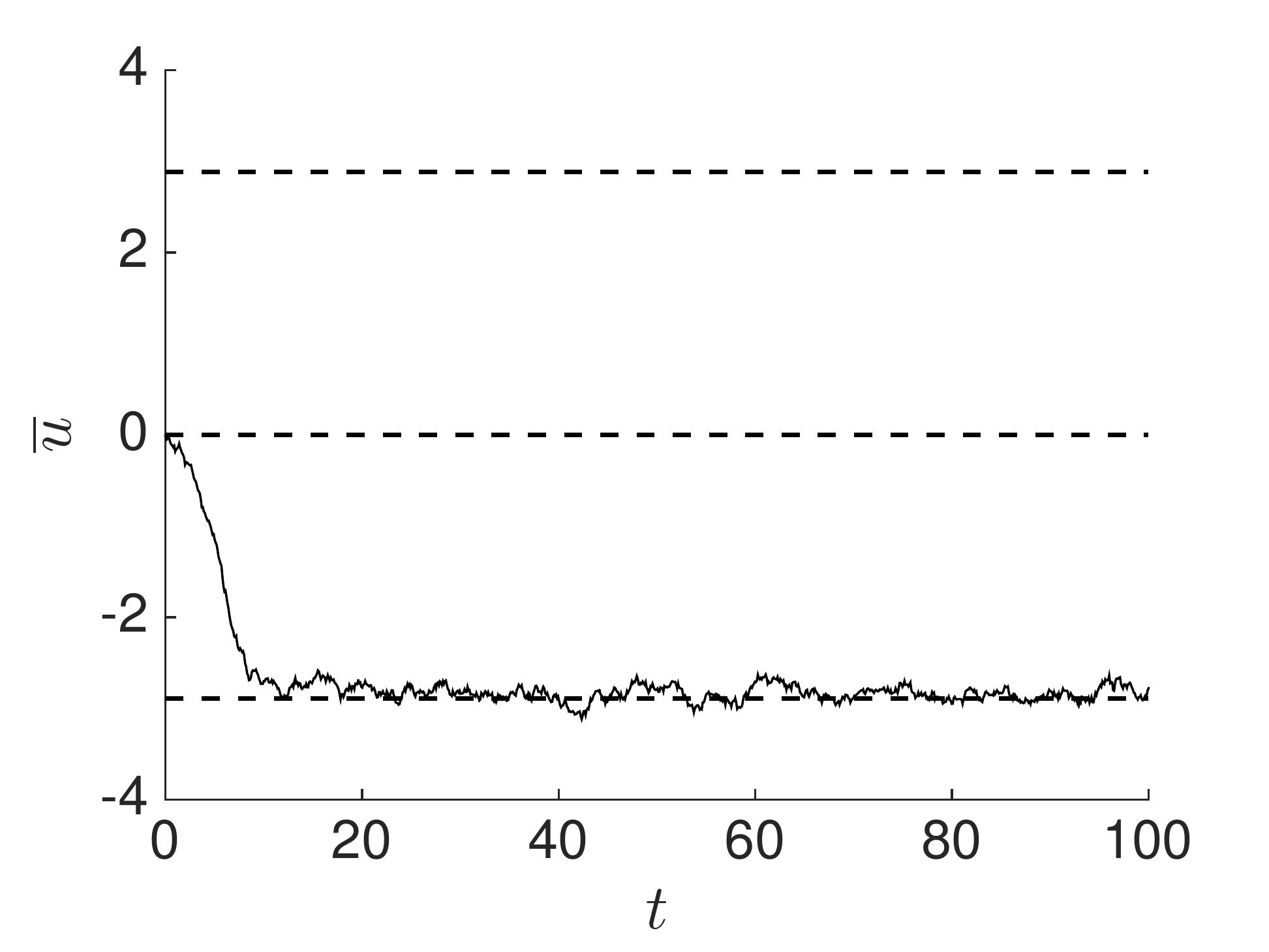}
	{\bf d)}
	\includegraphics[width=0.38\linewidth]{\chemin/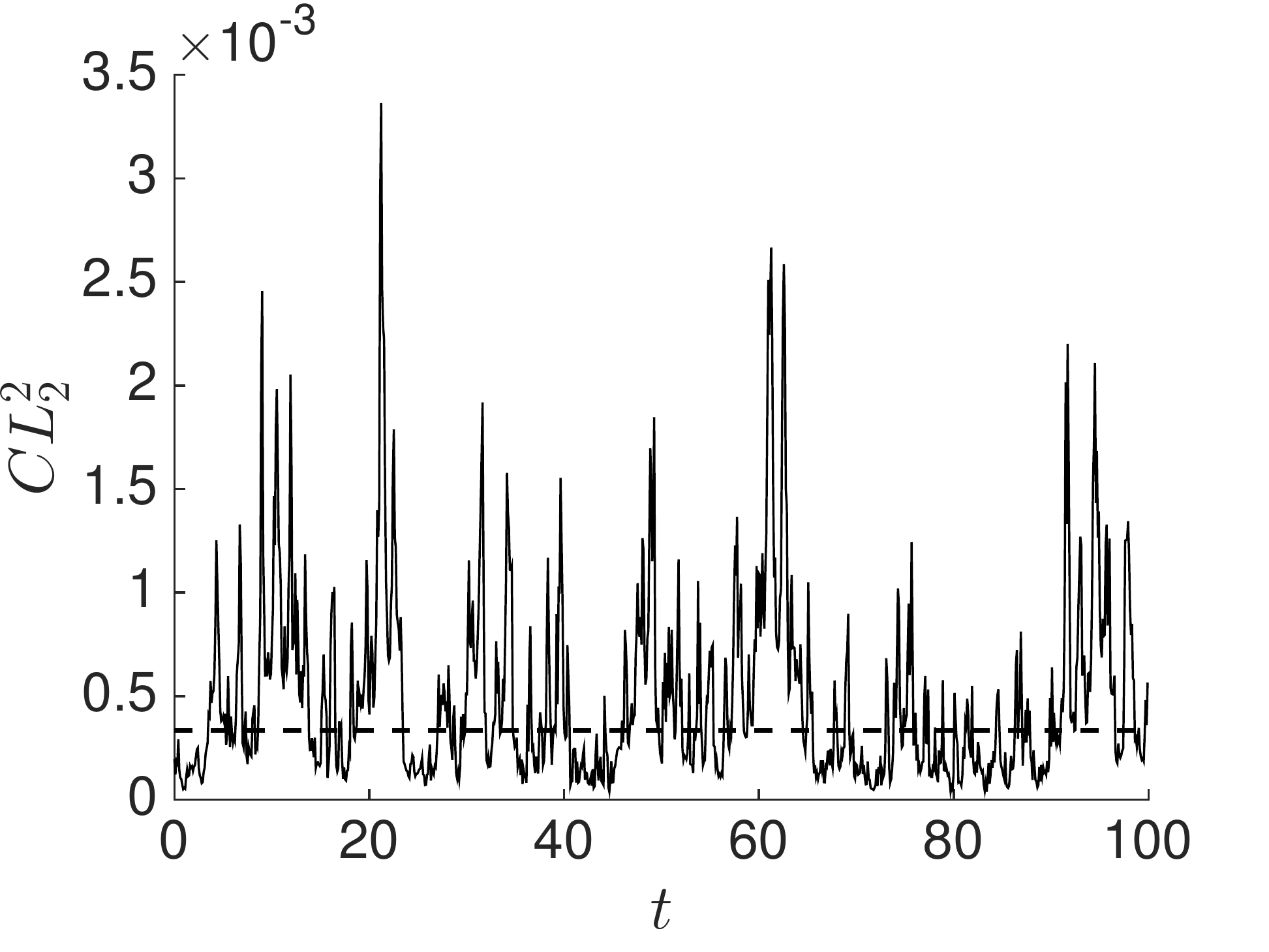}
	
	\caption{The empirical average velocity $\bar{u}^n$ and the  square centered $L^2$-discrepancy $CL_2^2(n)$ at each time step $t_n$ for $h=2$ (a-b) and $h=6$ (c-d). 
The dashed lines in the velocity plots are the solutions of $G(\xi)=\xi$.
The dashed lines in the discrepancy plots stand for the value (\ref{def:CL2u}) corresponding to a uniform sampling.
Here $\Delta t=0.1$, $N=500$, and $\sigma=2$.
One can see  that the spatial distribution is uniform and the velocity average is $0$ when $h=2$ and $- \xi_e$ when $h=6$.
}
	\label{fig:three equilibria}
\end{figure}

From Proposition \ref{prop:equilibrium}, the number of stationary states (\ref{eq:nonlinear Fokker-Planck}) 
is equal to the number of solutions $\xi$ for the compatibility condition (\ref{eq:compatibility condition}).
When $h \leq h_c=4$, the function $u-G(u)$ is increasing and  $0$ is the unique solution of the compatibility
equation (\ref{eq:compatibility condition}).
When $h>h_c$, the function $u-G(u)$ derives from a double well-potential and there are three solutions $0,\pm\xi_e$ 
to the compatibility equation (\ref{eq:compatibility condition}), with
\begin{equation}
\xi_e = 5 \sqrt{ \frac{h-4}{h}}.
\end{equation}

In Figure \ref{fig:three equilibria} we test the cases of $h=2$ and $h=6$ and plot the empirical average velocity 
\begin{equation}
\bar{u}^n = \frac{1}{N}\sum_{i=1}^{N}u_i^n
\end{equation}
and the square centered $L^2$-discrepancy:
\begin{equation}
CL_2^2(n) =  \int_0^1 |F_n(x) -x |^2 dx, \quad \quad F_n(x) = \frac{{\rm Card}(i =1,\ldots,N\, , \, x_i^n/L \leq x)}{N} .
\end{equation}
It measures the discrepancy between the empirical distribution of the positions and the uniform distribution over $[0,L]$.
By \cite{hickernell} it is given by
\begin{eqnarray}
\nonumber
CL_2^2(n) &=& \frac{13}{12} - \frac{1}{N} \sum_{i=1}^N  \Big( 2+  \Big|\frac{x_i^n}{L}-\frac{1}{2}\Big| -  \Big|\frac{x_i^n}{L}-\frac{1}{2}\Big| ^2  \Big)
\\
&& +\frac{1}{2 N^2} \sum_{i,j=1}^N \Big( 2+  \Big|\frac{x_i^n}{L}-\frac{1}{2}\Big|  + \Big|\frac{x_j^n}{L}-\frac{1}{2}\Big| -   \Big|\frac{x_i^n}{L}-\frac{x_j^n}{L}\Big|  \Big)  .
\end{eqnarray}
It is known that, if the positions are independent and identically distributed according to the uniform distribution over $[0,L]$,
then the square centered $L^2$-discrepancy has mean \cite{fang02}:
\begin{equation}
\label{def:CL2u}
CL_{2,u}^2  =\frac{1}{N} \Big( \frac{5}{4}-\frac{13}{12}\Big)  ,
\end{equation}
while its variance is of order $N^{-2}$.
Therefore, as long as the empirical distribution in positions stays uniform, 
the  square centered $L^2$-discrepancy stays close to the value~(\ref{def:CL2u}).

In Figure \ref{fig:three equilibria}
the initial locations and velocities $\{(x_i^0, u_i^0)\}_{i=1}^N$ are sampled from the distribution 
$\rho_0(x,u)dx du$, the stationary state in (\ref{eq:equilibrium}) with $\xi=0$.  That is, $\{x_i^0\}_{i=1}^N$ 
are uniformly sampled over $[0,L]$ and $\{u_i^0\}_{i=1}^N$ are sampled from the Gaussian distribution with mean $0$ and variance $\sigma^2/2$.\\
- 
When $h =2 < h_c=4$, the only solution of $\xi = G(\xi)$ is $\xi=0$,
which means that $\rho_0$ is the unique stationary state.
Moreover, $G'(0)=(h+1)/5<1$ and the noise level is strong enough to make it linearly stable.
Indeed, we can see in Figure~\ref{fig:three equilibria}a
that the empirical average velocity $\bar{u}^n$ oscillates around zero
and the square centered $L^2$-discrepancy oscillates around the value $CL_{2,u}^2 $,
which indicates that the empirical distribution of locations and velocities stays close to $\rho_0$.\\
- When $h=6>h_c=4$, the order states $\rho_{\pm\xi_e}$ exist. 
Moreover, $G'( \xi_e) = (13-2h)/5 <  (13-2h_c)=1$ and  the noise level is strong enough to make the order 
states $\rho_{\pm \xi_e}$ linearly stable. 
In addition, $G'(0)  = (h+1)/5 > (h_c+1)/5 = 1$, therefore, the disorder state $\rho_0 $ is an unstable state.
Under these circumstances, the initial empirical distribution that is close to $\rho_0$ quickly changes 
as can be seen in Figure~\ref{fig:three equilibria}c-d:  
$\bar{u}^n$ oscillates around $-\xi_e$ (it could have been $\xi_e$), while  the square centered $L^2$-discrepancy oscillates 
around the value $CL_{2,u}^2 $.
This indicates that the empirical distribution of locations and velocities becomes close to $\rho_{-\xi_e}$.

\subsection{Linear stability}

Here we  assume  that $h>h_c=4$ so that there are three stationary states
$\rho_0$, $\rho_{\pm\xi_e}$.  
We have $G'(0)  = (h+1)/5 > 1$, $G'(\pm \xi_e)   = (13-2h)/ 5< 1$,
and thus the $0$th order mode of $\rho_0$ is unstable while the $0$th-order modes of $\rho_{\pm\xi_e}$ are stable.
An immediate conclusion is that $\rho_0$ is 
an unstable state for the nonlinear Fokker-Planck equation (\ref{eq:nonlinear Fokker-Planck}).
The zeroth-order mode of $\rho_{\pm\xi_e}$ is linearly stable. 
However, the stability of $\rho_{\pm\xi_e}$ requires that all the nonzeroth-order modes of 
$\rho_{\pm\xi_e}$ are linearly stable.  From Proposition \ref{prop:stablity of mode k}, 
we know that for a sufficiently large $\sigma$, all the nonzeroth-order modes 
 are stable. Moreover the critical value of $\sigma$ ensuring the stability can be determined as explained 
 in subsection \ref{subsec:instaspa}.

\begin{figure}
	\centering
	{\bf a)}
	\includegraphics[width=0.38\linewidth]{\chemin/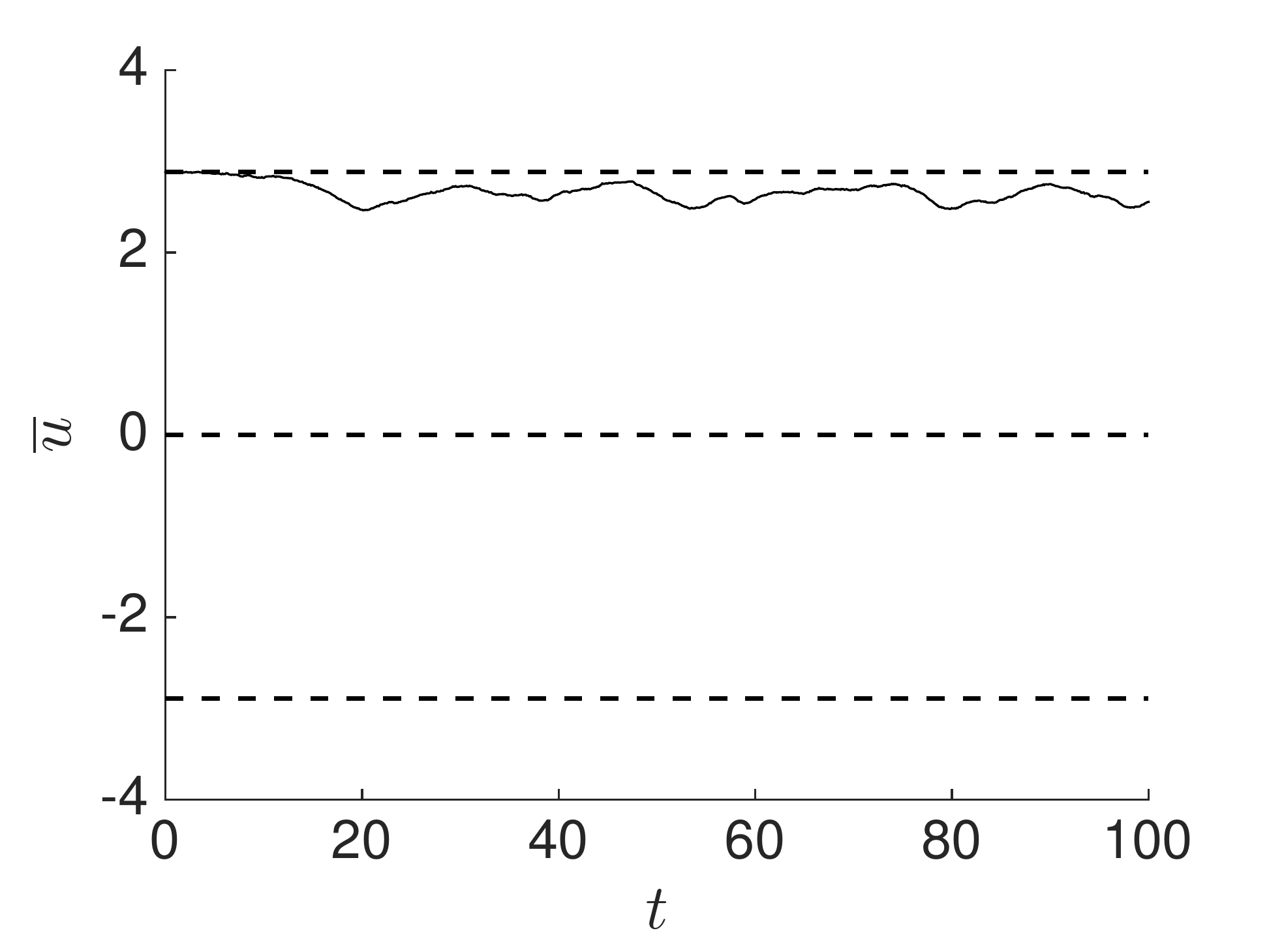}
	{\bf b)}
	\includegraphics[width=0.38\linewidth]{\chemin/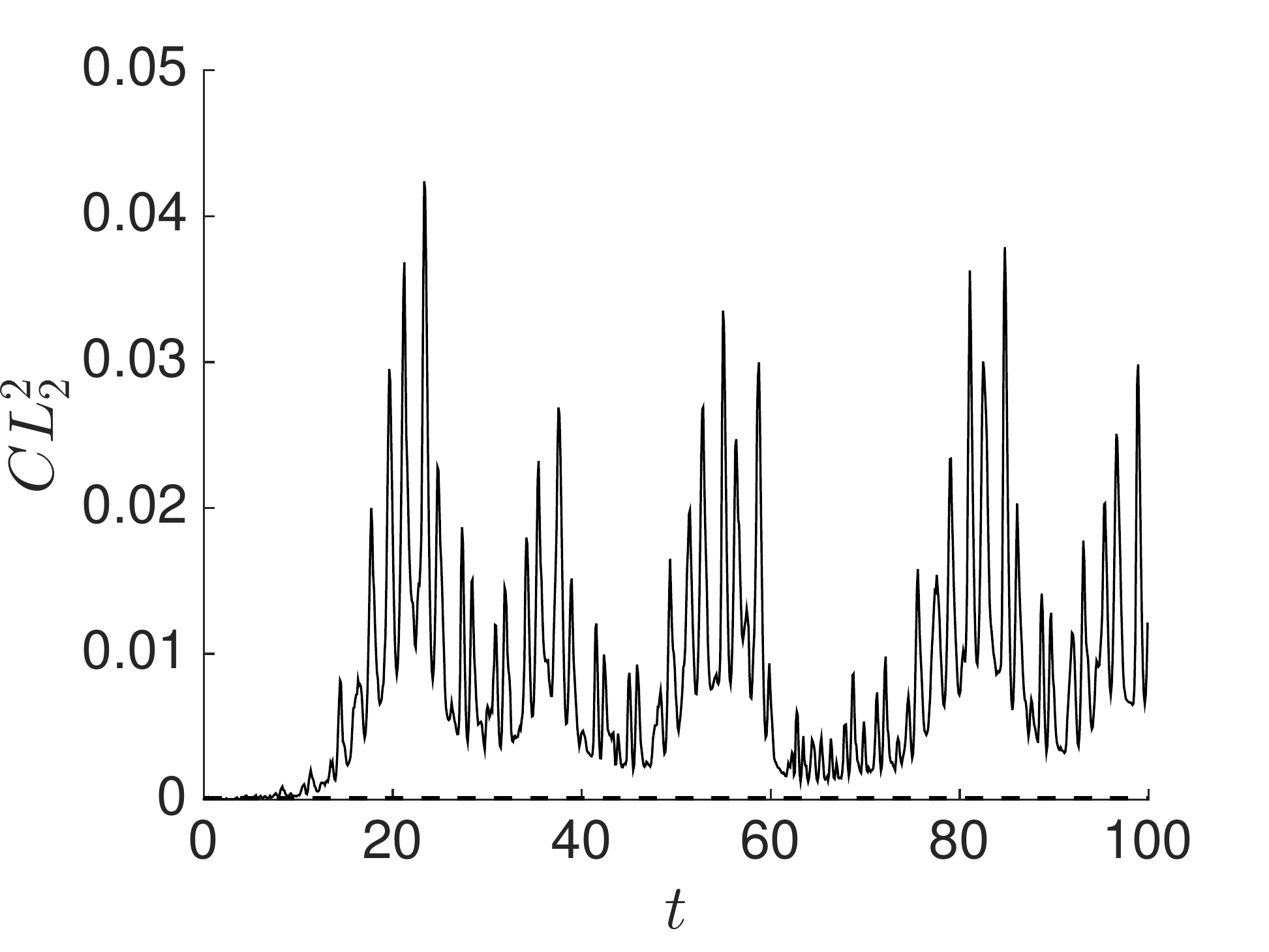}
	
	{\bf c)}
	\includegraphics[width=0.38\linewidth]{\chemin/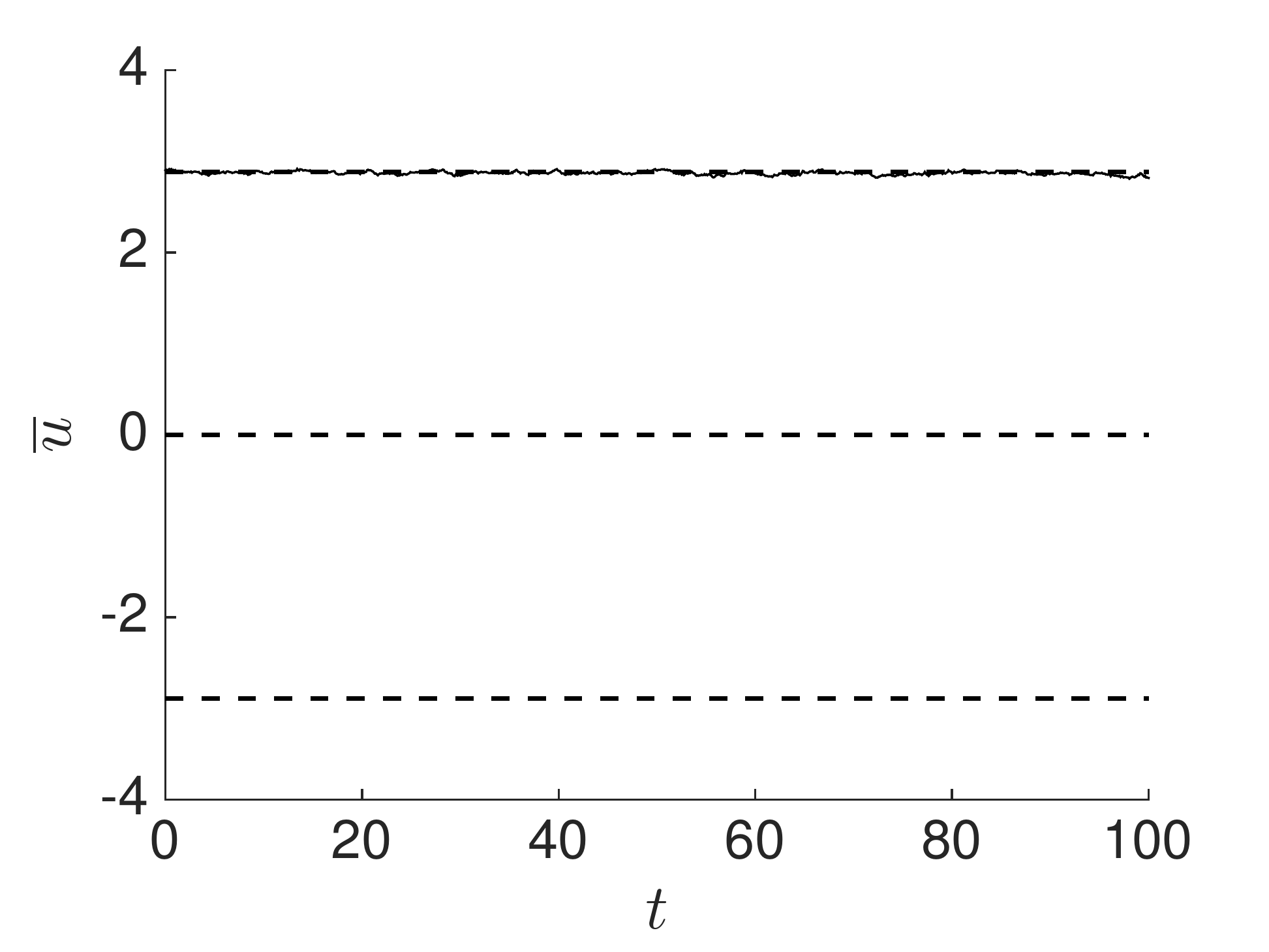}
	{\bf d)}
	\includegraphics[width=0.38\linewidth]{\chemin/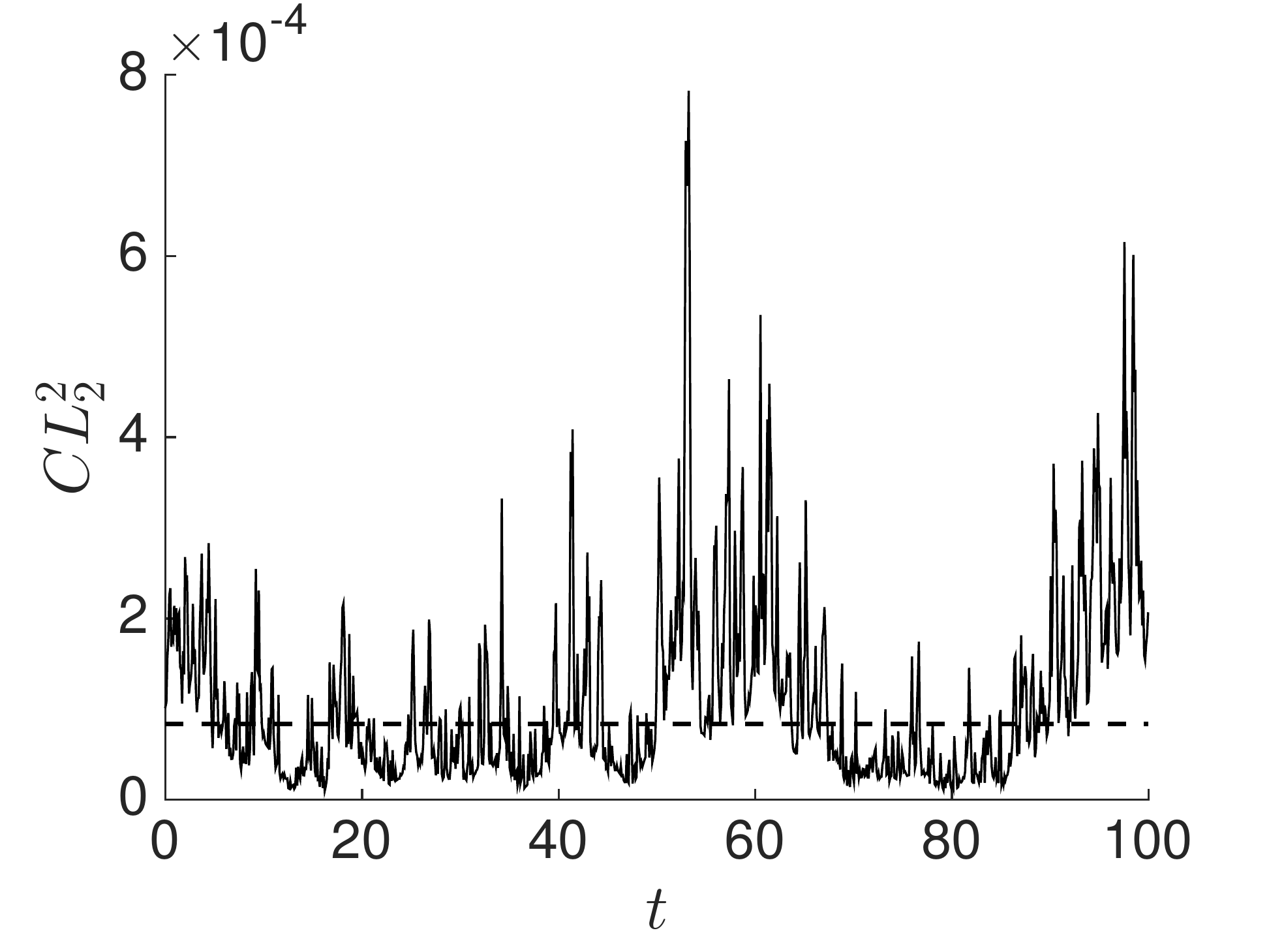}
	
	{\bf e)}
	\includegraphics[width=0.38\linewidth]{\chemin/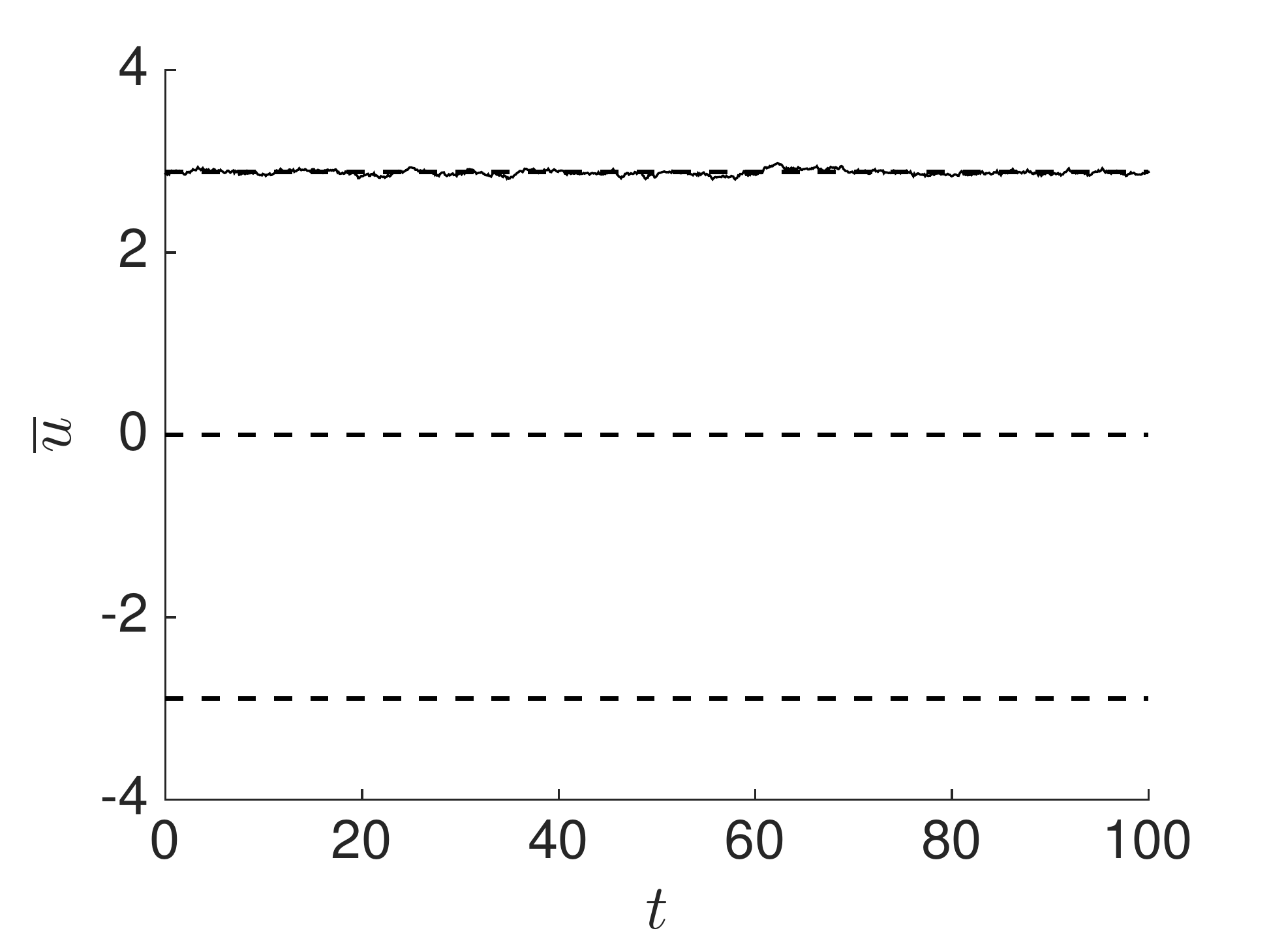}
	{\bf f)}
	\includegraphics[width=0.38\linewidth]{\chemin/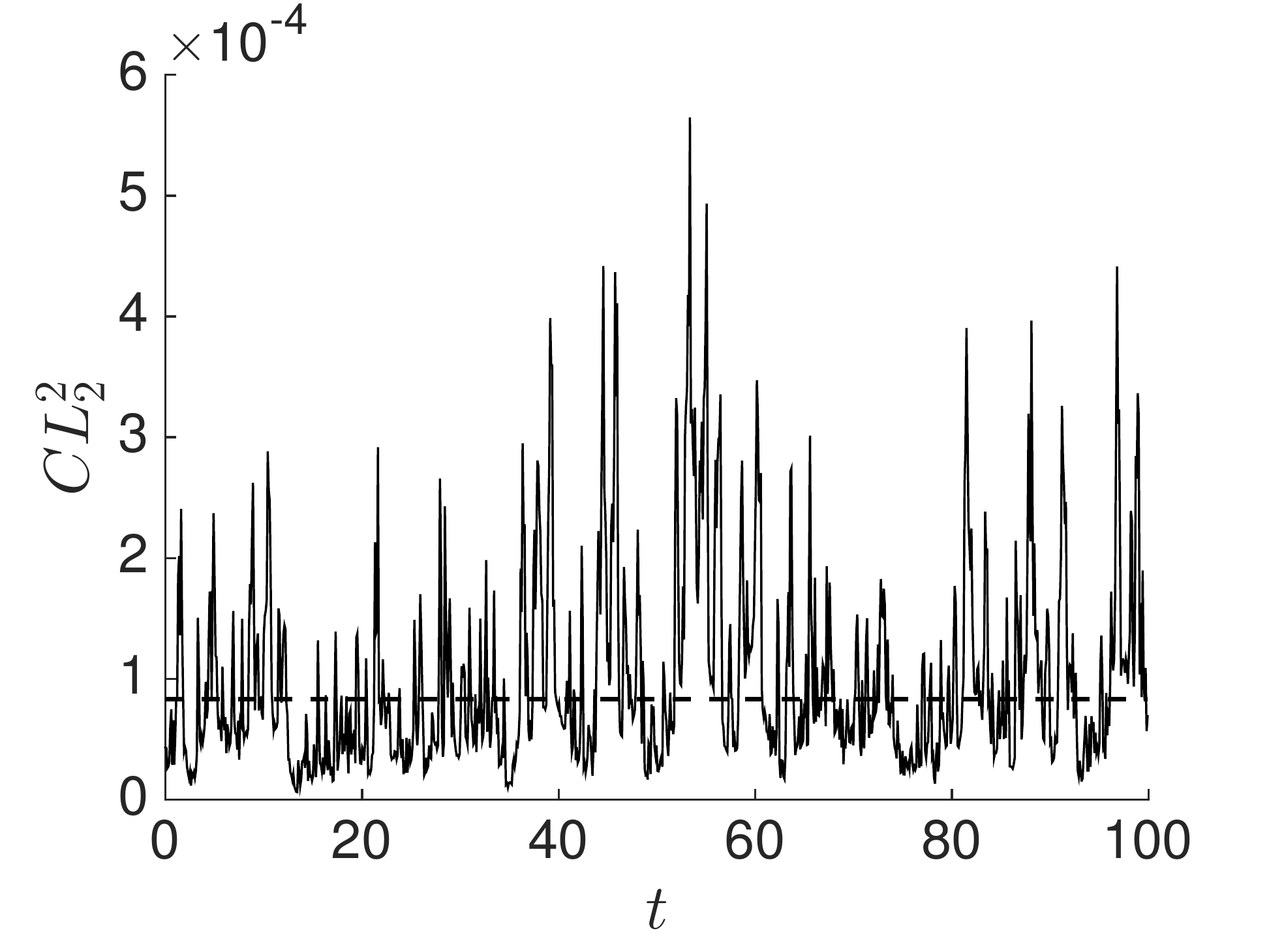}
	
	\caption{The empirical average velocity $\bar{u}^n$ and the  square centered $L^2$-discrepancy $CL_2^2(n)$ 
	for $\sigma=0.5$ (a-b), $\sigma=1$ (c-d), $\sigma=1.5$ (e-f).  
	 The initial positions $\{x_i^0\}_{i=1}^N$ are uniformly sampled over $[0,L]$ and the initial velocities $\{u_i^0\}_{i=1}^N$ are sampled from the Gaussian distribution with mean $\xi_e$ and variance $\sigma^2/2$.  
Here $\Delta t=0.1$, $N=2000$, and $h=6$.
One can see  that the average velocity is  not $\xi_e$ and the spatial distribution is not uniform when $\sigma=0.5$ 
while the average velocity is $\xi_e$ and the spatial distribution is uniform when $\sigma=1$ or $1.5$.
} 
	\label{fig:2}
\end{figure}

\begin{figure}
	\centering
	{\bf a)}
	\includegraphics[width=0.38\linewidth]{\chemin/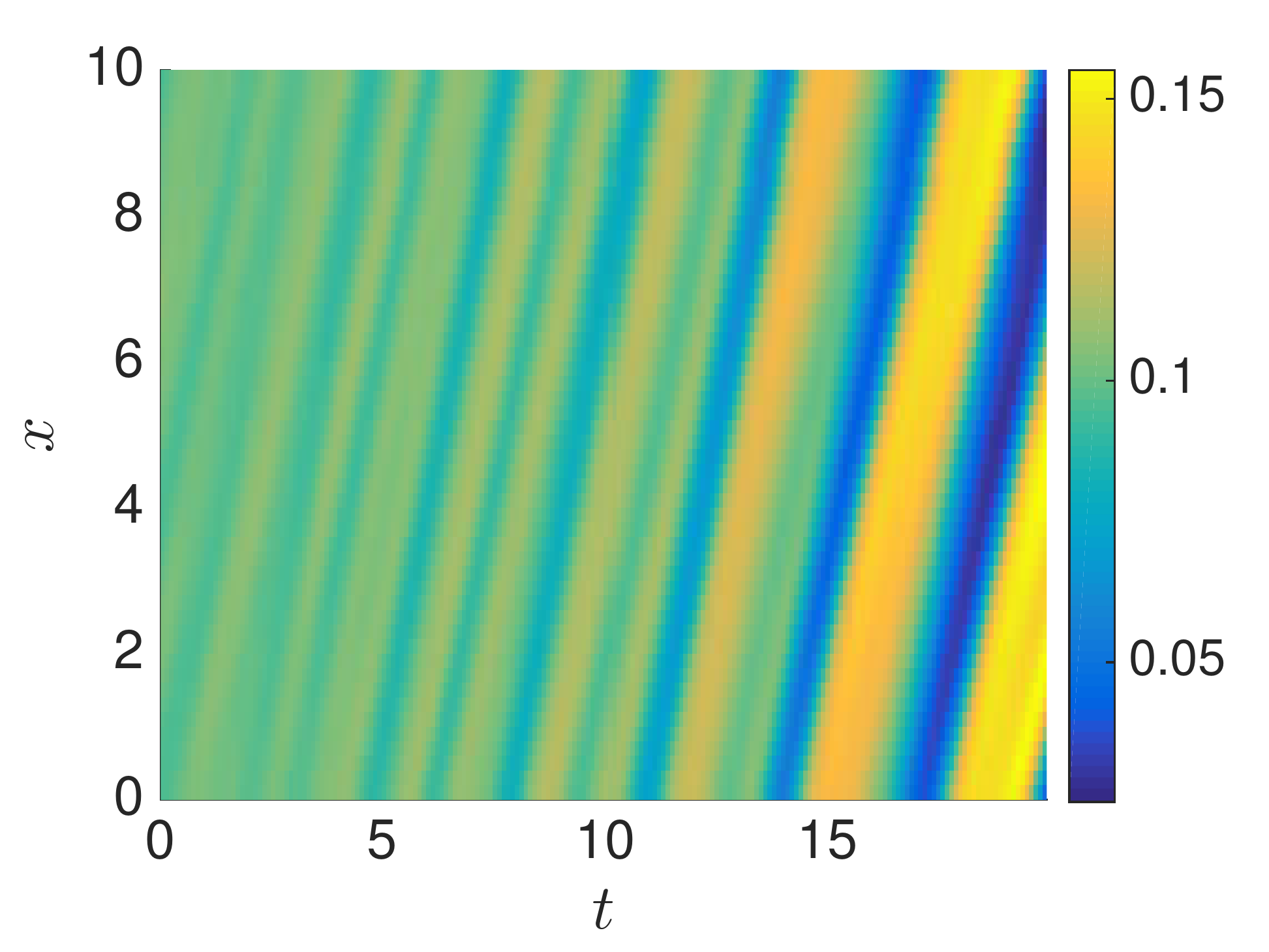}
	{\bf b)}
	\includegraphics[width=0.38\linewidth]{\chemin/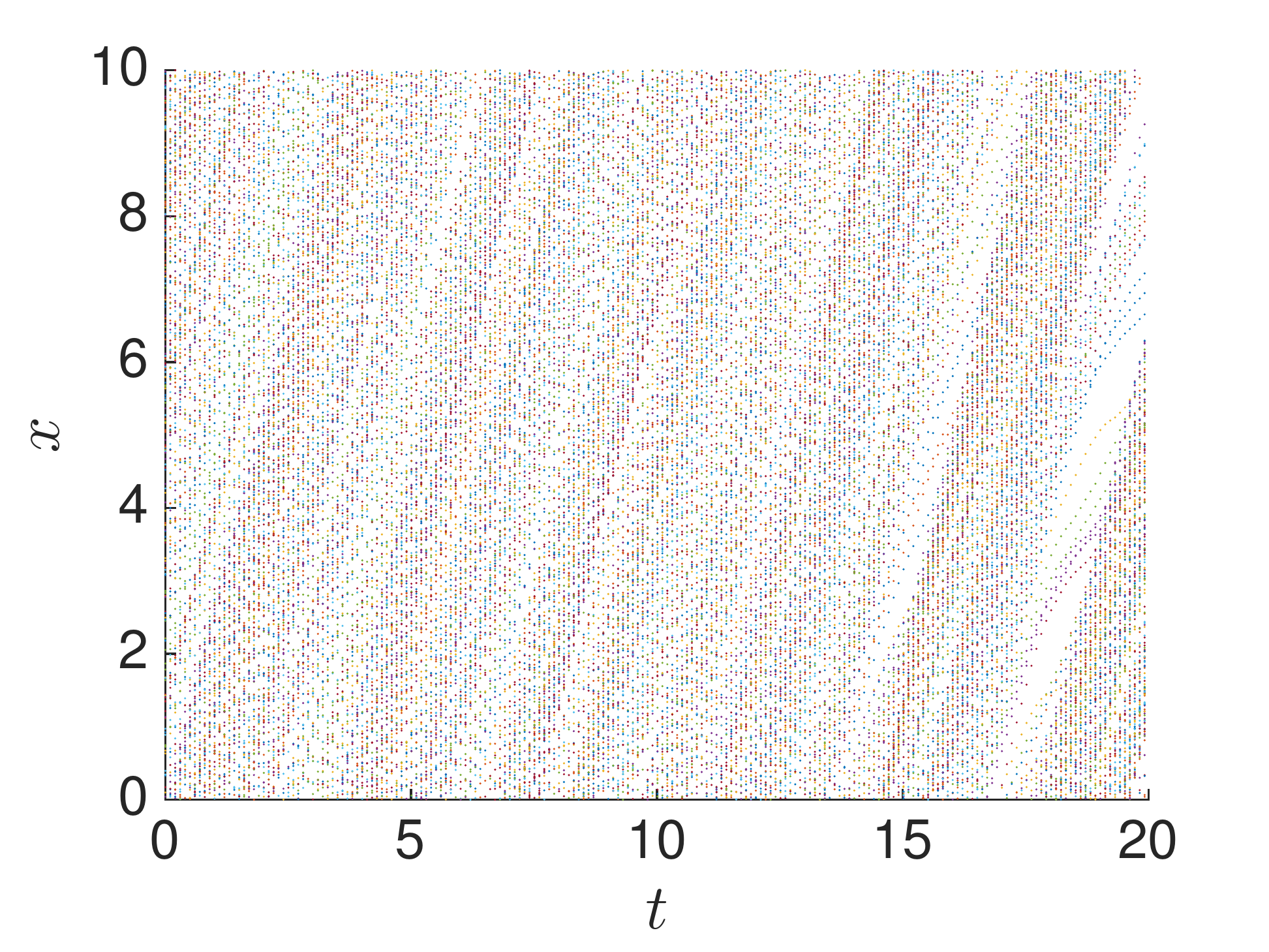}

	{\bf c)}
	\includegraphics[width=0.38\linewidth]{\chemin/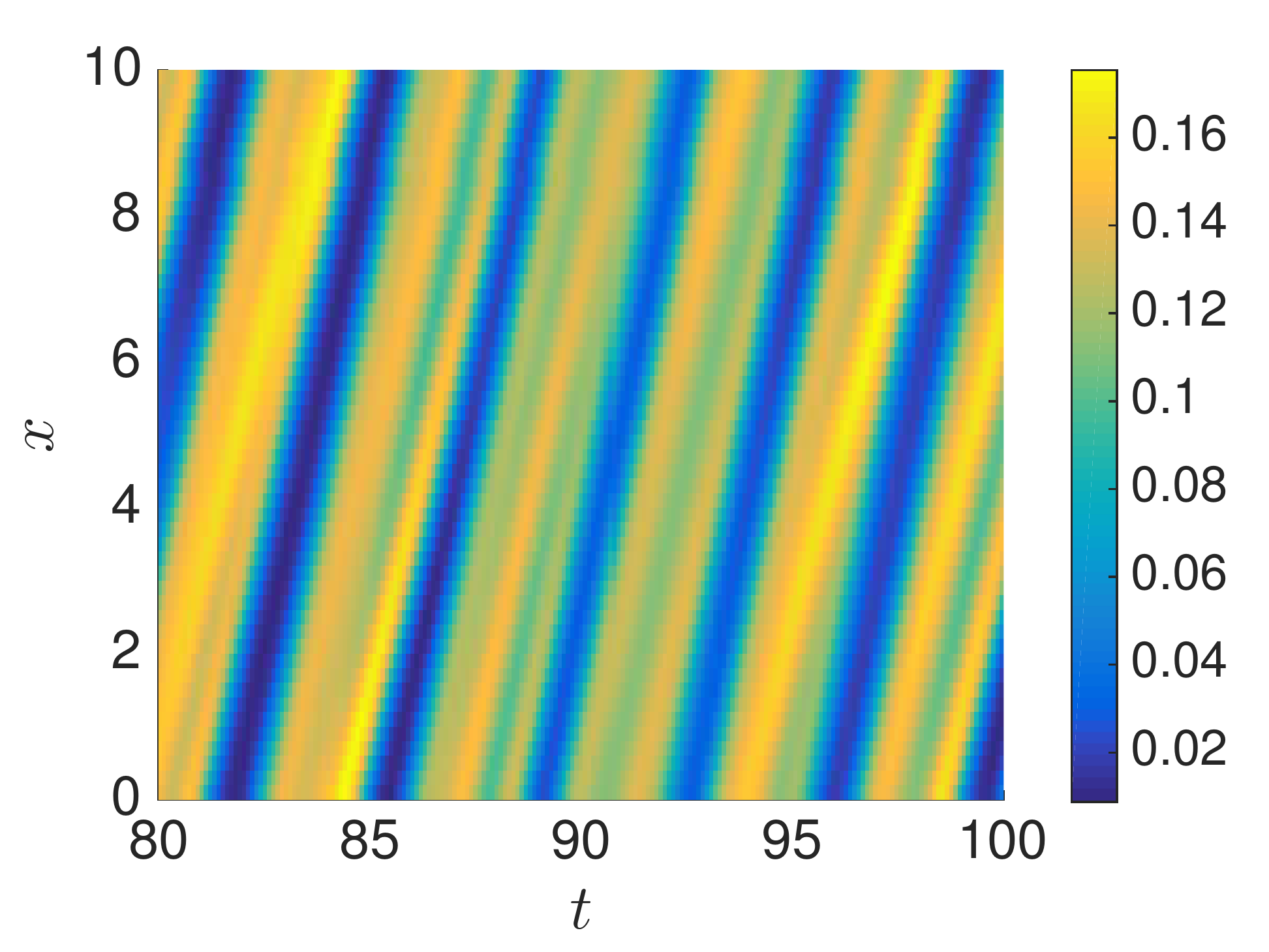}
	{\bf d)}
	\includegraphics[width=0.38\linewidth]{\chemin/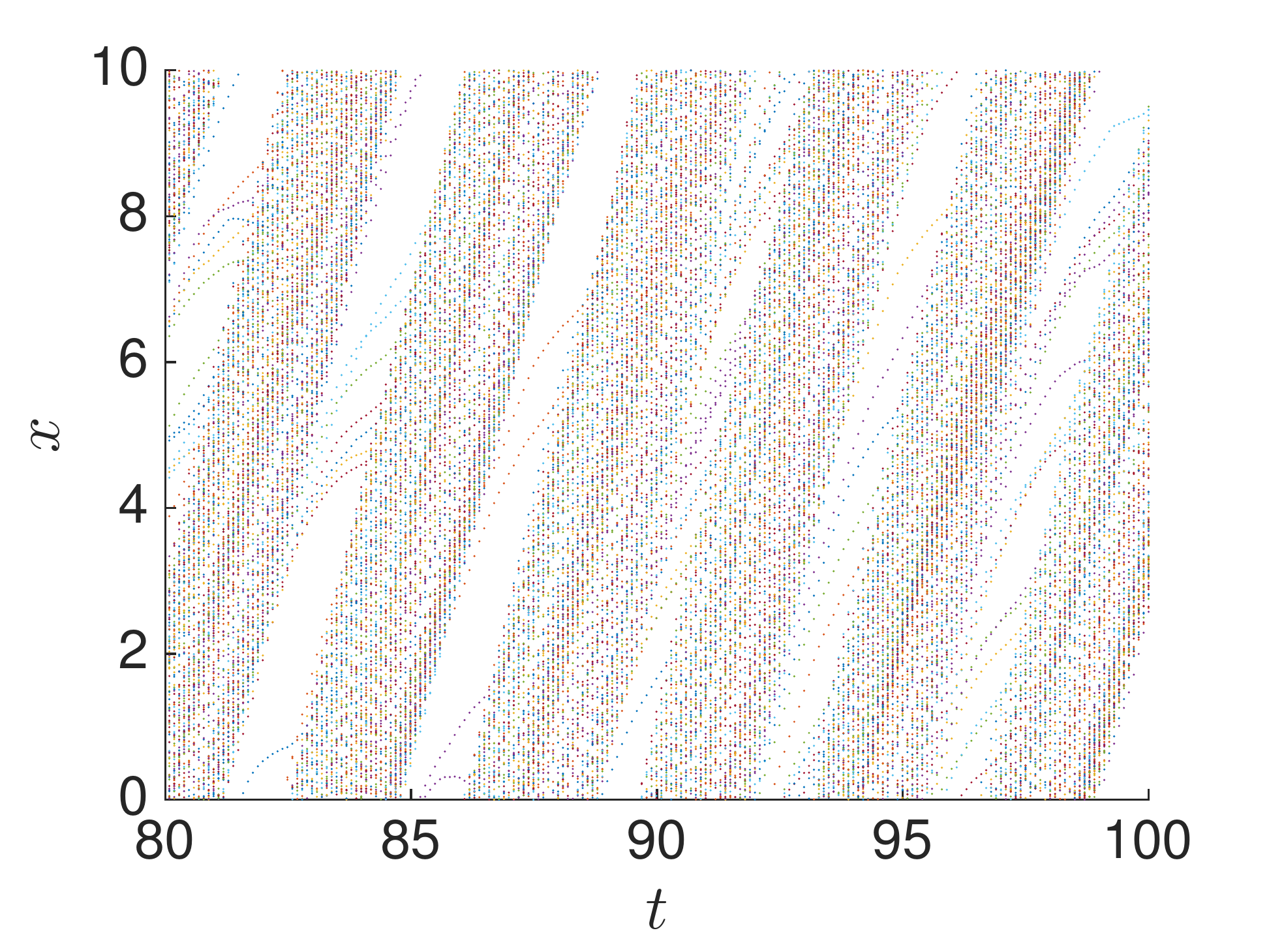}

	\caption{The empirical position distribution smoothed by kernel density estimation (a and c) and the first $200$ trajectories (b and d). 
	The initial positions $\{x_i^0\}_{i=1}^N$  are uniformly sampled over $[0,L]$ and the initial velocities $\{u_i^0\}_{i=1}^N$ are sampled from the Gaussian distribution with mean $\xi_e$ and variance $\sigma^2/2$.  
	Here $\Delta t=0.1$, $N=2000$, $\sigma=0.5$, and $h=6$. 
	One can see  that the spatial distribution is not uniform and a cluster is moving with an apparent velocity of $3.6$.
	}
	\label{fig:3}
\end{figure}

If $h=6$, then the theory in subsection \ref{subsec:instaspa}
predicts that the condition for the stability of the order states  $\rho_{\pm \xi_e}$  
is that $\sigma$ should be larger than $0.85$.
The simulations shown in Figure~\ref{fig:2} confirm these predictions:
$\rho_{\pm \xi_e}$ is stable when $\sigma$ exceeds the  threshold value $0.85$.
They also allow us to illustrate the instability mechanism exhibited above: 
when $\sigma$ is below the threshold value, 
the average velocity of the particles is linearly stable, but
the uniform distribution of the positions of the particles is not stable and a spatial modulation 
grows that gives rise to one moving cluster.
Note that Figure \ref{fig:2}a indeed shows 
that the uniform distribution of the positions is not stable when $\sigma=0.5$ since the square $L^2$-discrepancy takes values much larger than
(\ref{def:CL2u}).
The theory in subsection \ref{subsec:instaspa} predicts that
the most unstable mode is the first one $k_{\rm max}=1$,  which means that there should be one moving cluster.
Figure \ref{fig:growthrate1} plots the predicted growth rate $\gamma_r(1)$ of the first mode and the coefficient $\gamma_i(1)$.
For $\sigma=0.5$, we have $\gamma_i(1)\simeq 2.15$ and the apparent velocity for the cluster is predicted 
to be approximately equal to $\gamma_i(1)\times L / (2\pi) \simeq 3.4$ by (\ref{eq:predictvel}), 
which is close to, although slightly larger than, the average velocity of the particles in the stationary distribution $\xi_e=2.9$.
There is no contradiction as the apparent velocity of the cluster can be interpreted as a ``phase velocity"
while the average velocity of the particles  can be interpreted as a ``group velocity".
The growth of the first mode giving rise to the moving cluster can be observed in the simulations (see Figure \ref{fig:3}). 
Moreover, we can see in the right pictures of Figure \ref{fig:3} that the particles in the tail of the moving cluster
are transfered to the front of the ``next" cluster (which is the same one, by periodicity), 
which indeed allows the cluster to apparently move faster than the average velocity of the particles: in the simulations, the
 average velocity of the particles is about $2.8$, while the  velocity of the cluster is about $3.6$,
 very close to the predicted value~$3.4$.

\begin{figure}
	\centering
	{\bf a)}
	\includegraphics[width=0.38\linewidth]{\chemin/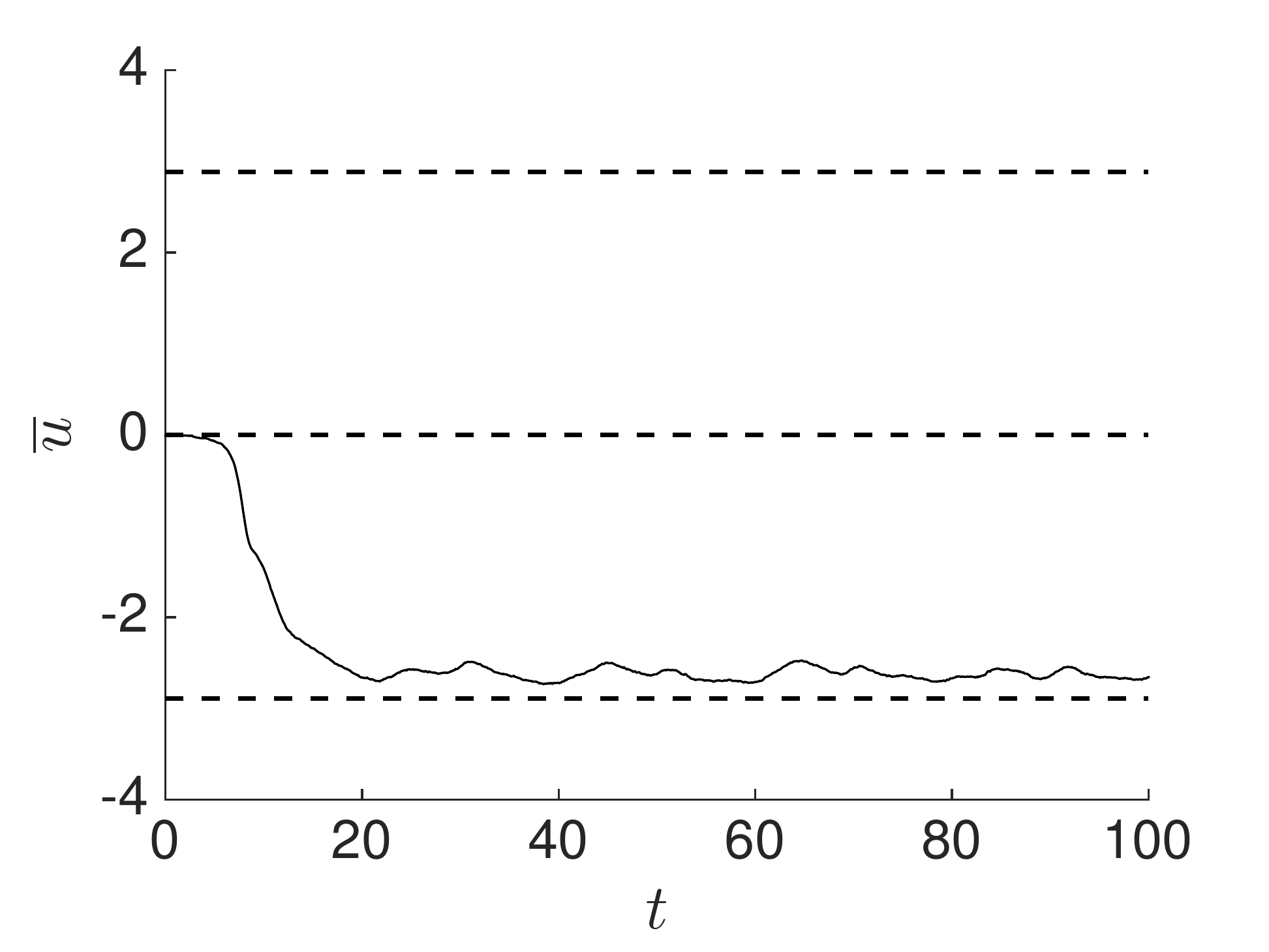}
	{\bf b)}
	\includegraphics[width=0.38\linewidth]{\chemin/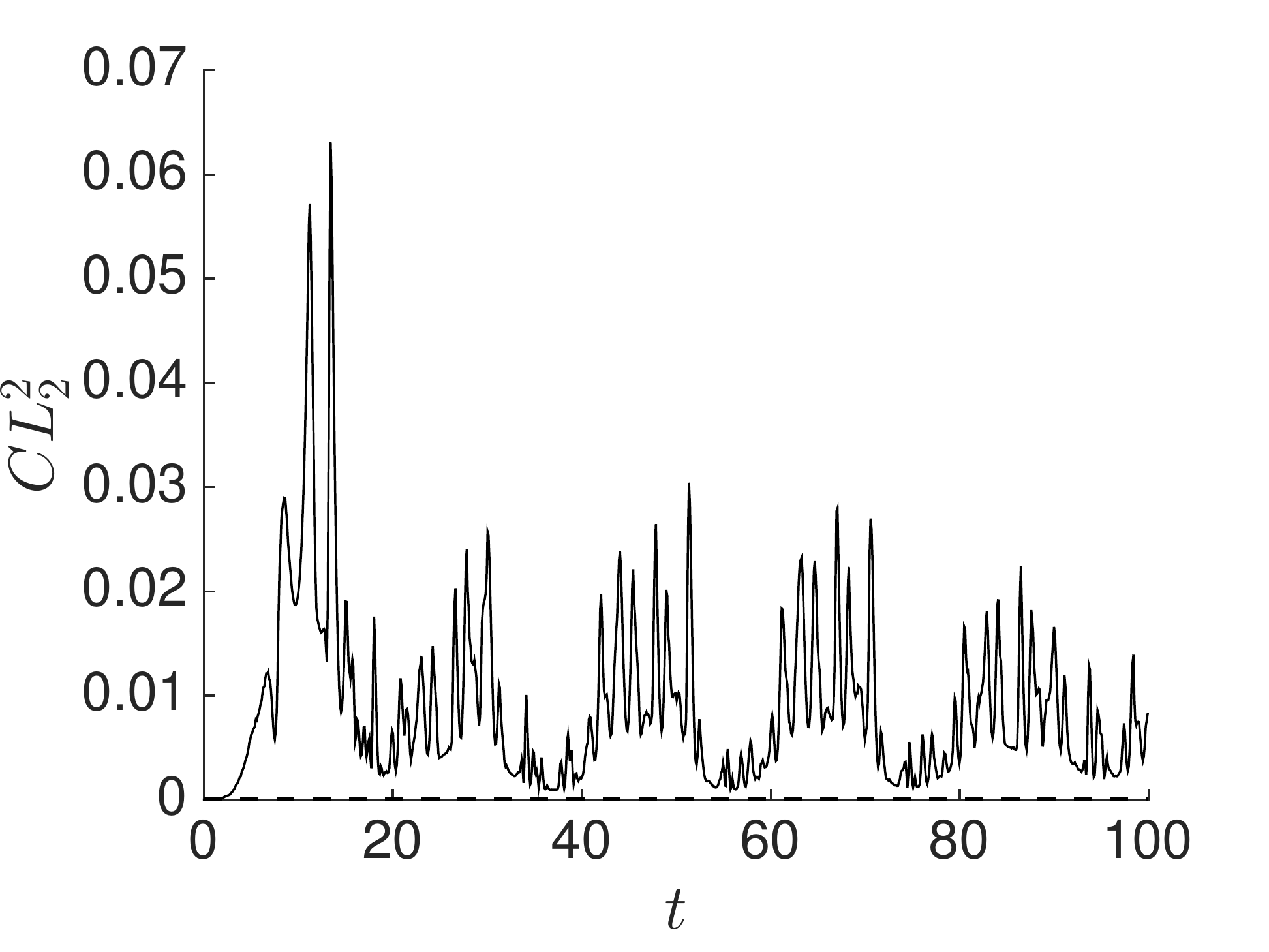}
	
	{\bf c)}
	\includegraphics[width=0.38\linewidth]{\chemin/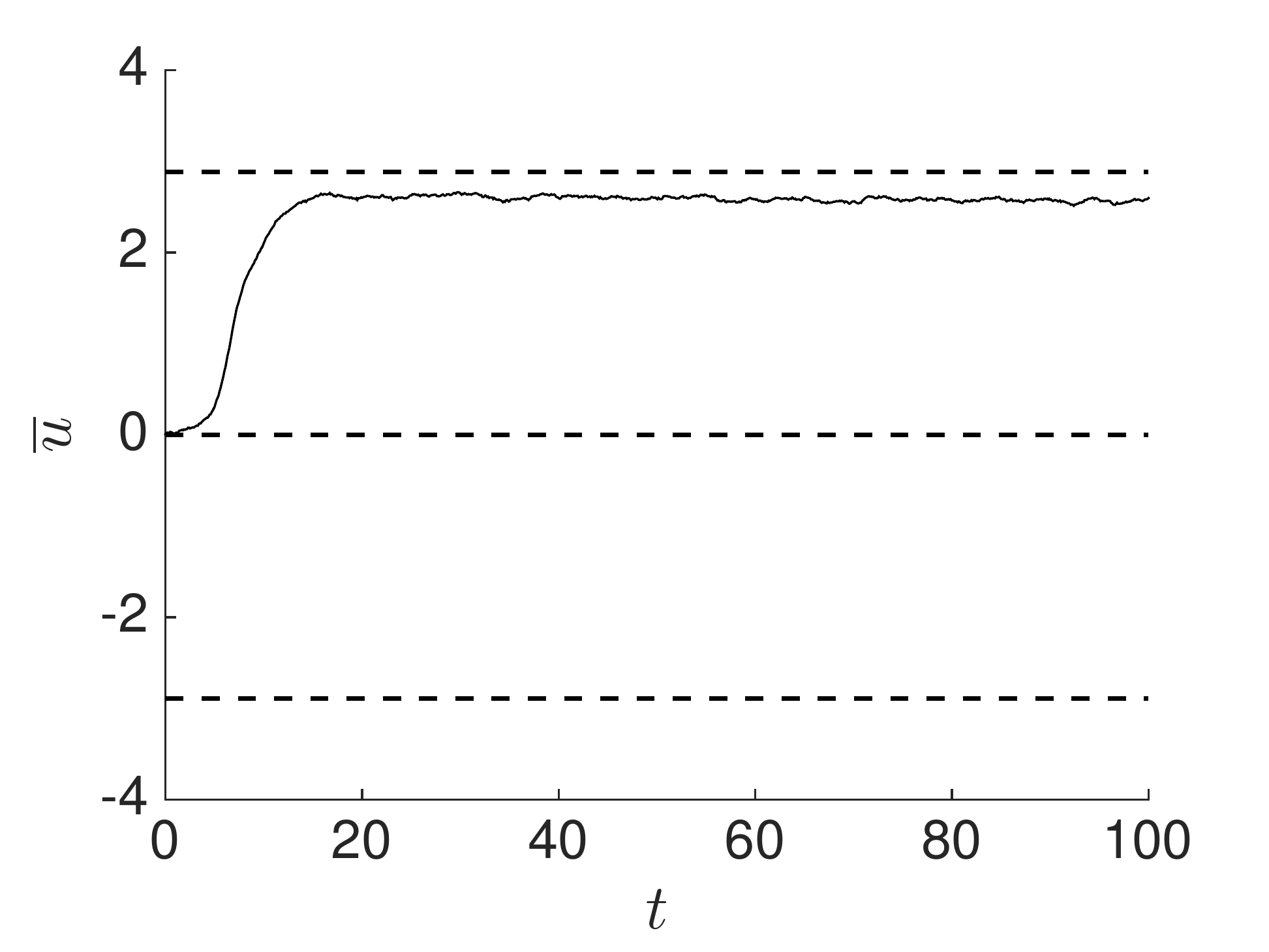}
	{\bf d)}
	\includegraphics[width=0.38\linewidth]{\chemin/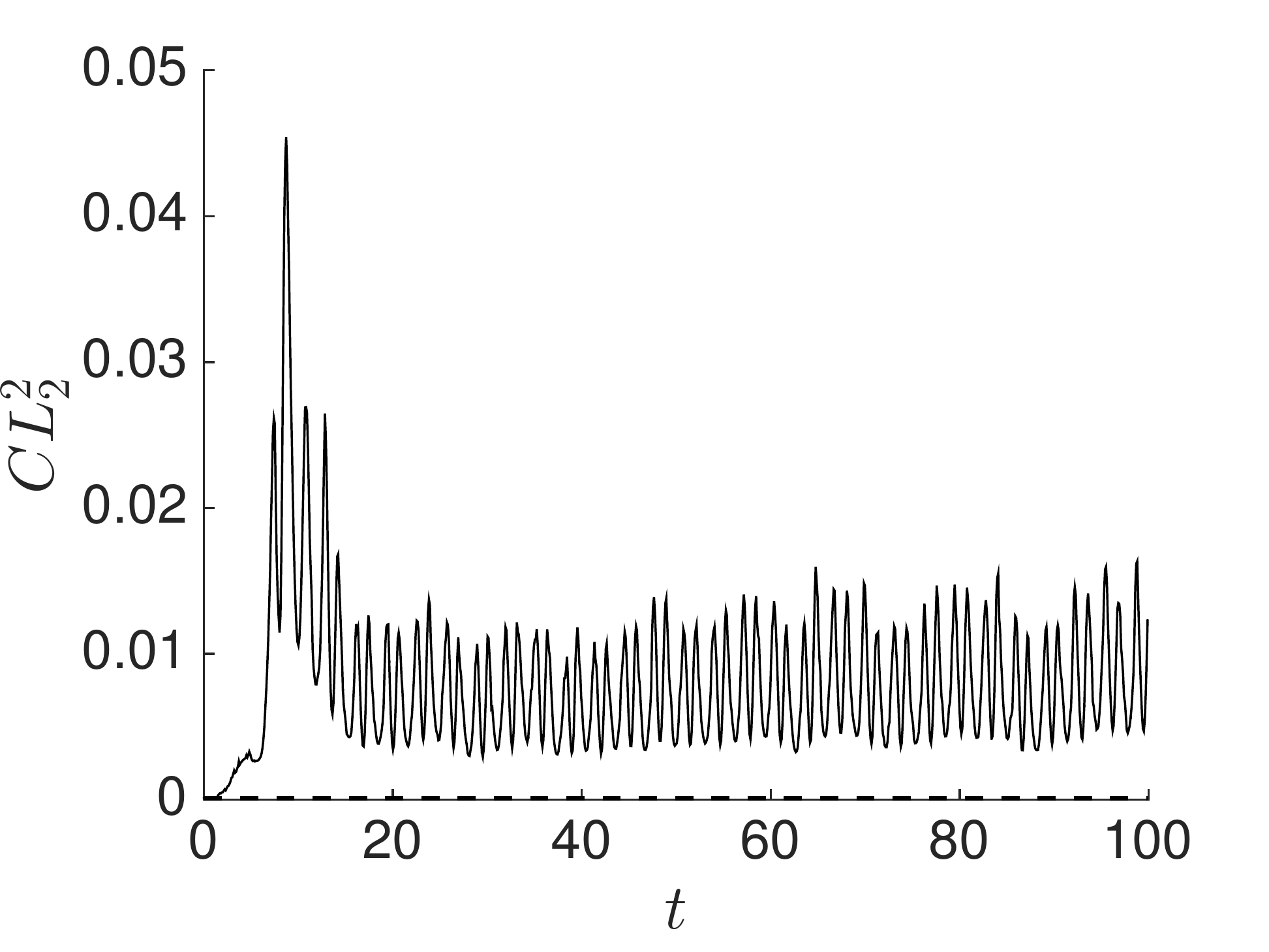}
	
	{\bf e)}
	\includegraphics[width=0.38\linewidth]{\chemin/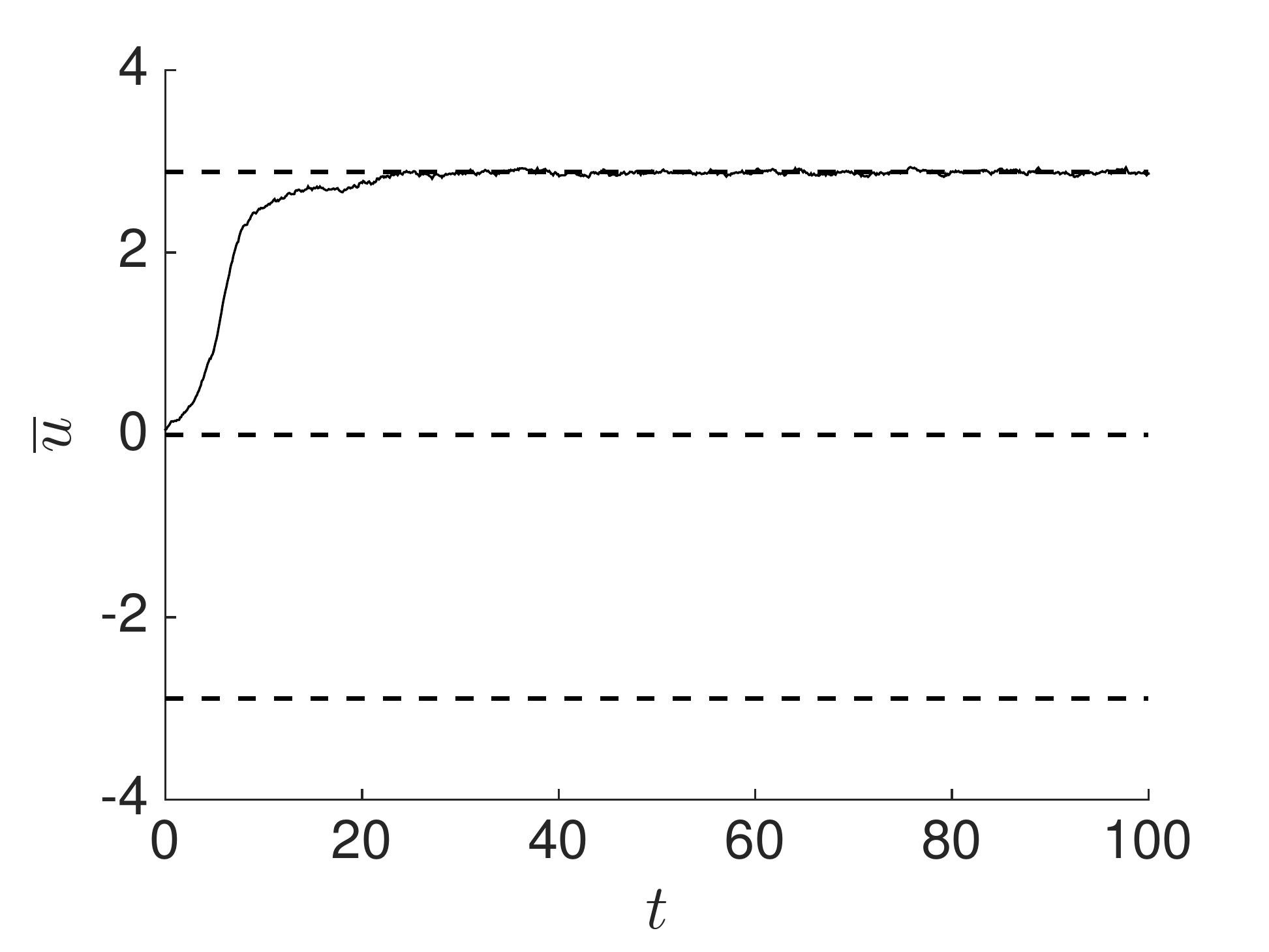}
	{\bf f)}
	\includegraphics[width=0.38\linewidth]{\chemin/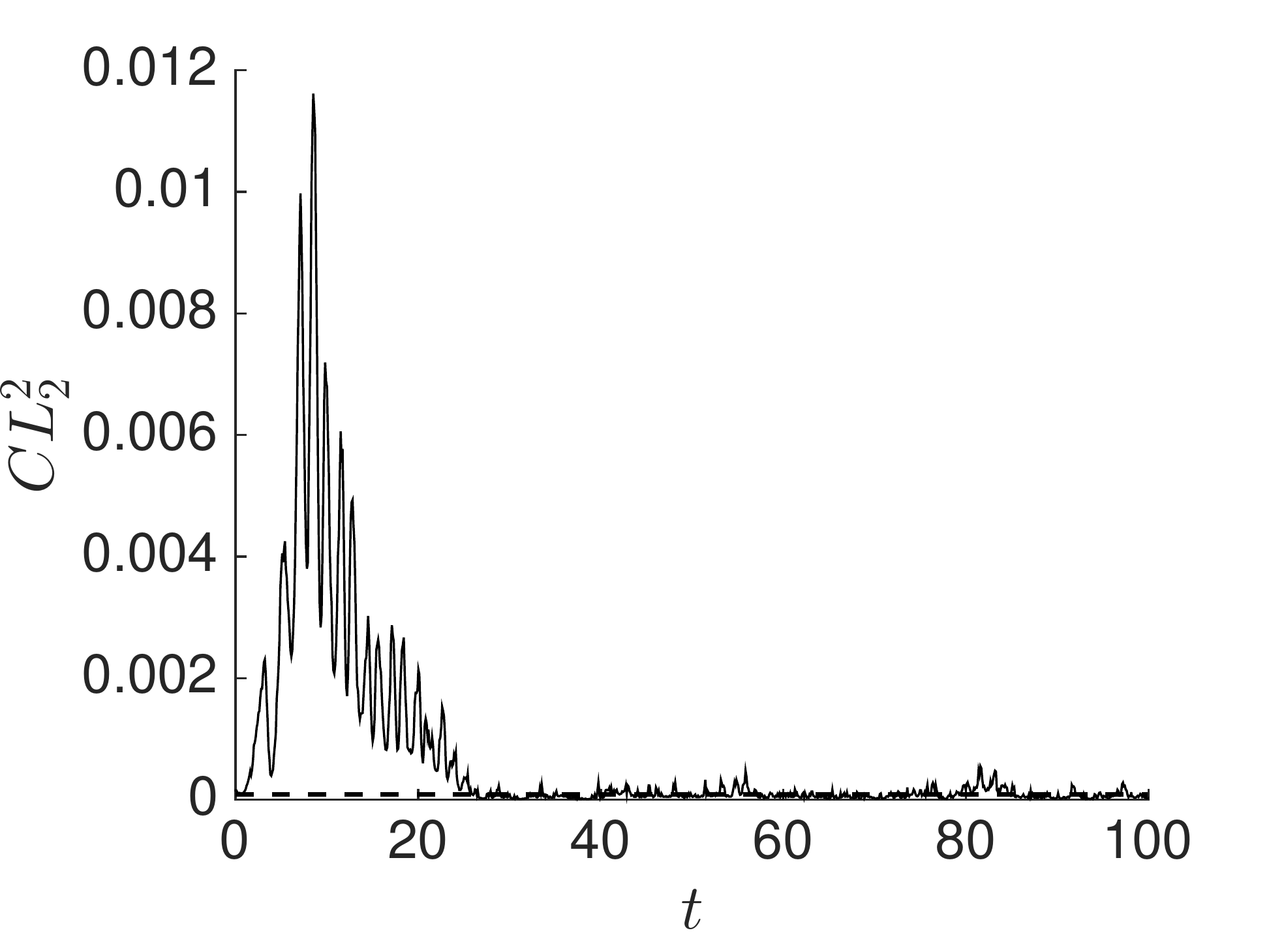}
	
	\caption{The empirical average velocity $\bar{u}^n$  and the  square centered $L^2$-discrepancy $CL_2^2(n)$ 
	for $\sigma=0.5$ (a-b), $\sigma=1$ (c-d), $\sigma=1.5$ (e-f). 
	 The initial positions $\{x_i^0\}_{i=1}^N$ 
are uniformly sampled over $[0,L]$ and the initial velocities $\{u_i^0\}_{i=1}^N$ are sampled from the Gaussian distribution with mean $0$ and variance $\sigma^2/2$.  Here $\Delta t=0.1$, $N=2000$, and $h=6$. 
One can see  that the average velocity is  not $\pm \xi_e$ and the spatial distribution is not uniform when $\sigma=0.5$ or $1$
while the average velocity is $\xi_e$ and the spatial distribution is uniform when $\sigma=1.5$.
}
	\label{fig:4}
\end{figure}

The previous results illustrate the fact that the order states $\rho_{\pm \xi_e}$ are linearly stable 
when $\sigma$ exceeds the threshold value.
This means that, if the initial conditions are close to $\rho_{\pm \xi_e}$, then 
the empirical distribution remains close to it. This is what we see in Figure \ref{fig:2}c-f.
When the initial conditions are far from $\rho_{\pm \xi_e}$ and when the noise level is significantly larger than the threshold value, then 
the empirical distribution quickly converges to one of the two stationary distributions $\rho_{\pm \xi_e}$,
as we can see in Figure \ref{fig:4}e-f:
after a transient period, the empirical distribution in positions becomes uniform and the mean velocity becomes equal to $\pm \xi_e$.
When the initial conditions are far from the 
two stable stationary states and when the noise level is larger than but close to the threshold value, 
then the empirical distribution reaches a state that is not stationary but looks like 
a moving cluster, as shown in Figure \ref{fig:4}a-d.

\begin{figure}
	\centering
	{\bf a)}
	\includegraphics[width=0.38\linewidth]{\chemin/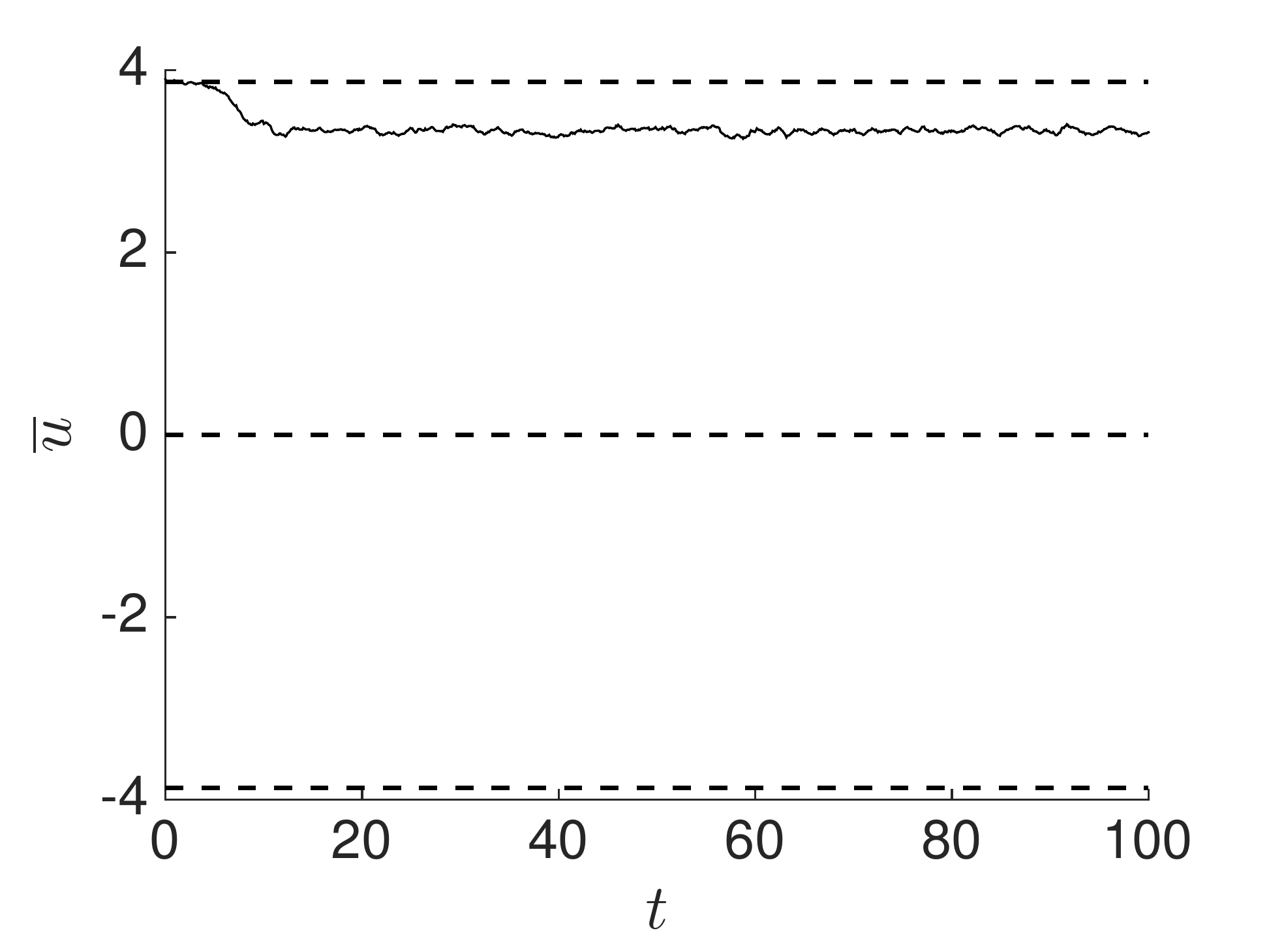}
	{\bf b)}
	\includegraphics[width=0.38\linewidth]{\chemin/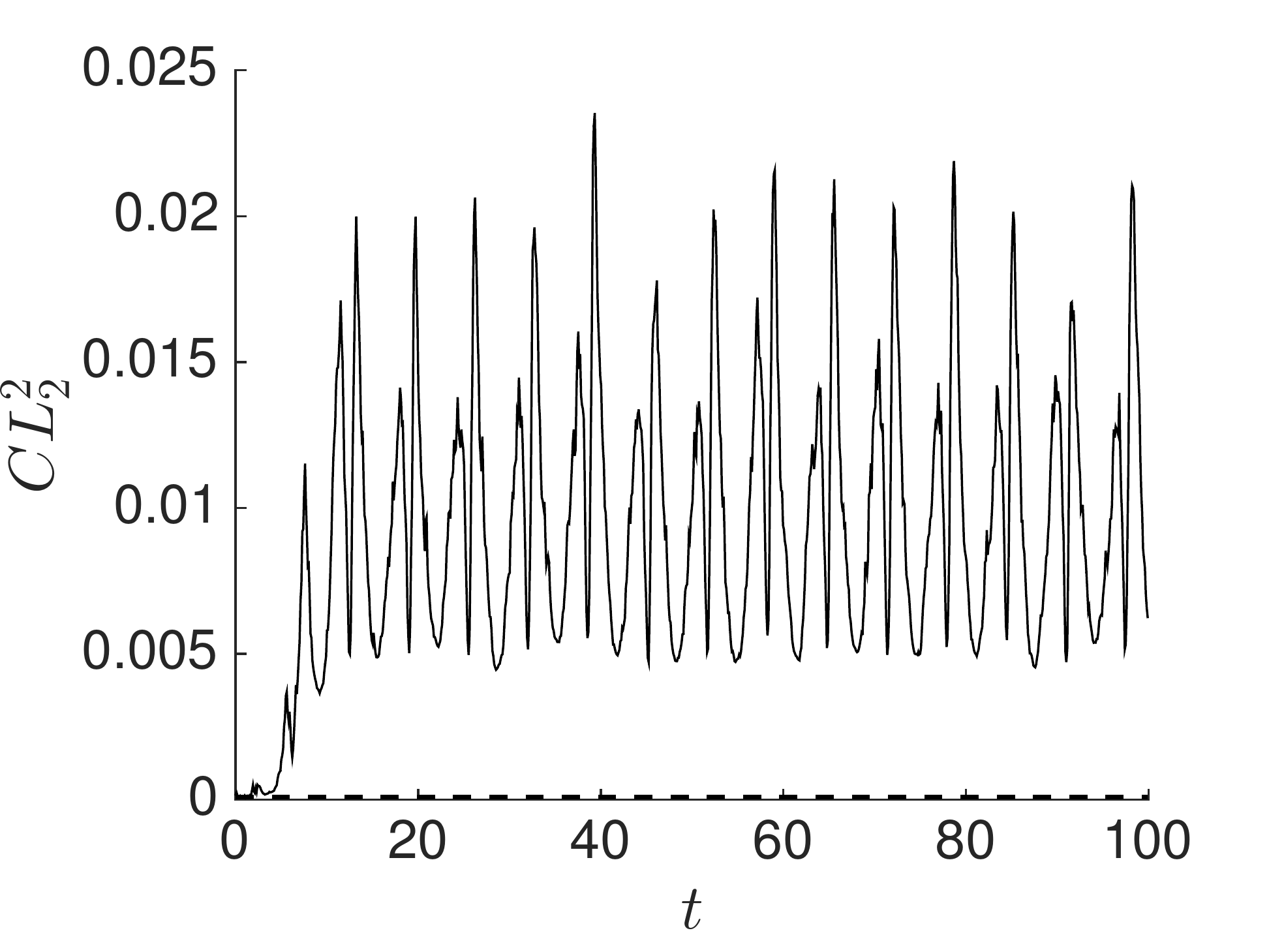}

	{\bf c)}
	\includegraphics[width=0.38\linewidth]{\chemin/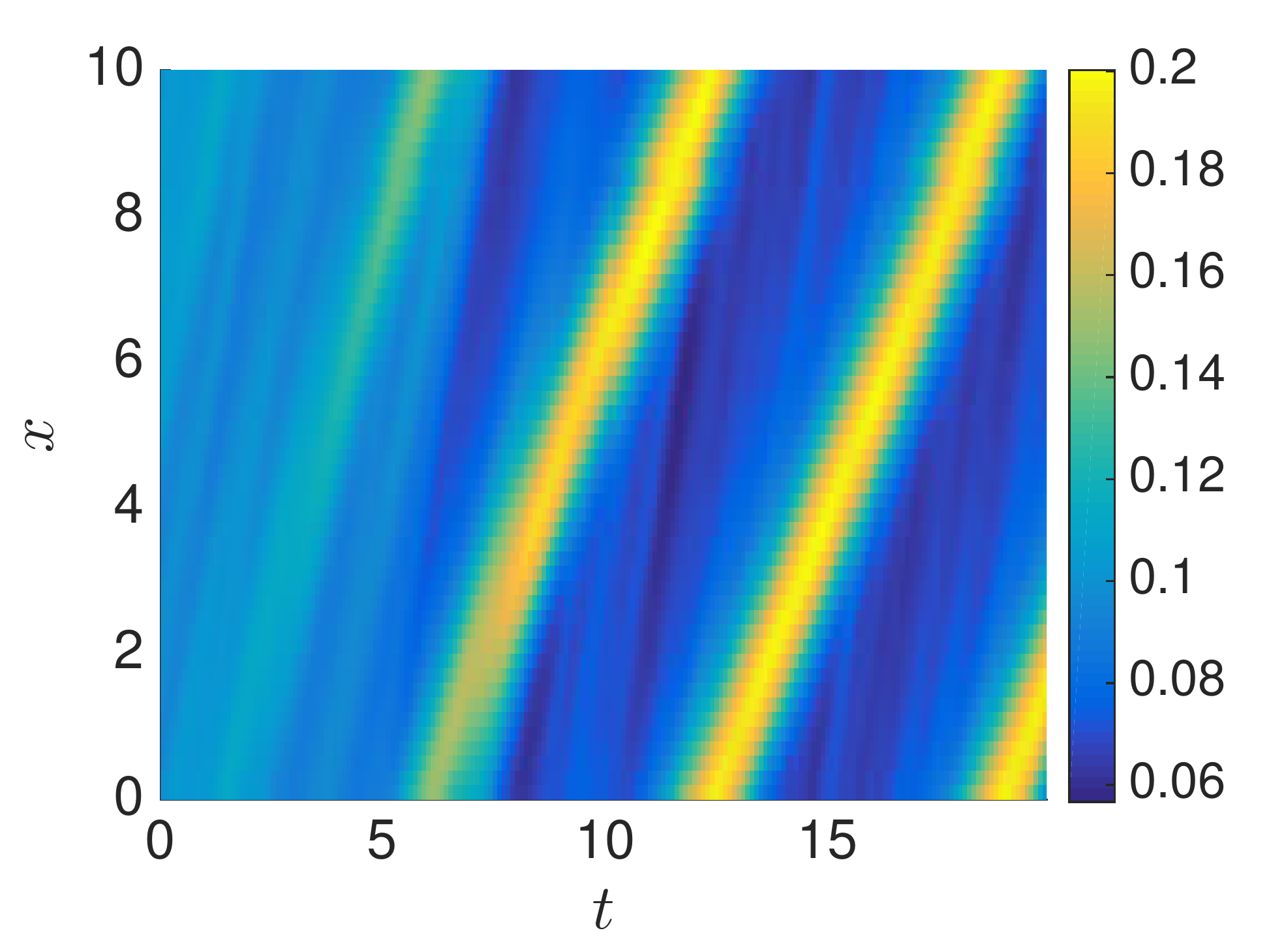}
	{\bf d)}
	\includegraphics[width=0.38\linewidth]{\chemin/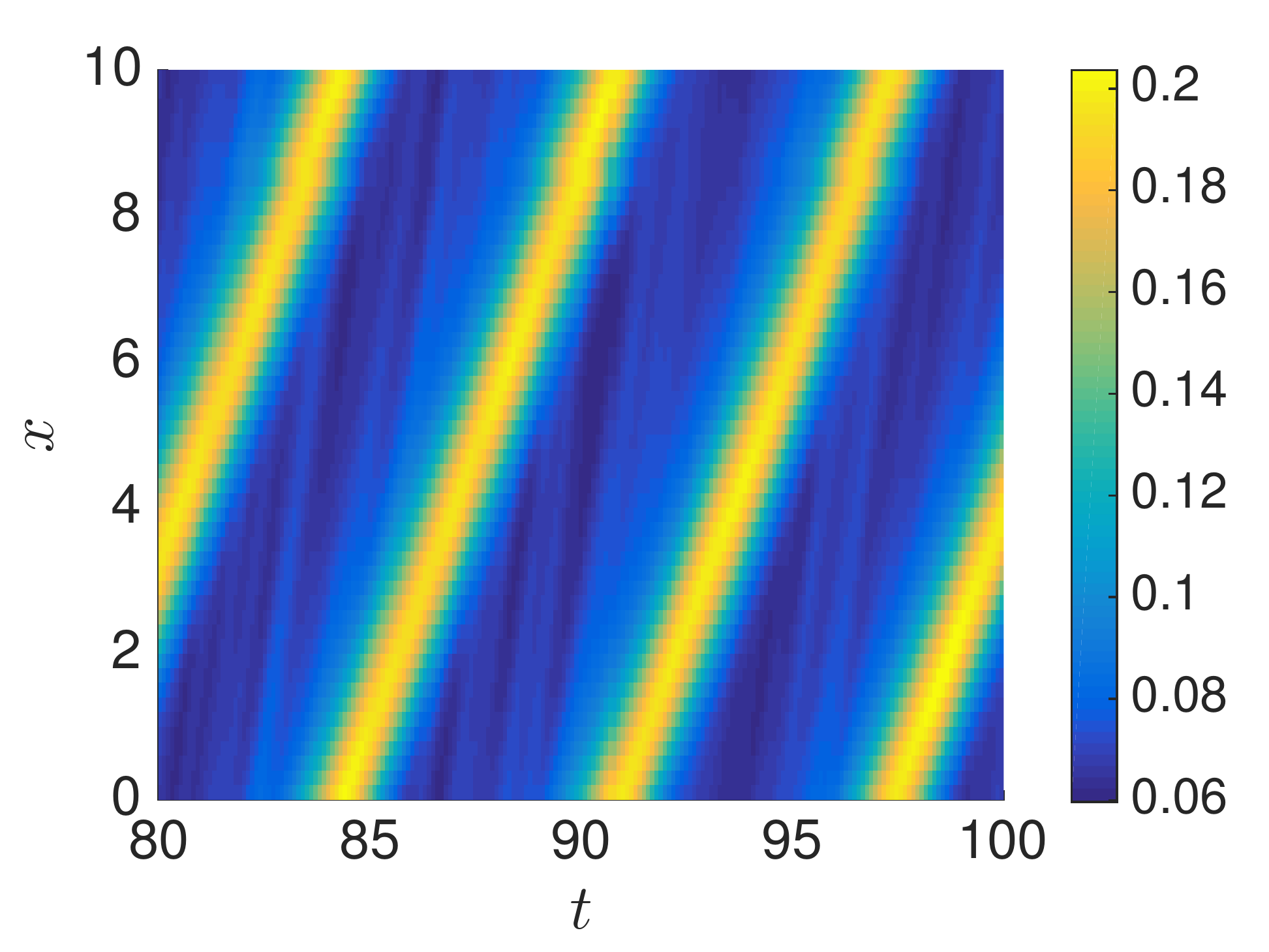}

	\caption{The  empirical average velocity $\bar{u}^n$ (a),
	 the  square centered $L^2$-discrepancy $CL_2^2(n)$ (b), and the empirical position distribution 
	 smoothed by kernel density estimation (c and d).
	The initial positions $\{x_i^0\}_{i=1}^N$  are uniformly sampled over $[0,L]$ and the initial velocities $\{u_i^0\}_{i=1}^N$ are sampled from the Gaussian distribution with mean $\xi_e$ and variance $\sigma^2/2$.  
	Here $\Delta t=0.1$, $N=2000$, $\sigma=1$, and $h=10$. 
	One can see  that the average velocity is  not $\xi_e$ and the spatial distribution is not uniform.
}
	\label{fig:5}
\end{figure}

\begin{figure}
	\centering
	{\bf a)}
	\includegraphics[width=0.38\linewidth]{\chemin/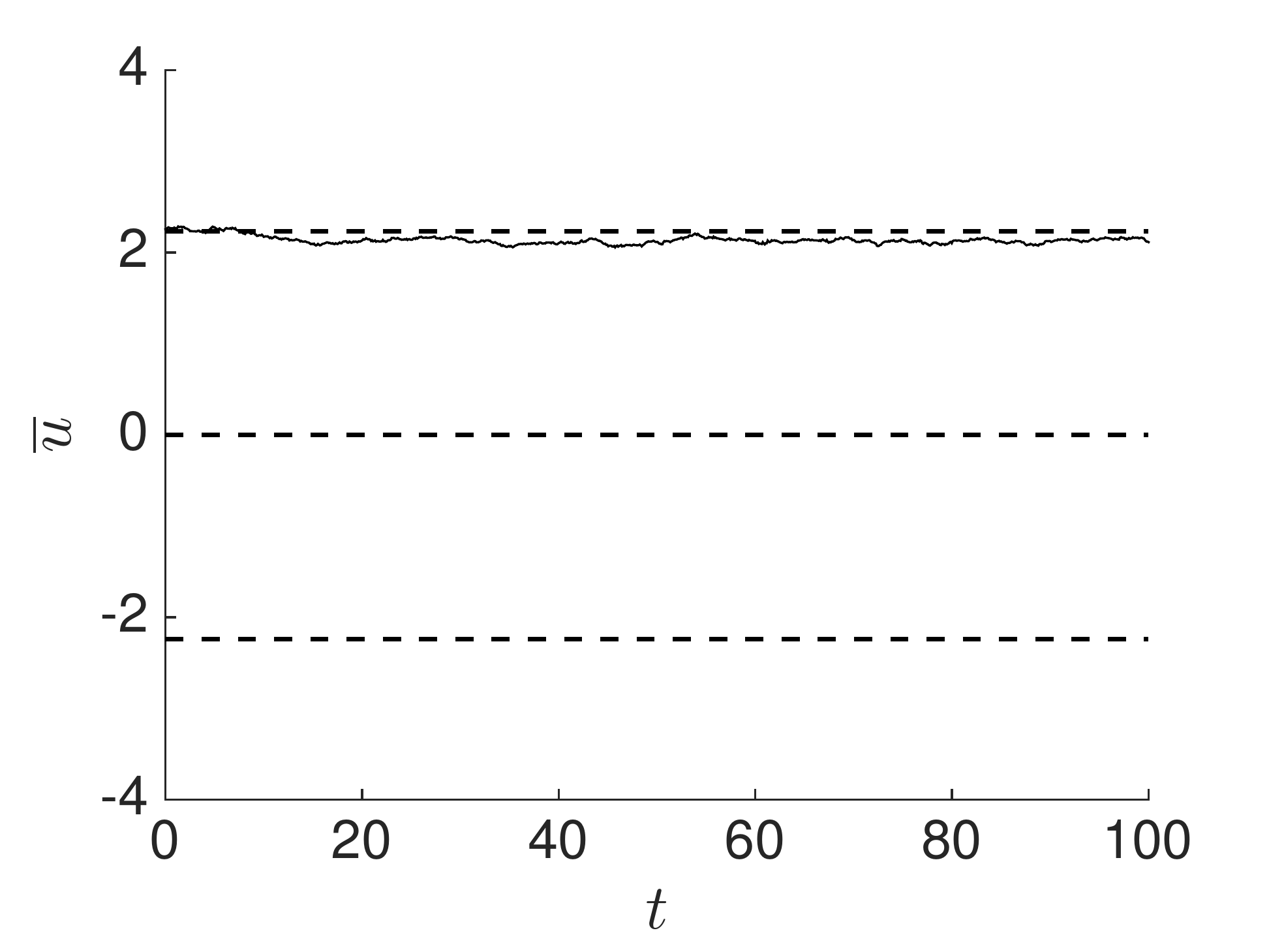}
	{\bf b)}
	\includegraphics[width=0.38\linewidth]{\chemin/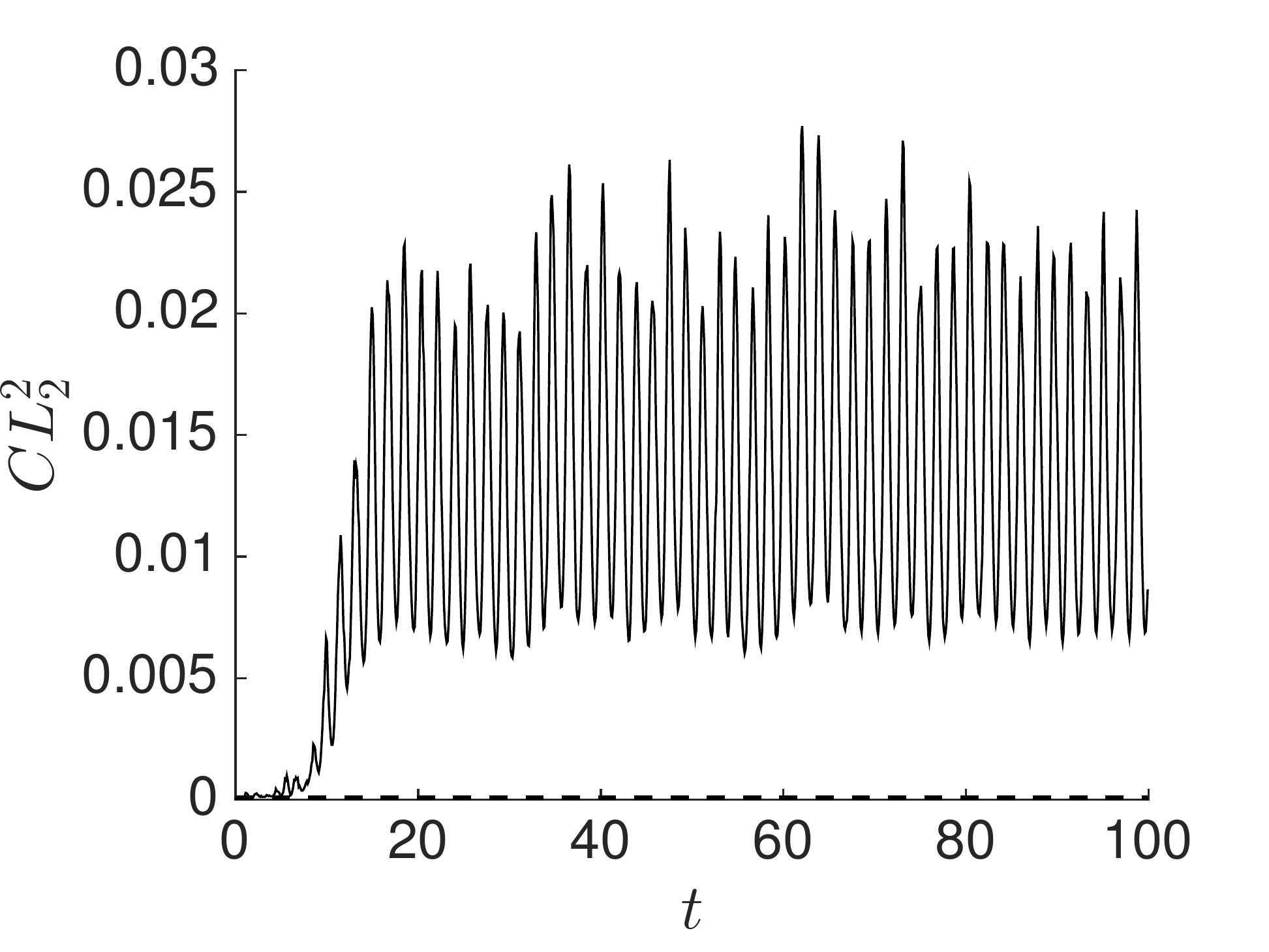}

	{\bf c)}
	\includegraphics[width=0.38\linewidth]{\chemin/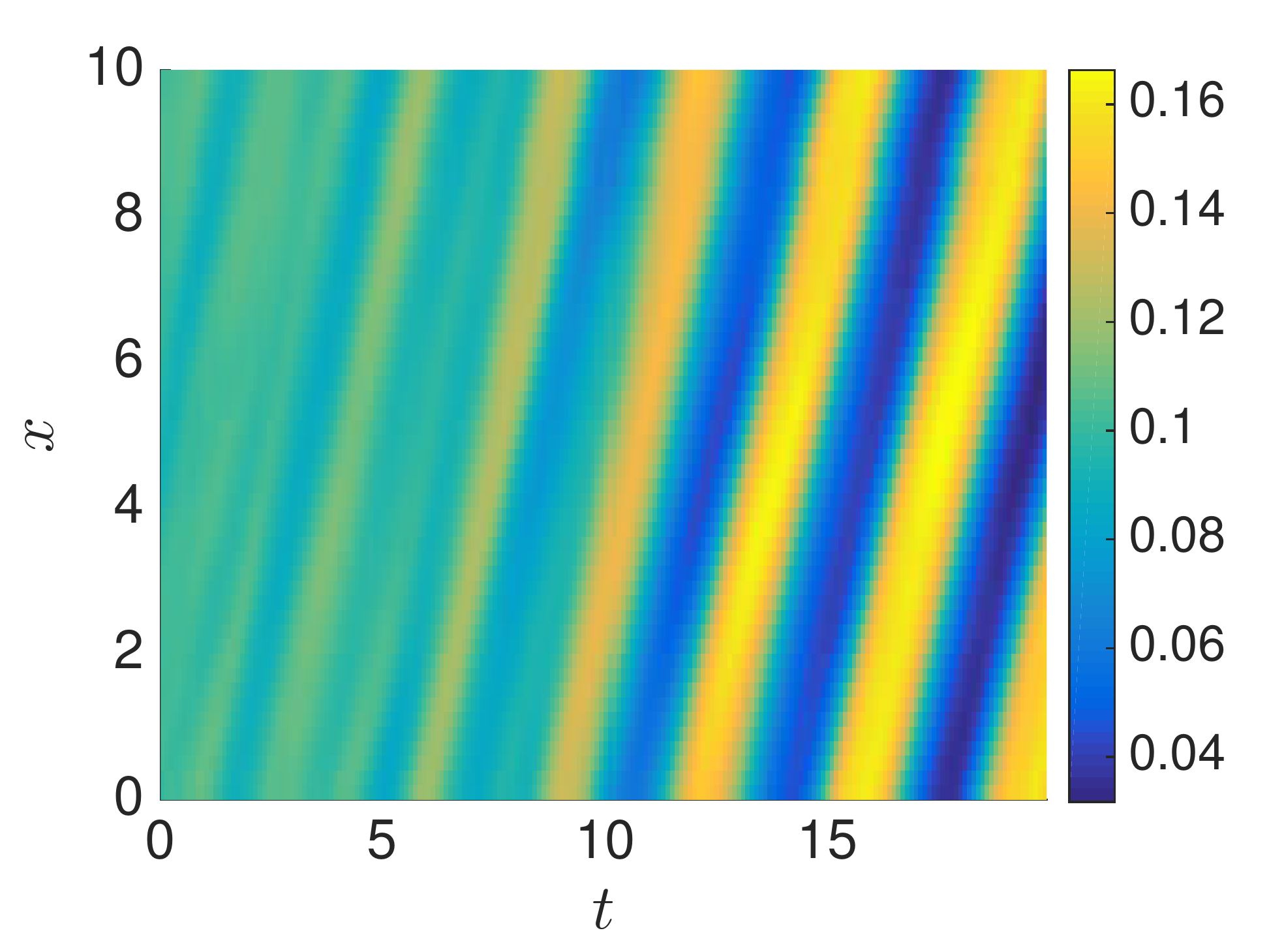}
	{\bf d)}
	\includegraphics[width=0.38\linewidth]{\chemin/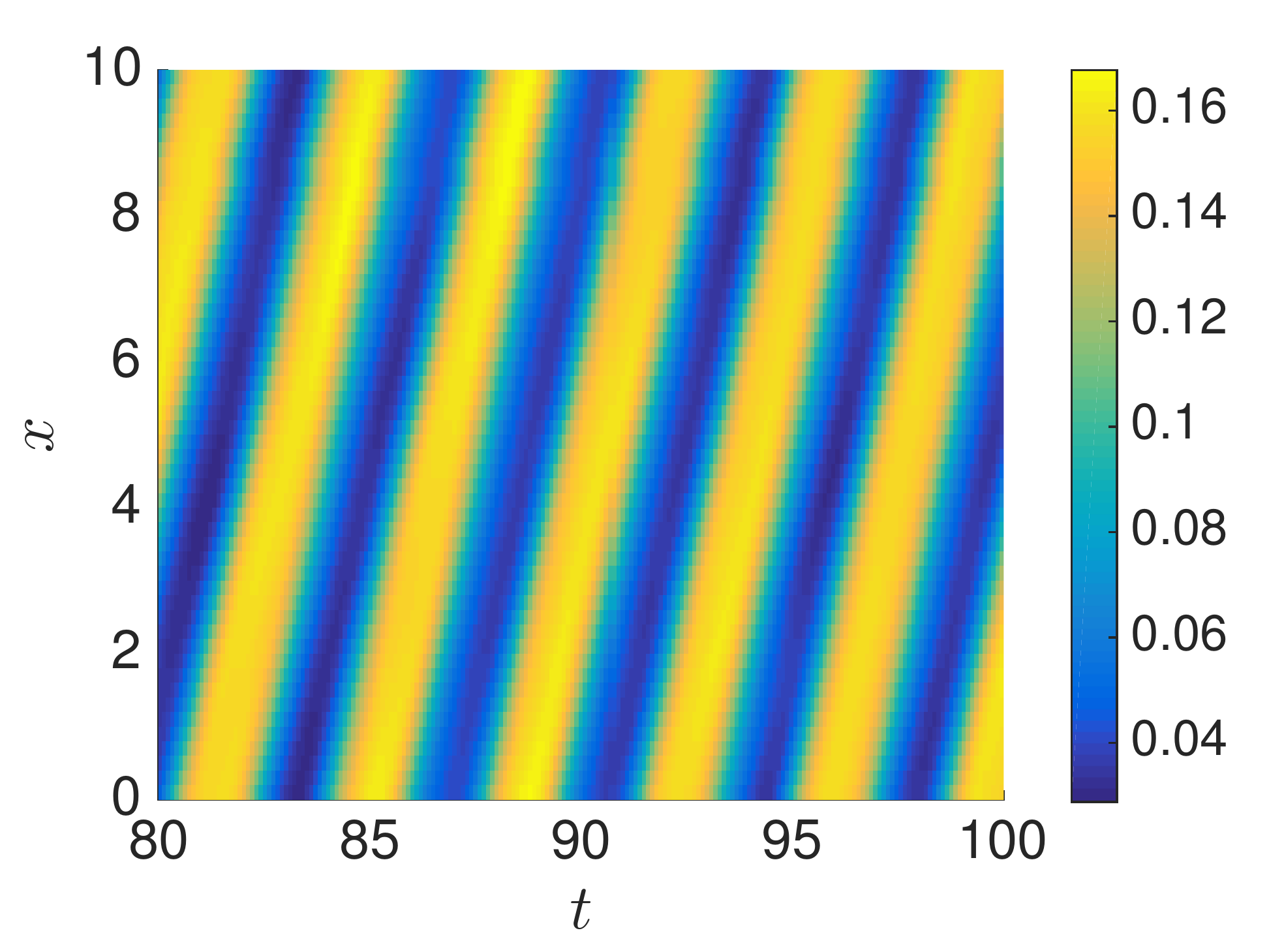}

	\caption{The  empirical average velocity $\bar{u}^n$ (a),
	 the  square centered $L^2$-discrepancy $CL_2^2(n)$ (b), and the empirical position distribution 
	 smoothed by kernel density estimation (c and d).
	The initial positions $\{x_i^0\}_{i=1}^N$  are uniformly sampled over $[0,L]$ and the initial velocities $\{u_i^0\}_{i=1}^N$ are sampled from the Gaussian distribution with mean $\xi_e$ and variance $\sigma^2/2$.  
	Here $\Delta t=0.1$, $N=2000$, $\sigma=1$, and $h=5$. 
		One can see  that the average velocity is  not $\xi_e$ and the spatial distribution is not uniform.
}
	\label{fig:6}
\end{figure}

As revealed by the linear stability analysis carried out in subsection \ref{subsec:instaspa},
the results are not monotoneous as a function of $h$.
Figure \ref{fig:5} and \ref{fig:6} show that the stationary states $\rho_{\pm \xi_e}$ are unstable when $\sigma=1.0$ 
for $h=5$ and for $h=10$ while they are stable for $h=6$ (see Figure \ref{fig:2}c-d),
as predicted by the theory in subsection~\ref{subsec:instaspa}.
This means that the two order states are stable when the depths of the wells of the double-well potential
are neither too large nor too small.

\subsection{Large deviations}
\begin{figure}
	\centering	
	{\bf a)}
	\includegraphics[width=0.38\linewidth]{\chemin/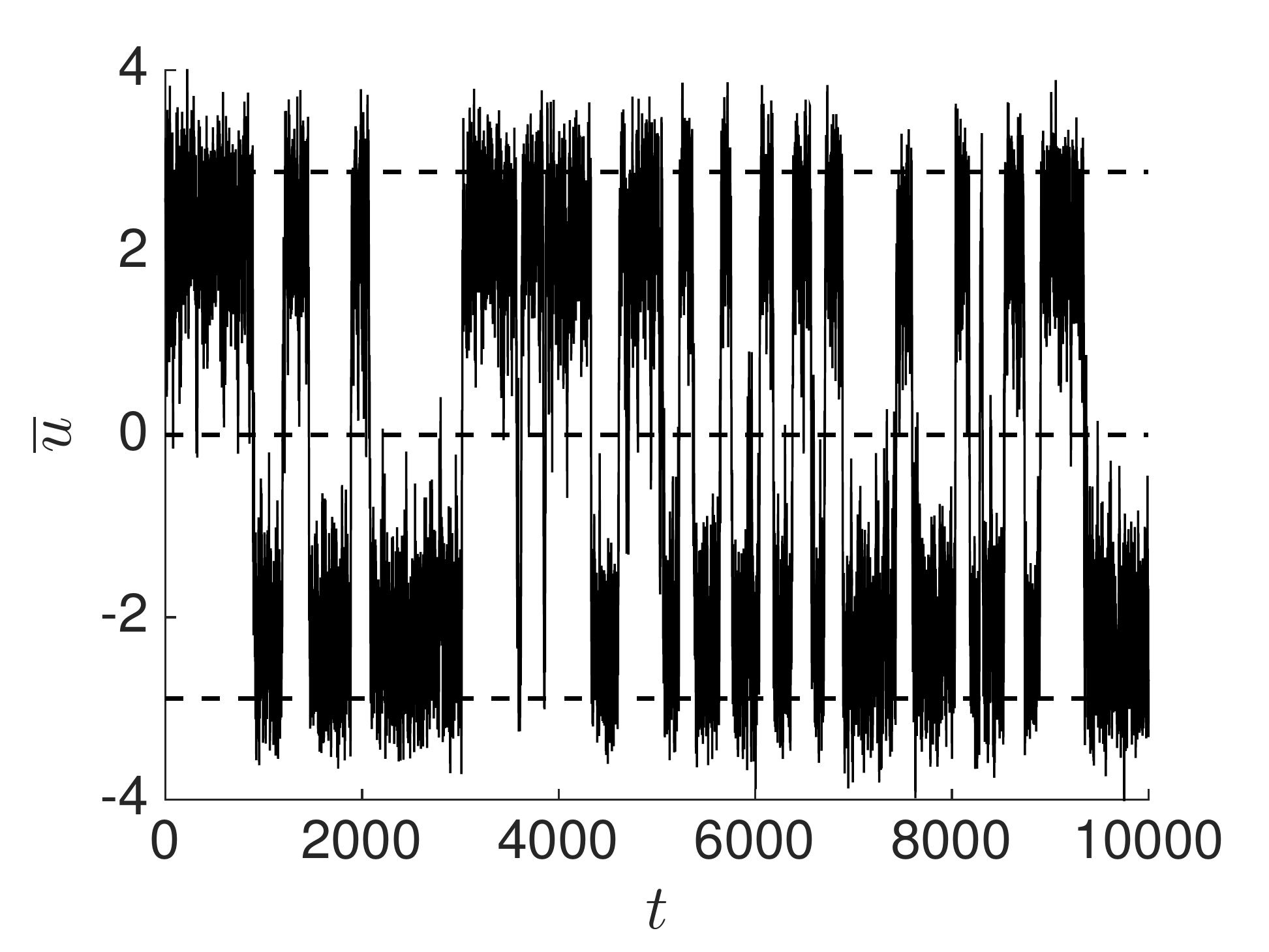}
	{\bf b)}
	\includegraphics[width=0.38\linewidth]{\chemin/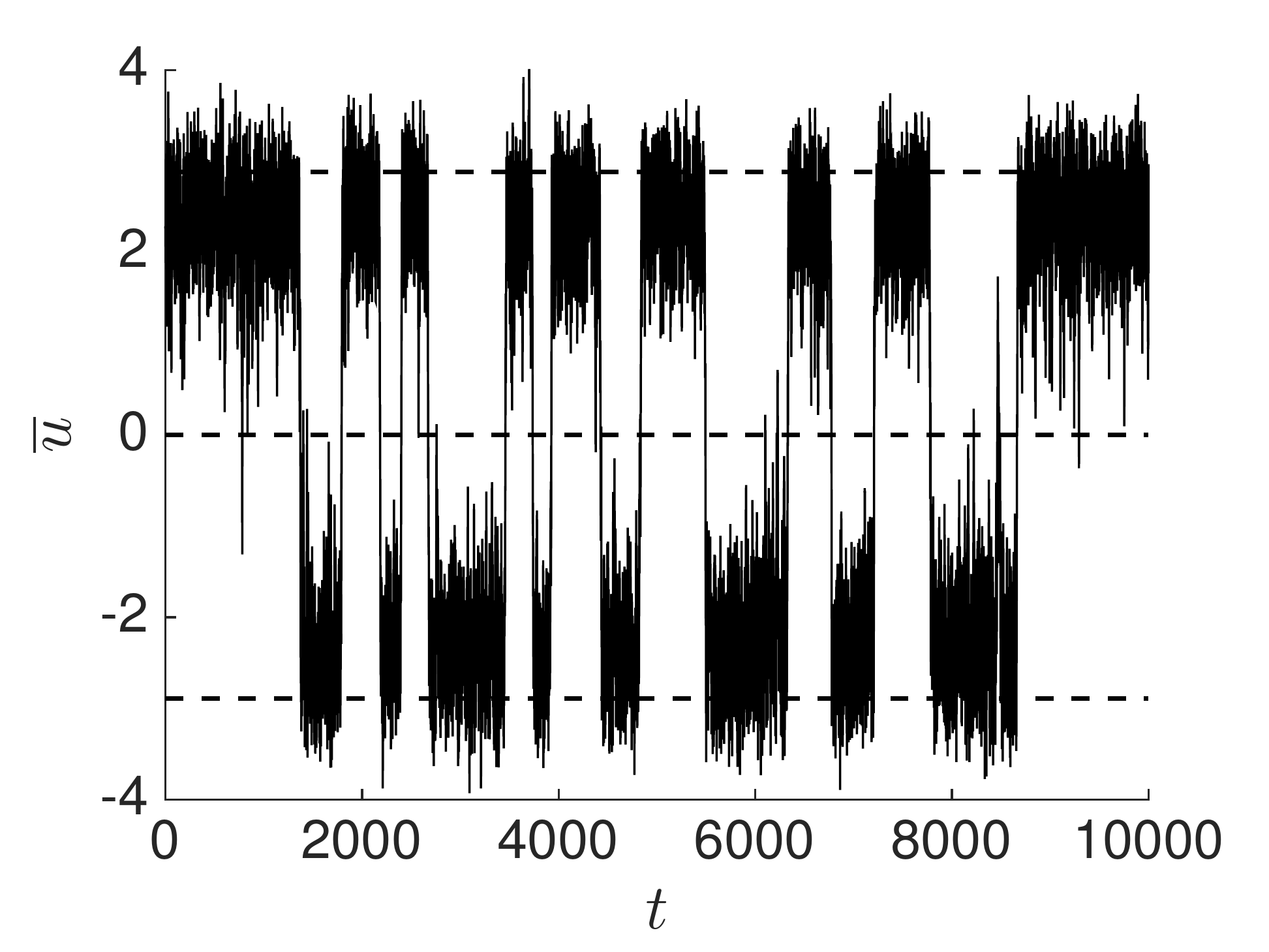}
	
	{\bf c)}
	\includegraphics[width=0.38\linewidth]{\chemin/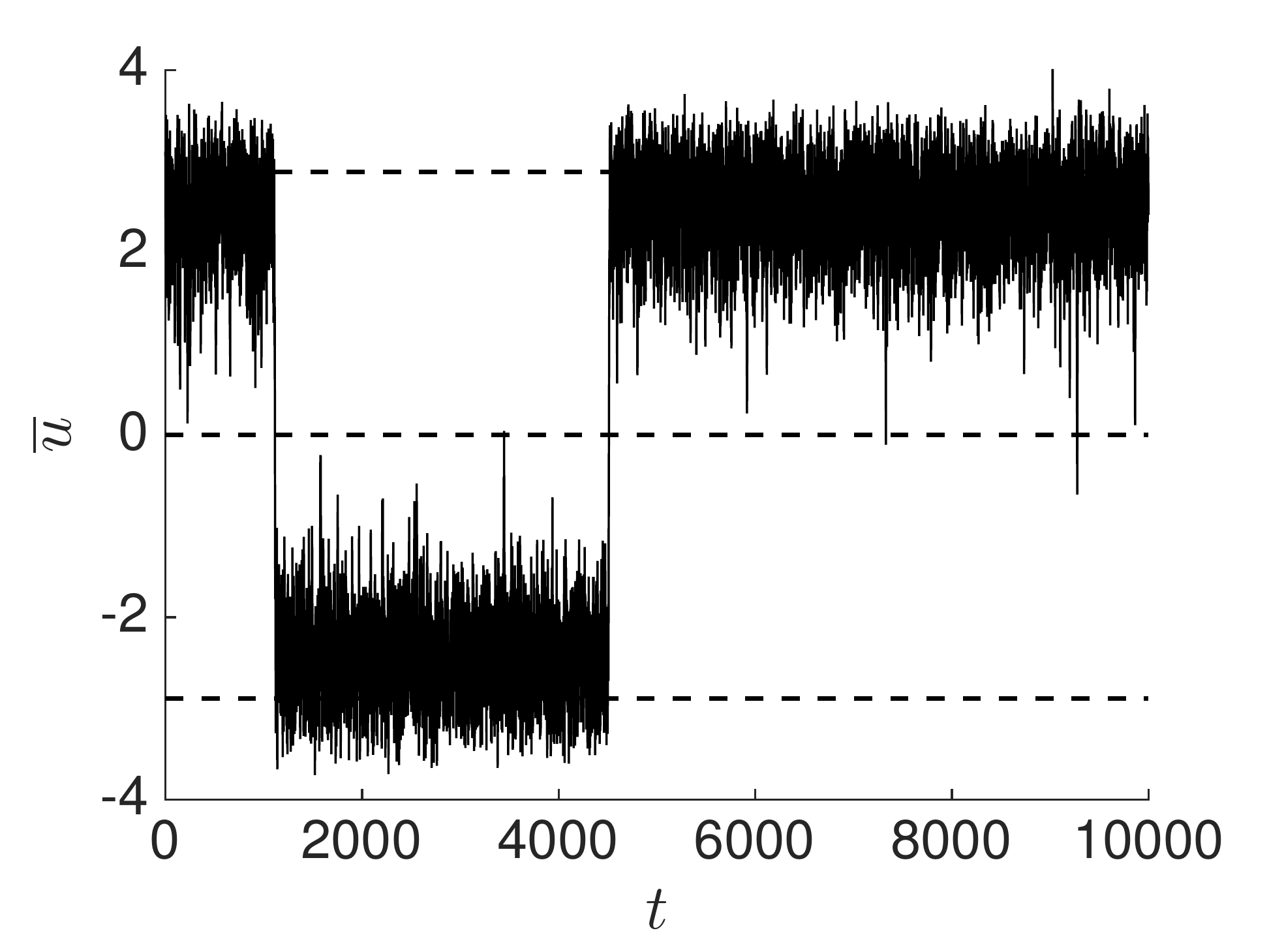}
	{\bf d)}
	\includegraphics[width=0.38\linewidth]{\chemin/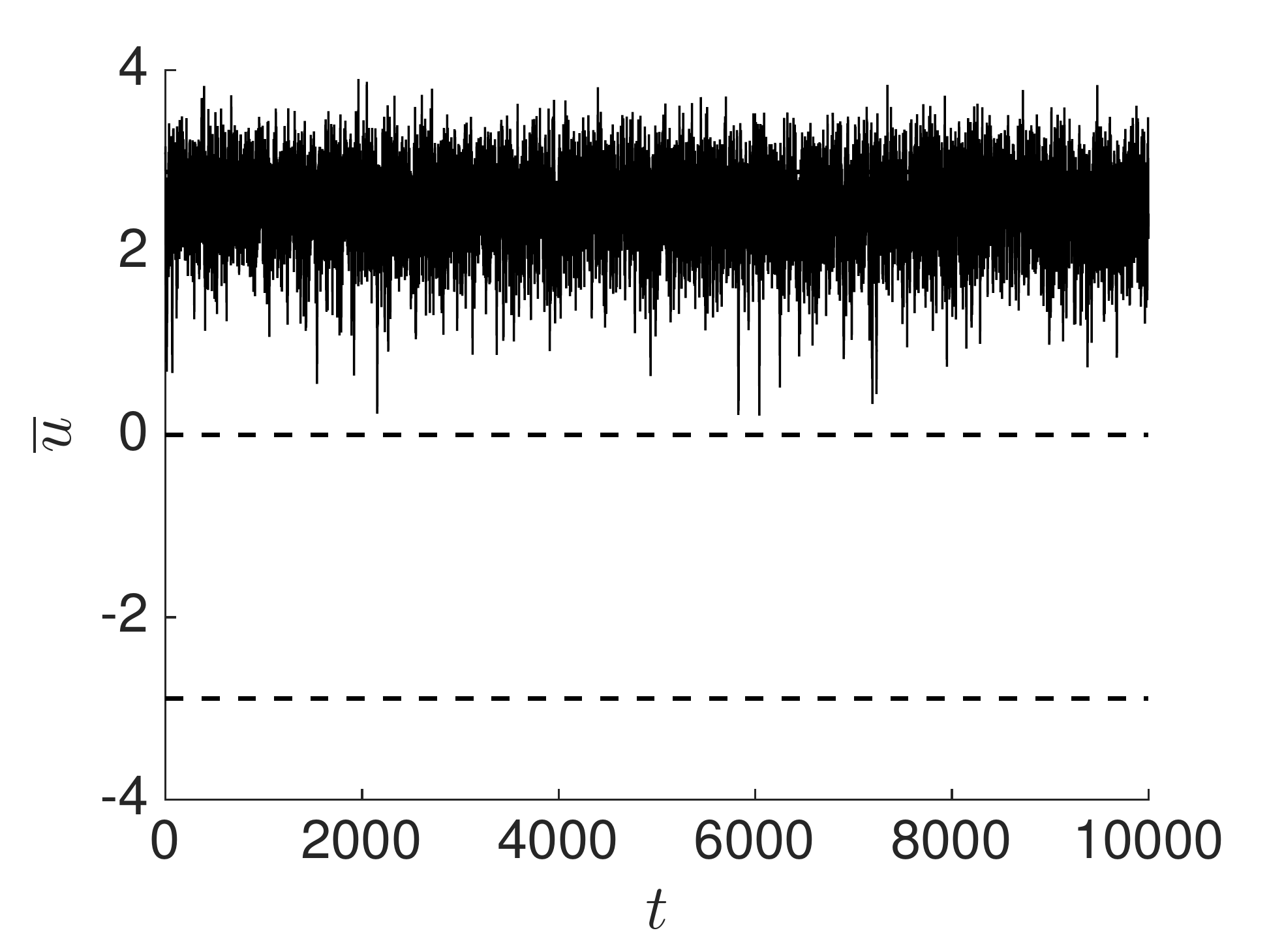}
	\caption{The empirical average velocity $\bar{u}^n$ at each time step $t_n$ for $N=80$ (a), $N=100$ (b), $N=120$ (c), and $N=140$ (d).  
		Here $\Delta t=0.1$, $h=6$ and $\sigma=5$. 
		 The frequencies of the transitions between the two stable order states decay when $N$ increases. 
The system has less transitions with a higher number of agents.}
	\label{fig:LD for N}
\end{figure}

\begin{figure}
	\centering	
	{\bf a)}
	\includegraphics[width=0.38\linewidth]{\chemin/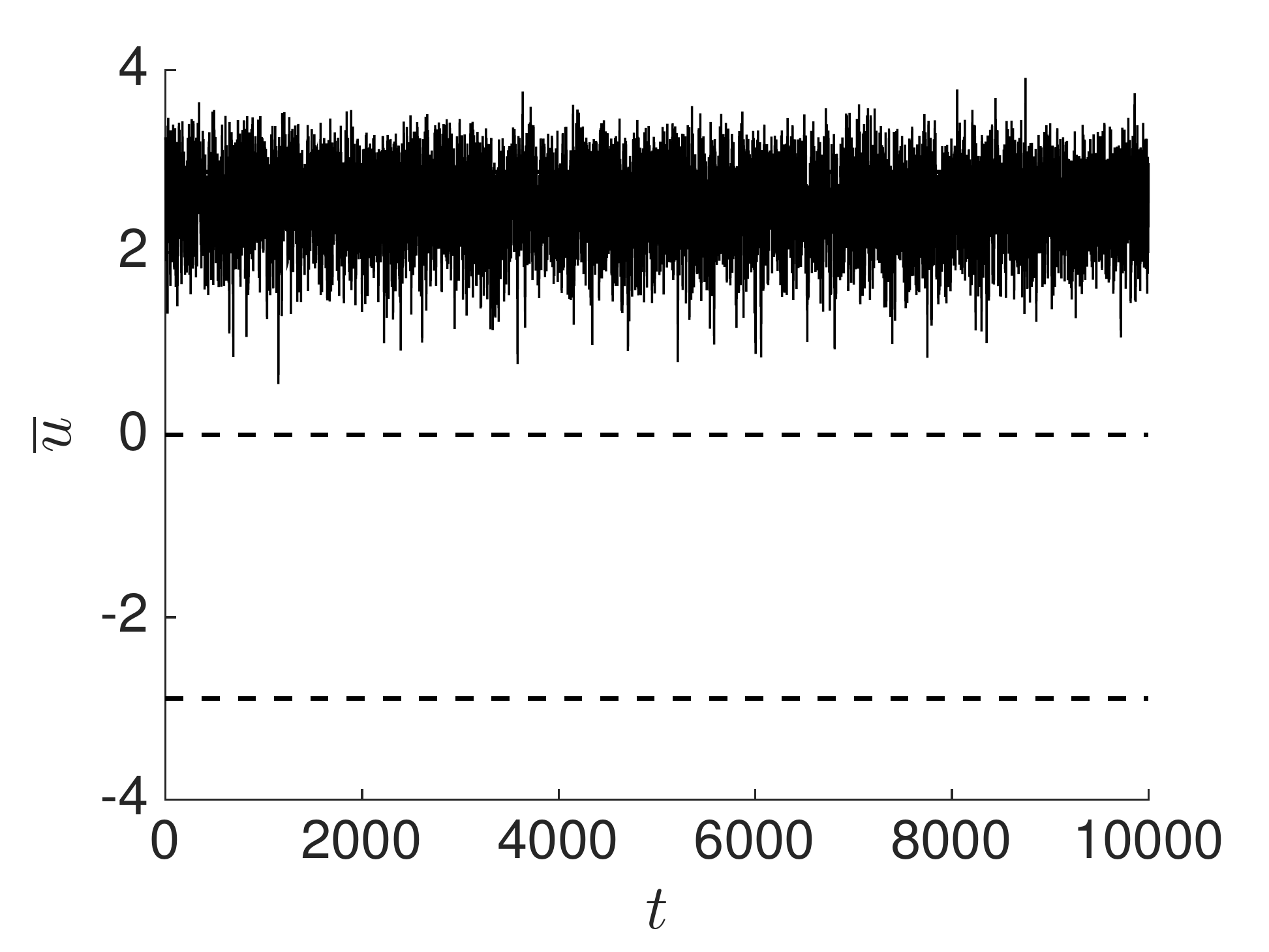}
	{\bf b)}
	\includegraphics[width=0.38\linewidth]{\chemin/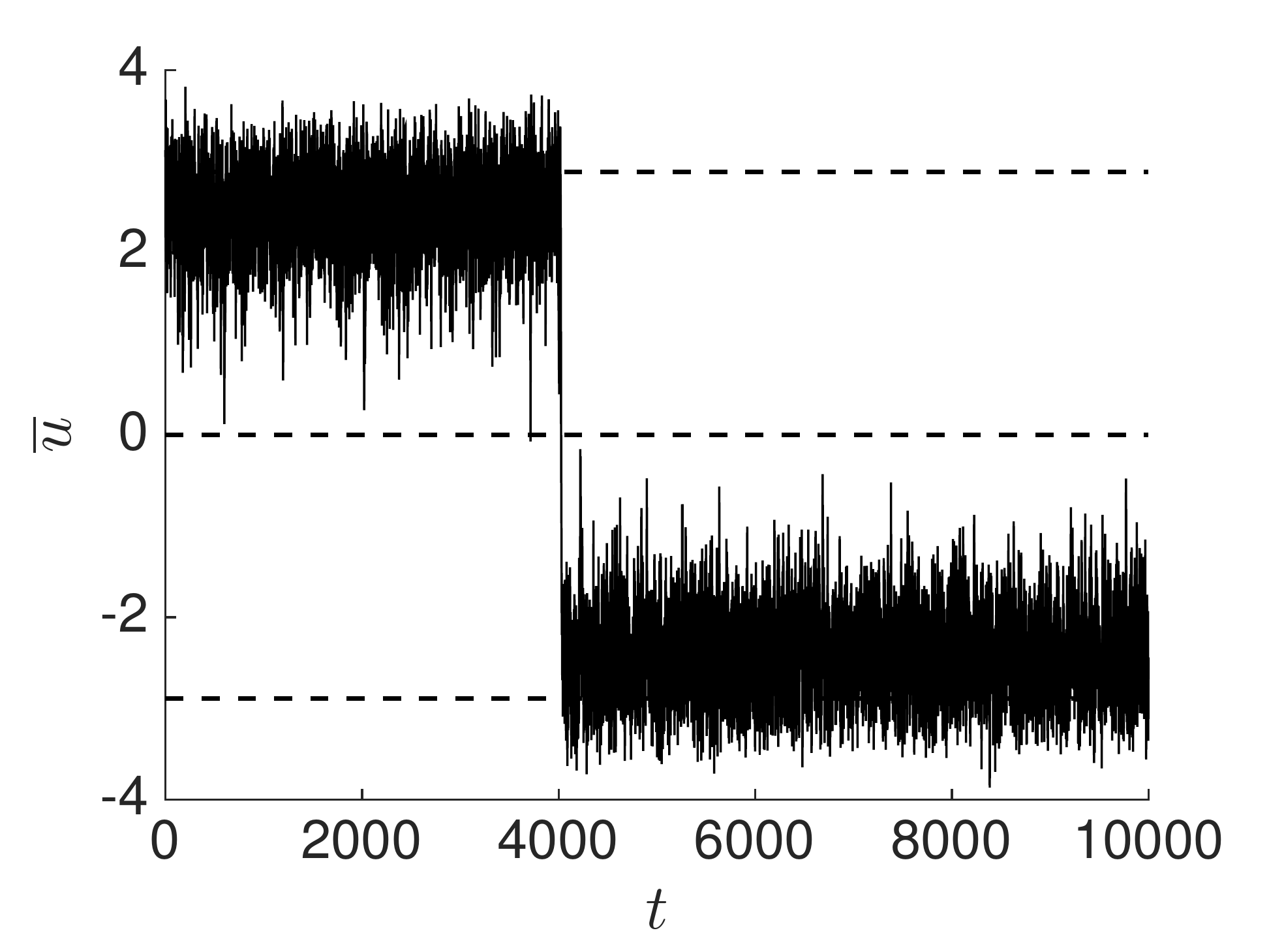}
	
	{\bf c)}
	\includegraphics[width=0.38\linewidth]{\chemin/fig_ld_n100_sigma5_h6.pdf}
	{\bf d)}
	\includegraphics[width=0.38\linewidth]{\chemin/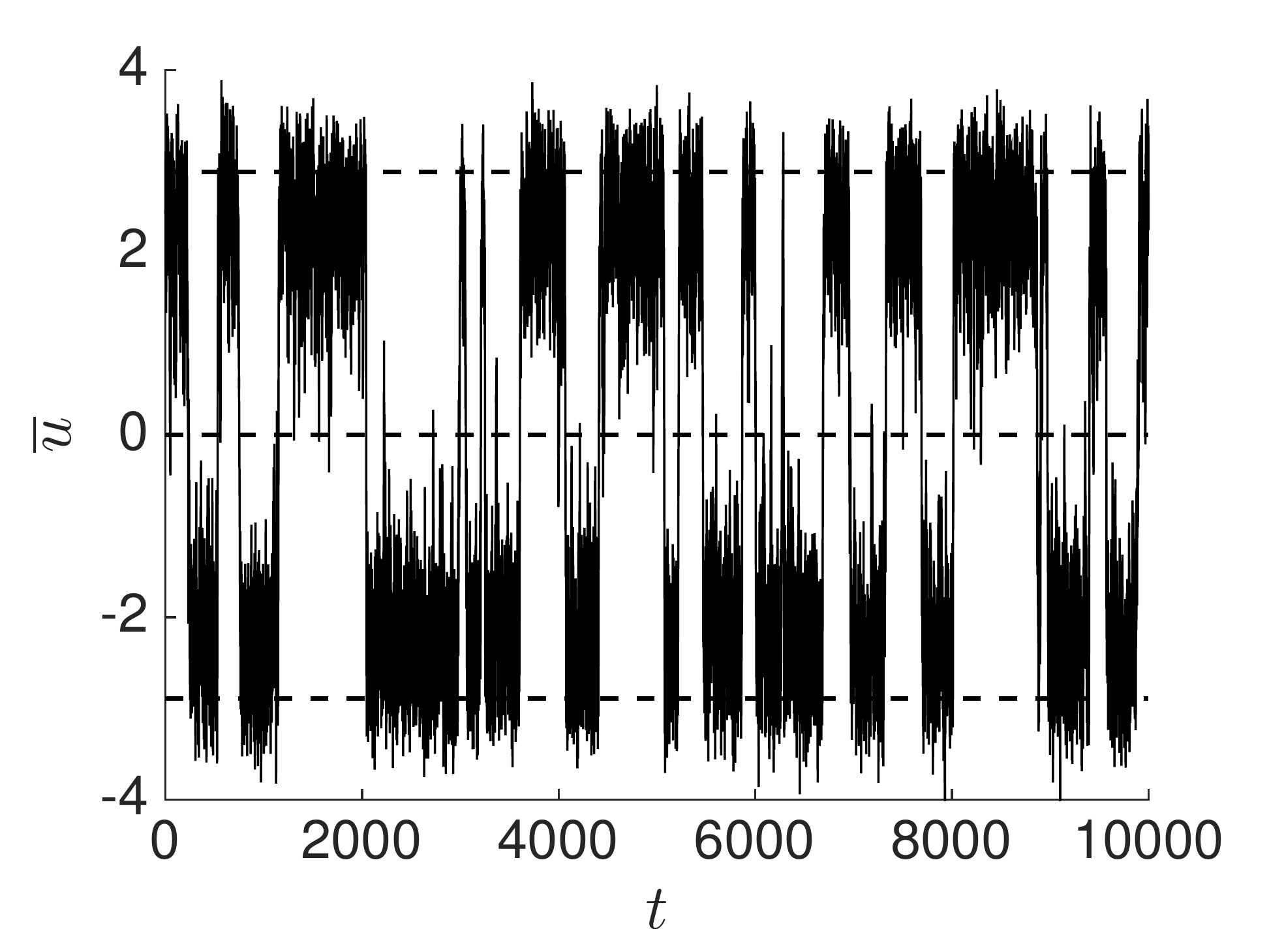}
	\caption{The empirical average velocity $\bar{u}^n$ at each time step $t_n$ for $\sigma=4$ (a), $\sigma=4.5$ (b), $\sigma=5$ (c), and $\sigma=5.5$ (d). 
		Here $\Delta t=0.1$, $N=100$, and $h=6$.  The frequencies of the transitions between the two stable order states increases when $\sigma$ increases. The system has more transitions with a higher $\sigma$.}
	\label{fig:LD for sigma}
\end{figure}

\begin{figure}
	\centering	
	{\bf a)}
	\includegraphics[width=0.38\linewidth]{\chemin/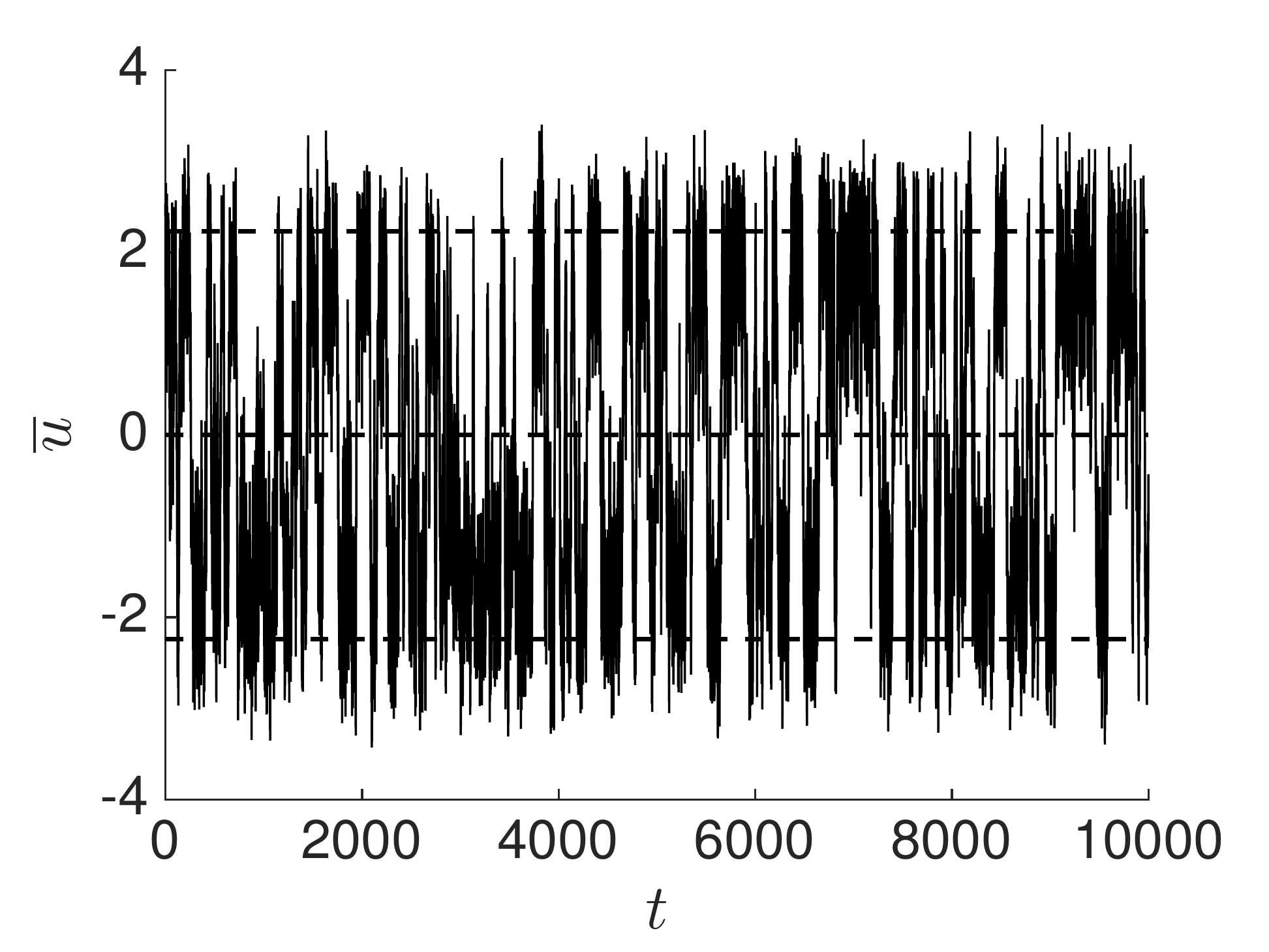}
	{\bf b)}
	\includegraphics[width=0.38\linewidth]{\chemin/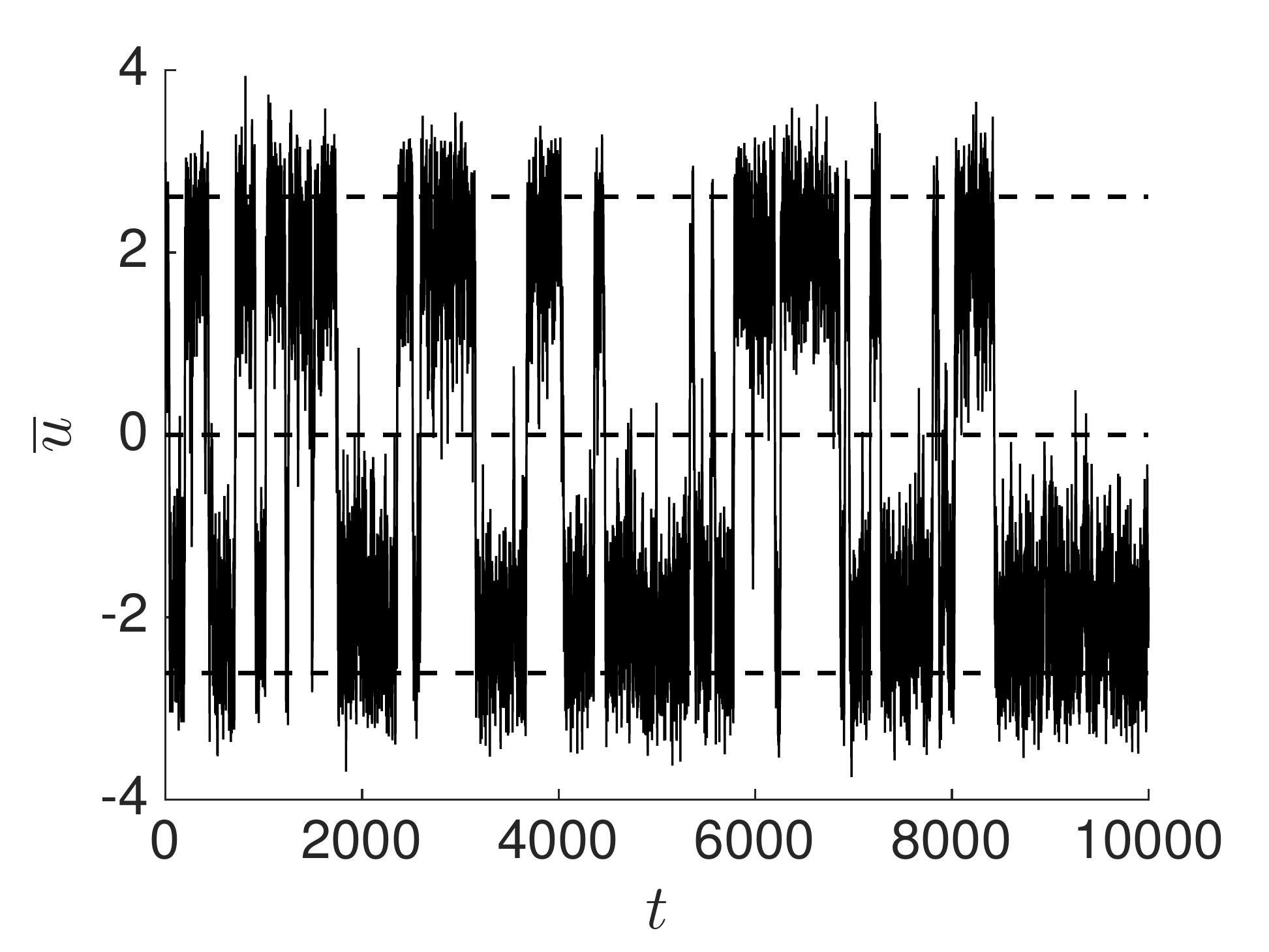}\\
	
	{\bf c)}
	\includegraphics[width=0.38\linewidth]{\chemin/fig_ld_n100_sigma5_h6.pdf}
	{\bf d)}
	\includegraphics[width=0.38\linewidth]{\chemin/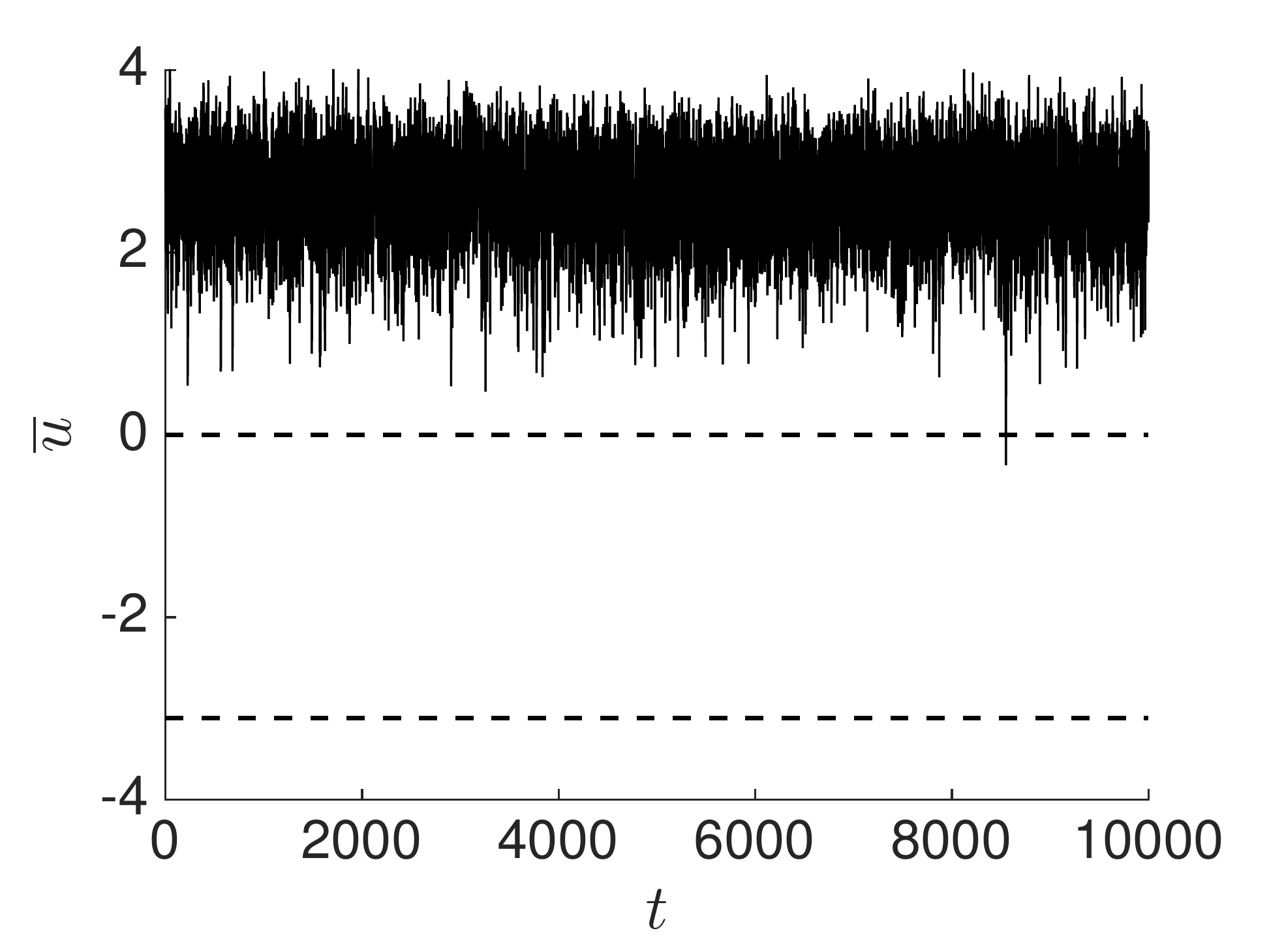}
	\caption{The empirical average velocity $\bar{u}^n$ at each time step $t_n$ for $h=5$ (a), $h=5.5$ (b), $h=6$ (c), and $h=6.5$ (d). 
	Here $\Delta t=0.1$, $N=100$,  and $\sigma=5$. The frequencies of the transitions between the two stable order states  decays with $h$.}
	\label{fig:LD for h}
\end{figure}

In this subsection we assume that the conditions are fulfilled  so that $\rho_{\xi_e}$ and $\rho_{-\xi_e}$ exist and are stable.
We consider the transition probability from one stable state to the other:
\begin{equation*}
	\mathbf{P}(\mu_N \in A) \overset{N\gg 1}{\approx} \exp\left(-N\inf_{\mu\in A} I(\mu)\right),
\end{equation*}
where the rare event $A$ is the set (\ref{eq:event of transitions}) of all possible transitions from $\rho_{\xi_e} $ at $t=0$ to $\rho_{-\xi_e} $ at $t=T$,
and the rate function $I$ is defined in (\ref{eq:rate function}).

In Figures \ref{fig:LD for N}, \ref{fig:LD for sigma}, and \ref{fig:LD for h}, we qualitatively examine the effects of $N$, $\sigma$, and $h$ to the transition 
probability $\mathbf{P}(\mu_N \in A)$ by checking the frequencies of transitions of the empirical average velocity 
$\bar{u}^n $ between the two stable order states 
$\rho_{\pm \xi_e}$.  In Figure \ref{fig:LD for N}, we can observe that the system has less transitions as $N$ becomes larger.  In other words, the probability of transition 
$\mathbf{P}(\mu_N \in A)$ becomes smaller as $N$ becomes larger; this is consistent with the fact that the large deviation principle predicts that 
$\mathbf{P}(\mu_N \in A)$ is an exponential decay function of $N$.  In Figure \ref{fig:LD for sigma}, we can observe that the system experiences more transitions as 
$\sigma$ becomes larger.  Therefore, $\mathbf{P}(\mu_N \in A)$ becomes larger as $\sigma$ becomes larger and this is confirmed by the fact that the rate function is 
approximately proportional to $1/\sigma^2$ (note, however, that $\sigma$ is also in the operator ${\mathcal L}^*$).  
In Figure \ref{fig:LD for h}, we find that the probability of transition decreases as $h$ increases. 
This is qualitatively consistent with the fact that $h$ determines the height of the potential barrier between the two potential wells,
in analogy with the classical problem of diffusion exit from a domain \cite[Section 5.7]{Dembo1998}.

\section{The Nonsymmetric Case}
\label{sec:nonsymmetry}

Let us briefly revisit the previous analysis when the average velocity has the form \cite{Tadmor2011}
\begin{equation}
	\left<u\right>_i = 
	\begin{cases}
		\frac{1}{N_i} \sum_{j=1}^{N} u_j\phi(\|x_j-x_i\|), &\text{if $N_i>0$,}\\
		0, &\text{if $N_i=0$,}
	\end{cases}
\end{equation}
where $N_i$ is the weighted number of agents in the neighborhood of the $i$th agent:
\begin{equation}
	N_i = \sum_{j=1}^N \phi(\|x_j-x_i\|).
\end{equation}
If $\mu_N(0,dx,du)$ converges to a deterministic measure  $\bar{\rho}(x,u)dxdu$, then $\mu_N(t,dx,du) $ converges to the deterministic measure
$\rho(t,x,u)dxdu$ whose density is the solution of the nonlinear Fokker-Planck equation 
\begin{equation}
	\label{eq:nonlinear Fokker-Planck ns}
	\frac{\partial \rho}{\partial t} = 
	- u\frac{\partial \rho}{\partial x}
	- \frac{\partial}{\partial u}\left[\left(G\left(\frac{\iint u'\phi(\|x'\|) \rho(t,x-x',u') du'dx'}
	{\iint \phi(\|x'\|) \rho(t,x-x',u') du'dx'} \right) - u \right) \rho\right]
	+ \frac{1}{2} \sigma^2 \frac{\partial^2 \rho}{\partial u^2} ,
\end{equation}
starting from $\rho(t=0,x,u)=\bar{\rho}(x,u)$.

The nonlinear Fokker-Planck equation (\ref{eq:nonlinear Fokker-Planck ns}) has stationary states
that we describe in the following proposition.

\begin{proposition}
	\label{prop:equilibrium ns}
	Let $\Xi$ be the set of solutions of the compatibility condition equation
	\begin{equation}
		\label{eq:compatibility condition ns}
		\xi = G(\xi).
	\end{equation}
	For any $\xi \in \Xi$, the state
	\begin{equation}
		\label{eq:equilibrium ns}
		\rho_\xi(x,u) = \frac{1}{L} F_\xi(u),  \quad F_\xi(u)= \frac{1}{\sqrt{\pi\sigma^2}} \exp\left(- \frac{(u-\xi)^2}{\sigma^2} \right)
	\end{equation}
	is a stationary solution of (\ref{eq:nonlinear Fokker-Planck ns}).
\end{proposition}
Note, however, that we could not prove that any stationary state is of the form (\ref{eq:equilibrium ns}).
It is possible to prove that any stationary state that is uniform in space is of the form (\ref{eq:equilibrium ns}),
but we could not prove that a stationary state must be uniform in space.

The linear stability analysis of the stationary states listed in Proposition \ref{prop:equilibrium ns}
 follows the same lines as in the symmetric case (\ref{def:aveu}).
Let $\xi \in \Xi$ and consider 
\begin{equation}
	\rho(t,x,u) = \rho_\xi(x,u)+\rho^{(1)}(t,x,u) = \frac{1}{L}F_\xi(u) + \rho^{(1)}(t,x,u),
\end{equation}
for small perturbation $\rho^{(1)}$.  We linearize the nonlinear Fokker-Planck equation (\ref{eq:nonlinear Fokker-Planck ns})
\begin{align}
	\label{eq:linearized Fokker-Planck ns}
	\frac{\partial \rho^{(1)}}{\partial t} &= 
	-u\frac{\partial \rho^{(1)}}{\partial x} 
	-\frac{\partial}{\partial u}\left[(\xi-u) \rho^{(1)}\right]\\
\notag	&\quad + \frac{1}{L }G'(\xi)\left[\iint (\xi-u') \phi(\|x'\|) \rho^{(1)}(t, x-x', u') du'dx' \right] F'_\xi(u)
	+ \frac{1}{2} \sigma^2 \frac{\partial^2 \rho^{(1)}}{\partial u^2}.
\end{align}
The perturbation  $\rho^{(1)}$ is periodic in $x$ and can be expanded as (\ref{eq:expand:rho1k}).
The Fourier coefficients $\hat{\rho}^{(1)}_k$ satisfy uncoupled equations
\begin{align}
	\notag
	\frac{\partial\rho_k^{(1)}}{\partial t} &=
	\frac{i 2\pi k}{L} u \rho_k^{(1)}  
	-\frac{\partial}{\partial u}\left[(\xi-u) \rho_k^{(1)}\right]\\
	&\quad + G'(\xi) \phi_k  \left[\int (\xi- u') \rho_k^{(1)}(t,u') du'\right] F'_\xi(u) + \frac{1}{2} \sigma^2 \frac{\partial^2 \rho_k^{(1)}}{\partial u^2}.
	\label{eq:modelink ns}
\end{align}
The necessary and sufficient condition for the linear stability of the $0$th order mode is
\begin{equation}
G'(\xi) < 1 .
\end{equation}
The sufficient condition for the linear stability of the other modes is that $|G'(\xi) \phi_k|<1$ for all $k \neq 0$ 
and the noise level $\sigma$ 
should be larger than a threshold value $\sigma_c$, which is, however different from the threshold value
of the symmetric case (\ref{eq:model}-\ref{def:aveu}).

\section{Conclusion}
\label{sec:conclusion}

We have analyzed the Czir\'ok model for the collective motion of locusts.  The mean-field theory is used to obtain a nonlinear 
Fokker-Planck equation as the number of agents tends to infinity.  We analyze the phase transition between the disorder and order states by the existence condition for the order states.  
We then perform the linear stability analysis on the stationary states, and prove that the order states are stable 
for a sufficiently  large noise level and when the interaction between the velocities of the particles is neither too small nor too strong. 
 We provide the fluctuation analysis and the large 
deviations principle.  
For a large but finite system we calculate the asymptotic, exponentially small transition probability from one order state to the other.
 Our analytical findings are verified by the extensive numerically simulations and are found in agreement
 with the experimental observations. 

\appendix

\section{Proof of Lemma \ref{prop:stablity of w_k}}
\label{app:A}
We have
\begin{align*}
	g_k(t,0) &= \frac{\sigma^2 }{4} D_k^2 (1-e^{-2t}) + (i\xi -\sigma^2 D_k) D_k(1-e^{-t}) + ( -i\xi D_k + \frac{\sigma^2}{2}D_k^2)t\\
	&\simeq ( -i\xi D_k + \frac{\sigma^2}{2}D_k^2)t,\quad \text{as $t\to\infty$,}
\end{align*}
\begin{equation*}
	\partial_\eta g_k(t,0) = \frac{\sigma^2 }{2}D_k(1-e^{-2t}) + (i\xi - \sigma^2 D_k)(1-e^{-t}) 
	\simeq (i\xi - \frac{\sigma^2}{2}D_k),\quad \text{as $t\to\infty$.}
\end{equation*}
Since $\hat{\rho}_k^{(1)} (0,\cdot) $ and $\partial_\eta \hat{\rho}_k^{(1)} (0,\cdot)$ are bounded we find that $\psi_k(t)$ decays exponential 
in time as $\exp ( -\sigma^2 D_k^2 t/2)$, which gives the first item of the Lemma.

We can show that $R_k(t) \in L^1[0,\infty)$ because it is the product of a bounded function with $e^{-g_k(t,0)}$ which decays exponentially as $t \to \infty$.
Note that it is important that $k\neq 0$ and $\sigma>0$ to ensure the decay.
Recall that $D_k = {2\pi k}/{L}$ and $\beta_k=D_k(1-e^{-t})$.  Then we can write $g_k(t,0)$, $\partial_\eta g_k(t,0)$, $H_\xi(-\beta_k)$, 
and $H_\xi'(-\beta_k)$ in terms of $D_k$ and $\beta_k$:
\begin{align*}
	g_k(t,0) &= -\frac{1}{4}\sigma^2\beta_k^2 - \frac{1}{2}\sigma^2 D_k\beta_k + i\xi\beta_k - i\xi D_kt + \frac{1}{2}\sigma^2 D_k^2 t,\\
	\partial_\eta g_k(t,0) &= \frac{1}{2}\sigma^2(e^{-t}-1)\beta_k + i\xi(1-e^{-t}),\\
	H_\xi(-\beta_k) &= -\beta_k \hat{F}_\xi(-\beta_k) = -\beta_k\exp\Big(i\xi\beta_k-\frac{1}{4}\sigma^2\beta_k^2\Big), \\
	H_\xi'(-\beta_k) &= \Big( 1+i\xi\beta_k-\frac{1}{2}\sigma^2\beta_k^2\Big)\exp\Big(i\xi\beta_k-\frac{1}{4}\sigma^2\beta_k^2\Big).
\end{align*}
Therefore,
\begin{equation*}
	-\partial_{\eta}g_{k}(t,0)H_{\xi}\left(-\beta_k\right)+e^{-t}H_{\xi}'(-\beta_k) 
	= \Big(-\frac{1}{2}\sigma^2\beta_k^2 + i\xi\beta_k+e^{-t}\Big)
	\exp\Big(i\xi\beta_k-\frac{1}{4}\sigma^2\beta_k^2\Big),
\end{equation*}
and we can also write $R_k(t)$ in terms of $D_k$ and $\beta_k$:
\begin{align*}
	R_k(t) &= G'(\xi)\phi_k \left[-\partial_\eta g_k(t,0)H_\xi(-\beta_k) + e^{-t}H_\xi'(-\beta_k)\right]e^{-g_k(t,0)}\\
	&= G'(\xi)\phi_k \Big(-\frac{1}{2}\sigma^2\beta_k^2 + i\xi\beta_k+e^{-t}\Big)
	\exp\Big[\frac{\sigma^2}{2}D_k\beta_k + i\xi D_kt - \frac{\sigma^2}{2}D_k^2 t\Big].
\end{align*}
The $L^1$-norm of $R_k $ can be bounded by
\begin{equation}
	\label{eq:L1 norm of Rk}
	\int_0^\infty |R_k(t)| dt \leq |G'(\xi)\phi_k|\int_0^\infty 
	\Big|-\frac{1}{2}\sigma^2\beta_k^2 + i\xi\beta_k+e^{-t}\Big|\exp\Big[\frac{\sigma^2}{2}D_k\beta_k - \frac{\sigma^2}{2}D_k^2 t\Big] dt.
\end{equation}
Because $\beta_k=D_k(1-e^{-t})$, we use the following bounds for (\ref{eq:L1 norm of Rk}):
\begin{align*}
	\Big|-\frac{1}{2}\sigma^2\beta_k^2 + i\xi\beta_k+e^{-t}\Big| &\leq \frac{\sigma^2}{2}D_k^2(1-e^{-t})^2 + |\xi||D_k|(1-e^{-t}) + e^{-t}\\
	&\leq 
	\begin{cases}
		\frac{\sigma^2}{2}D_k^2 t^2 + |\xi||D_k| t + 1, & 0\leq t\leq 1,\\
		\frac{\sigma^2}{2}D_k^2 + |\xi||D_k| + e^{-t}, & 1<t,
	\end{cases}
\end{align*}
\begin{equation*}
	\exp\Big[\frac{\sigma^2}{2}D_k\beta_k - \frac{\sigma^2}{2}D_k^2 t\Big] =
	\exp\Big[\frac{\sigma^2}{2}D_k^2(1-e^{-t}-t)\Big]
	\leq 
	\begin{cases}
		e^{\frac{\sigma^2}{2}D_k^2(-\frac{1}{4}t^2)}, & 0\leq t\leq 1,\\
		e^{\frac{\sigma^2}{2}D_k^2(1-t)}, & 1<t.
	\end{cases}
\end{equation*}
We have the following estimates:
\begin{equation*}
	\int_0^1 \Big(\frac{\sigma^2}{2}D_k^2 t^2 + |\xi||D_k| t + 1\Big) 
	e^{\frac{\sigma^2}{2}D_k^2(-\frac{1}{4}t^2)} dt
	\leq \frac{\sqrt{2\pi}}{\frac{\sigma}{2}|D_k|}\Big(1 + \sqrt{\frac{2}{\pi}}\frac{|\xi|}{\sigma} + \frac{1}{2}\Big),
\end{equation*}
\begin{equation*}
	\int_1^\infty \Big(\frac{\sigma^2}{2}D_k^2 + |\xi||D_k| + e^{-t}\Big) 
	e^{\frac{\sigma^2}{2}D_k^2(1-t)}dt
	= 1 + \frac{|\xi|}{\frac{\sigma^2}{2} |D_k|} + \frac{e^{-1}}{1+\frac{\sigma^2}{2}D_k^2}.
\end{equation*}
Therefore a sufficient condition to ensure that $\int_0^\infty |R_k(t)|dt<1$ is
	\begin{equation}
		\label{eq:stablity of mode k}
		|G'(\xi)\phi_k|
		< \left(1 + \frac{3\sqrt{2\pi}}{\sigma |D_k|} + \frac{3|\xi|}{\frac{\sigma^2}{2}|D_k|} + \frac{e^{-1}}{1+\frac{\sigma^2}{2}D_k^2} \right)^{-1}.
	\end{equation}
By the fact that $|D_k| = 2\pi |k|/L$ is increasing in $|k|$, 
we can have a simplified sufficient condition from (\ref{eq:stablity of mode k}):
	\begin{equation}
		\label{eq:simplified stablity of mode k}
		|G'(\xi) \phi_k |
		< \left(1 + \frac{3\sqrt{2\pi}}{\sigma D_1} + \frac{3|\xi|}{\frac{\sigma^2}{2}D_1} + \frac{e^{-1}}{1+\frac{\sigma^2}{2}D_1^2} \right)^{-1}.
	\end{equation}
If $| G'(\xi) \phi_k| <1$ for all nonzero $k$, then this condition is satisfied  for all nonzero $k$ if $\sigma$ is large enough.
This completes the proof of the second item of the Lemma.


\end{document}